\documentclass[trsc,nonblindrev]{informs3} 

\OneAndAHalfSpacedXI 

\usepackage{natbib}
 \bibpunct[, ]{(}{)}{,}{a}{}{,}%
 \def\bibfont{\small}%
\TheoremsNumberedThrough     

\EquationsNumberedThrough    

\usepackage{mathtools}
\usepackage{multirow}
\usepackage{tikz}
\usepackage{subcaption}
\usepackage{comment}
\usepackage{accents} 
\usepackage{algorithmic} 

\usepackage[labelformat=simple]{subcaption}

\usepackage{bm}
\newcommand{\vc}[1]{\bm{#1}}

\def\Re{{\mathbb R}}

\def\tr{^{\intercal}}
\def\pol{^{\circ}}

\def\Re{{\mathbb R}}

\def\proj{\mathop{\rm proj}}

\def\Re{{\mathbb R}}

\def\S{\mathcal{S}}
\def\LP{\mathrm{LP}}

\def\BB{B\&B }

\newcommand{\black}[1]{\textcolor{black}{#1}}

\newcommand{\suchthat}{\;\ifnum\currentgrouptype=16 \middle\fi|\;}

\usepackage[ruled,vlined,linesnumbered]{algorithm2e}

\begin{document}


\RUNAUTHOR{Rajabalizadeh and Davarnia}

\RUNTITLE{Solving a Class of CGLPs via ML}

\TITLE{Solving a Class of Cut-Generating Linear Programs via Machine Learning}

\ARTICLEAUTHORS{%
\AUTHOR{Atefeh Rajabalizadeh, Danial Davarnia}
\AFF{Department of Industrial and Manufacturing Systems Engineering, Iowa State University, Ames, IA 50011, \EMAIL{alizadeh@iastate.edu}, \EMAIL{davarnia@iastate.edu}}
} 

\ABSTRACT{%
Cut-generating linear programs (CGLPs) play a key role as a separation oracle to produce valid inequalities for the feasible region of mixed-integer programs.
When incorporated inside branch-and-bound, the cutting planes obtained from CGLPs help to tighten relaxations and improve dual bounds.
However, running the CGLPs at the nodes of the branch-and-bound tree is computationally cumbersome due to the large number of node candidates and the lack of a priori knowledge on which nodes admit useful cutting planes.
As a result, CGLPs are often avoided at default settings of branch-and-cut algorithms despite their potential impact on improving dual bounds. 
In this paper, we propose a novel framework based on machine learning to approximate the optimal value of a CGLP class that determines whether a cutting plane can be generated at a node of the branch-and-bound tree.
Translating the CGLP as an indicator function of the objective function vector, we show that it can be approximated through conventional data classification techniques. 
We provide a systematic procedure to efficiently generate training data sets for the corresponding classification problem based on the CGLP structure.
We conduct computational experiments on benchmark instances using classification methods such as logistic regression.
These results suggest that the approximate CGLP obtained from classification can improve the solution time compared to that of conventional cutting plane methods.
Our proposed framework can be efficiently applied to a large number of nodes in the branch-and-bound tree to identify the best candidates for adding a cut.		             
}%


	\KEYWORDS{cutting planes; cut-generating linear programs; machine learning; data classification; function approximation}
\HISTORY{}

\maketitle

%


\section{Introduction}
\label{sec:introduction}

Solving optimization problems with a large number of constraints can be computationally challenging due to explicit representation of all constraints in the base model.
In this situation---if the problem structure allows---a cut-generating linear program (CGLP) is employed as a separation oracle to add cutting planes successively to a relaxation of the model.
This oracle is particularly useful in computing the projection of a set onto a lower-dimensional space, which is a common technique for tightening relaxations of mixed integer programs (MIPs). 
The core structure of this projection procedure is given in the following proposition; see \cite{conforti:co:za:2014} for a detailed account on the role of CGLPs in MIPs.

\begin{proposition} \label{prop:projection}
	Consider a set $\mathcal{S} = \{(\vc{\alpha},\vc{\beta}) \in \Re^{p+q} \, | \, A\vc{\alpha} + B\vc{\beta} \leq \vc{\gamma}\}$, where $A$, $B$, and $\vc{\gamma}$ are parameter matrices of matching dimensions.
	Then, the collection of inequalities $\bar{\vc{x}} B\vc{\beta} \leq \bar{\vc{x}} \vc{\gamma}$ for all feasible solutions $\bar{\vc{x}} \in \mathcal{C} := \{\vc{x} \in \Re_+^{n} \, | \, \vc{x} A = \vc{0}\}$ describes the projection of $\mathcal{S}$ onto the $\vc{\beta}$-space, denoted by $\proj_{\beta} (\mathcal{S})$.
	Further, for any point $\bar{\vc{\beta}} \in \Re^q$, if $\max\{\vc{x}(B\bar{\vc{\beta}} - \vc{\gamma}) \, | \, \vc{x} \in \mathcal{C}\} = 0$, then $\bar{\vc{\beta}} \in \proj_{\beta} (\mathcal{S})$.
	Otherwise, if $\max\{\vc{x}(B\bar{\vc{\beta}} - \vc{\gamma}) \, | \, \vc{x} \in \mathcal{C}\} > 0$ with \black{an optimal ray} $\bar{\vc{x}}$, then $\bar{\vc{\beta}}$ can be separated from $\proj_{\beta} (\mathcal{S})$ by $\bar{\vc{x}} B\vc{\beta} \leq \bar{\vc{x}} \vc{\gamma}$.
	\Halmos
\end{proposition}

In view of Proposition~\ref{prop:projection}, the problem 
$\max\{\vc{x}(B\bar{\vc{\beta}} - \vc{\gamma}) \, | \, \vc{x} \in \mathcal{C}\}$ is referred to as the CGLP associated with the projection of $\mathcal{S}$.
This CGLP is a critical tool for solving both special-structured and general MIPs in lower-dimensional spaces.
For instance, CGLPs form the basis of disjunctive programming \cite{Bal79} as a predominant method to obtain strong cutting planes; see \cite{NemWol99}.
In this approach, a convex hull description of the MIP is constructed in a higher-dimensional space through convexifying the disjunctive union of a finite number of polyhedra, which is then projected onto the space of original variables through solving a CGLP.
The resulting inequalities are added to the relaxations of the problem as separating cutting planes.

In the CGLP of Proposition~\ref{prop:projection}, different values of the point $\bar{\vc{\beta}}$---which is desired to be separated---lead to different objective functions, while the feasible region defined by the polyhedral cone $\mathcal{C}$ remains the same.
Depending on the objective function, the optimal value of the CGLP can be either $0$ or $\infty$, indicating whether or not the point can be separated from the projection. 
In MIP applications, this oracle is often invoked at the nodes of the branch-and-bound (B\&B) tree to obtain cutting planes that separate the optimal solution of the LP relaxation to improve the dual bound.
Even though the CGLP is a linear program, solving it repeatedly at the nodes of the \BB tree can be computationally prohibitive even for moderate-size problems.
Various studies in the literature are devoted to finding an approximate CGLP by reducing variables and constraints with the goal of reducing the solution time at the price of producing potentially weaker cuts; see \cite{BALAS:per:2002}.
Despite the time improvement achieved by such approximate models, their repetitive invocation at the massive scales required in branch-and-cut methods could still render them practically expensive.

A recent study by \cite{davarnia:ra:ho:2022} shows that a new class of cutting planes that targets \textit{consistency} of the partial assignments corresponding to the nodes of the \BB tree can be more effective than the traditional cutting planes in reducing memory allocation and solution time.
It is shown in this work that consistency cuts can be generated based only on the optimal value of the CGLP, unlike traditional cutting planes that require the optimal solution of the CGLP in addition to its optimal value.
To generate consistency cuts, we only need the characterization of the \textit{indicator function} of the CGLP that has two possible outcomes: $0$ when the optimal value of the CGLP is $0$; and $1$ when the optimal value of the CGLP is $\infty$; see \cite{Rockafellar:97} for properties of indicator functions of convex programs.

A common occurrence when implementing CGLPs inside of a branch-and-cut scheme is that the CGLP is solved numerous times over a \textit{fixed} cone at a given layer of the \BB tree.
This observation raises the question: whether one can find an approximation of the CGLP indicator function \textit{over} the entire domain of the objective function vectors, so that for any given vector, we can quickly evaluate the approximate indicator function without the need to solve an LP or an approximation model?
Even though such an approximate function will not guarantee the precise output value of a given input, it can identify a priori best node candidates in the \BB tree that are \textit{most likely} to admit cutting planes.
When used as a preprocessing technique in the B\&B, this approach can lead to a substantial time save for large-scale repetitive invocation of CGLPs.
In this paper, we propose a new framework to approximate the CGLP indicator function through the lens of function approximation and machine learning.

Function approximation has a rich history in mathematics and computer science; see \cite{rivlin:1969} for an introduction.
One of its most common applications is concerned with situations where the explicit form of the underlying function is difficult to obtain, and evaluating the function at a point is expensive.
As a result, a subset of \textit{critical} points in the domain of the function is selected, for which the function value is calculated.
Then, an approximate function is computed by minimizing the residual error between the estimate value and the actual value of the selected points in the domain.
This approach is sometimes referred to as \textit{interpolation} in the literature.
Viewing the selected points as \textit{input} and the actual function value as \textit{response}, the interpolation technique shares core principles with machine learning methods in finding the best fit for the underlying function when considering the least-squared error as the objective \citep{busoniu:ba:de:er:2017}.
For instance, the neural network at the core of deep learning is a classical example of function approximation \citep{wasserman:2006}. 
Non-parametric machine learning models such as $k$-nearest neighbor and support-vector machines are also directly used for function approximation \citep{hammer:ge:2003}.
Even parametric models such as regression can be used for interpolation when the objective is to find estimators that minimize the least-squared error \citep{friedman:94}.

Using the above analogy, we employ standard machine learning and, specifically, classification methods to approximate the indicator function of the CGLP.
The use of machine learning in the branch-and-cut process has rapidly grown over the past few years due to their potential in assisting with \BB strategies such as node selection and branching order \citep{khalil:le:so:ne:di:2016,alvarez:lo:we:2017,zarpellon:jo:lo:be:2020}.
In cutting plane domains, machine learning methods have been used to select valid inequalities from a pool of candidate cuts \cite{tang:ag:fa:2020,balcan:pr:sa:vi:2021}.
The main limitation of the above approaches is that they work best when the target value of the output is \textit{categorical}, i.e., they select candidates (node, branch, or cut) among a given pool of finitely many alternatives.
As a result, such techniques are not viable for producing new cutting planes that have not been previously identified, and hence machine learning methods have never been used in combination with cut-generating efforts inside B\&B.
We address this gap in this paper by providing a novel perspective to model the CGLP through function approximation and use classification methods to produce new cutting planes.

The contributions of this paper are as follows.
We design a novel framework based on machine learning to approximate CGLPs.
To our knowledge, this is the first work that takes advantage of data classification to help with cut-generating efforts in the branch-and-cut process.
This framework is applicable to a broad range of optimization problems beyond CGLP such as the convex hull membership problem as a fundamental problem in computer science and mathematics.
Further, we develop a theoretical and systematic procedure to efficiently generate different types of training data sets for our classification problem.
We show how these methods can be effectively implemented inside branch-and-cut. 
Preliminary computational experiments suggest that the resulting approximation improves the solution time of the cut-generating efforts compared to that of the conventional cutting plane methods that rely on solving the CGLP to produce cuts.

The remainder of the paper is organized as follows.
In Section~\ref{sec:classification}, we introduce structural properties of the CGLP that allow for transforming it into a classification problem.
In addition, we develop a theoretical framework to identify class-0 and class-1 data points that will be used as training sets for classification.
We discuss how the developed framework can be applied to a CGLP problem during \BB in Section~\ref{sec:consistency}.
Section~\ref{sec:computation} evaluates the performance of the proposed approach through computational experiments.
Concluding remarks are given in Section~\ref{sec:conclusion}.
Additional results and discussions are provided in Appendices that can be found in the online supplement.

\noindent
\textbf{Notation.} We use bold letters to denote vectors. To simplify notation, we do not use the transposition symbol when representing the scalar products between vectors. For instance, the product $\vc{x}\vc{y}$ implies that $\vc{x}$ is a row vector and $\vc{y}$ is a column vector of matching dimensions.

\section{Formulation and Classification} \label{sec:classification}

Consider a general polyhedral cone $\mathcal{C} = \{\vc{x} \in \Re^n \, | \, \vc{a}^i \vc{x} \leq 0, \, \forall i \in I \}$, where $\vc{a}^i$ is a row vector of appropriate dimension, and $I = \{1, \dotsc, m, m+1, \dotsc, m+n\}$ represents the index set of constraints, including the non-negativity bounds on variables as the last $n$ constraints, i.e., $\vc{a}^{m+j} = -\vc{e}^j$ for $j \in N = \{1, \dotsc, n\}$ where $\vc{e}^j$ is the $j$th unit vector in $\Re^n$.
It is easy to verify that the projection cone of the CGLP of Proposition~\ref{prop:projection} can be formulated in the form of $\mathcal{C}$.

For any given $\vc{c} \in \Re^n$, we are interested in solving $z^*(\vc{c}) = \max \{\vc{c} \vc{x} \, | \, \vc{x} \in \mathcal{C} \}$, which is referred to as the \textit{support function} of $\mathcal{C}$ in the literature \citep{Rockafellar:97}.
This support function can have two outcomes $z^*(\vc{c}) \in \{0, \infty\}$, which can be translated into an \textit{indicator function} with values 0 if $z^*(\vc{c}) = 0$, and 1 if $z^*(\vc{c}) = \infty$.
The goal is to use data classification techniques to approximate the value of this indicator function without solving the corresponding optimization problem.
In particular, we aim to classify vectors $\vc{c} \in \Re^n$ into class-0 if $z^*(\vc{c}) = 0$ and class-1 if $z^*(\vc{c}) = \infty$.
To this end, we propose to train a machine learning model based on some given pairs of input and response variables of the form $(\vc{c}^k, z^*(\vc{c}^k))$ for $k \in K$, where $K$ is a data set.

To build a training set that is representative of the actual geometry of the input vector space for $z^*(\vc{c})$, we require a balanced data set composed of both classes.
In the sequel, we first identify the structure of vectors $\vc{c}$ that belong to each class, and then choose an appropriate subset of such vectors in our training set with an emphasis on the efficiency of generation.
Define the polar of $\mathcal{C}$ as $\mathcal{C}\pol = \{\vc{y} \in \Re^n \, | \, \vc{y} \vc{x} \leq 0, \forall \vc{x} \in \mathcal{C}\}$.
It is clear from this definition that $\mathcal{C}\pol$ can be described as the cone generated by vectors $\vc{a}^i$ for $i \in I$.

\begin{proposition} \label{prop:polar}
	A vector $\vc{c} \in \Re^n$ belongs to class-0 if and only if $\vc{c} \in \mathcal{C}\pol$.
\end{proposition}

\proof{Proof.}
For the direct implication, assume that $\vc{c}$ belongs to class-0, i.e., $z^*(\vc{c}) = 0$.
By definition of $z^*(\vc{c})$, we have that $\vc{c} \vc{x} \leq 0$ for all $\vc{x} \in \mathcal{C}$.
For the reverse implication, assume that $\vc{c} \in \mathcal{C}\pol$, i.e., $z^*(\vc{c}) = \max \{\vc{c} \vc{x} \, | \, \vc{x} \in \mathcal{C} \} \leq 0$.
Since the origin is a feasible solution to this optimization problem, we conclude that $z^*(\vc{c}) = 0$. 
\Halmos \endproof

\begin{remark} \label{remark:CHM}
	In addition to representing the CGLP, the indicator function $z^*(\vc{c})$ can also be used to solve combinatorial problems such as the \textit{convex hull membership} (CHM) problem, which is the problem of deciding whether a given point belongs to the convex hull of a finite number of points.
	The CHM has applications in computer science \cite{karmarkar:84}, computational geometry \cite{toth:or:go:2017}, dynamic programming \cite{bertsikas:2017}, and decision diagrams \cite{davarnia:va:2020,davarnia:2021,salemi:da:2022}.
	To illustrate this relation, consider the set of points $\vc{x}^i \in \Re^n$ for $i=1,\dotsc,k$, and a point $\bar{\vc{x}} \in \Re^n$.
	It is easy to verify that, $\bar{\vc{x}}$ belongs to the convex hull of $\vc{x}^i$ for $i=1,\dotsc,k$ if and only if the extended vector $\vc{c} = (1, \bar{\vc{x}})$ belongs to the cone generated by the extended vectors $\vc{a}^i = (1, \vc{x}^i)$ $i=1,\dotsc,k$; see Section 2 of \cite{Rockafellar:97} for a detailed account.
	Such a cone can be viewed as a polar cone $\mathcal{C}\pol$ associated with the indicator function $z^*(\vc{c})$.
	As a result, the classification techniques we develop in this paper can also be applied to the CHM problem.
	We present computational experiments for such application in Section~\ref{subsec:chm}.
\end{remark}

\begin{figure}[!hbt]
	\centering
	\includegraphics[scale=0.25]{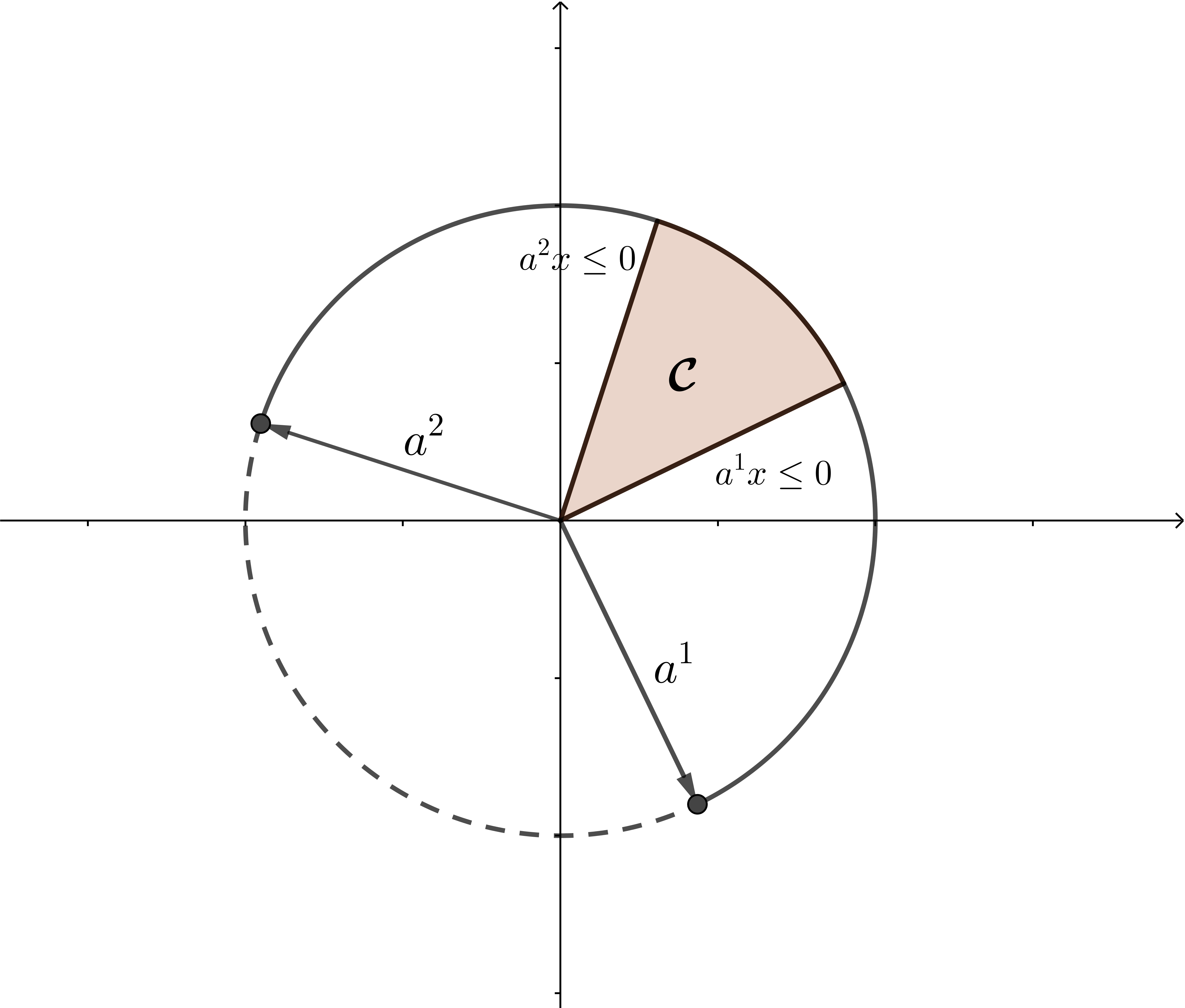}
	\caption{The CGLP projection cone and its polar cone.}
	\label{fig:P1}
\end{figure}

\begin{figure}[!hbt]
	\centering
	\includegraphics[scale=0.25]{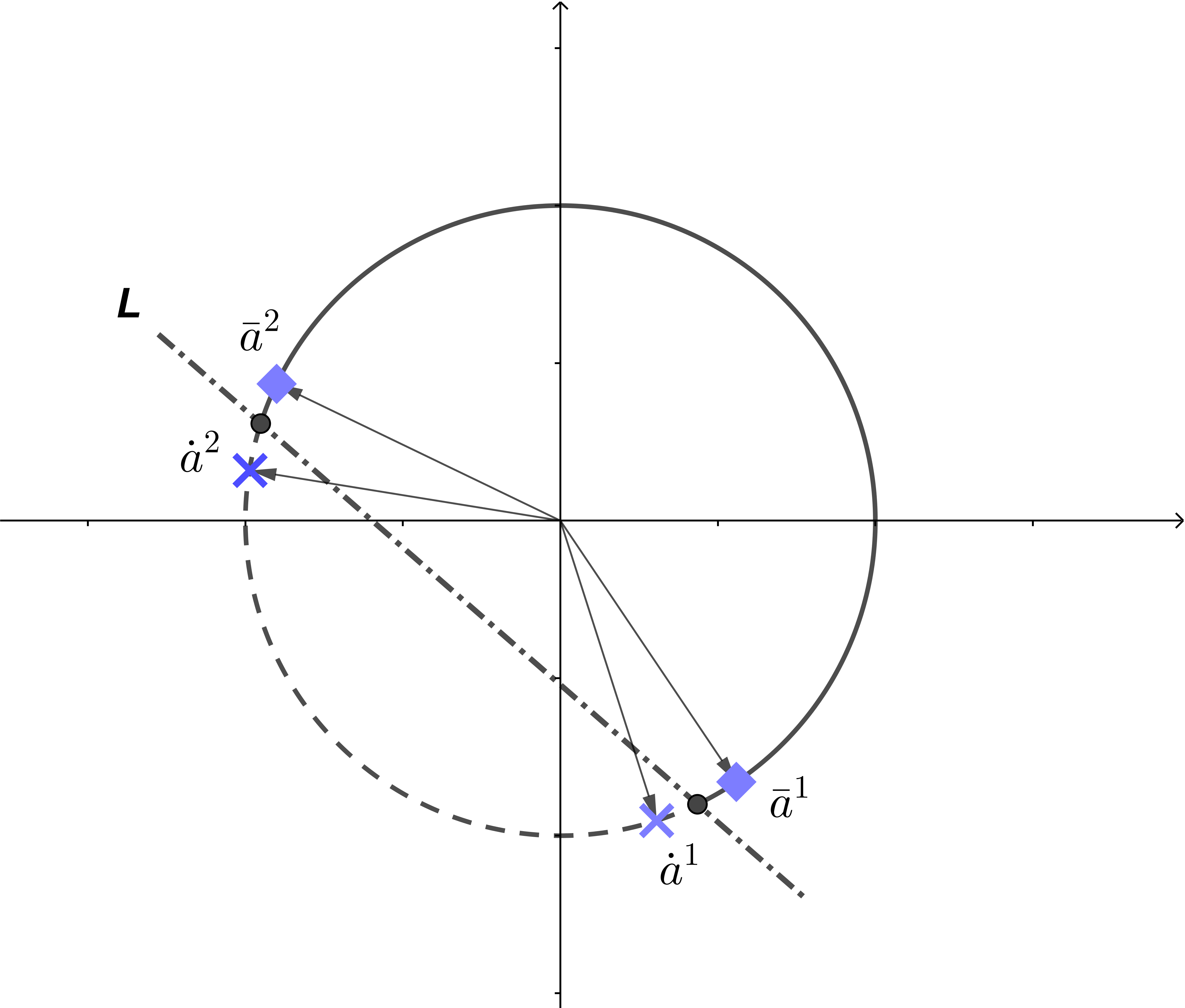}
	\caption{Classification of class-0 and class-1 with respect to the polar of the projection cone.}
	\label{fig:P2}
\end{figure}

The goal of classification is to find a classifier that separates class-0 and class-1 vectors.
It follows from Proposition~\ref{prop:polar} that these two classes are separated by the boundary of the polar cone $\mathcal{C}\pol$. 
As a result, we seek to train a classifier that approximates the boundary of the polar cone. 
As common in function approximation, selecting \textit{critical} points from the domain of the unknown function is key to obtain good approximations for that function.
Critical points are representative of main functional patterns over the domain of the function.
While characterizing such subset of points is difficult for general functions, they are often selected from the vicinity of the breakpoints that form the boundary of the graph of indicator functions, which maps to the boundary of $\mathcal{C}\pol$ in our problem.
Such a selection that contains points on both sides of the boundary improves the approximation accuracy as illustrated in Figures~\ref{fig:P1} and \ref{fig:P2}.
In particular, Figure~\ref{fig:P1} shows the restriction of the cone $\mathcal{C} = \{\vc{x} \in \Re^2: a^i \vc{x} \leq 0, i=1, 2\}$ to the unit disc.
The corresponding polar cone $\mathcal{C}^{\pol}$ is generated by vectors $a^1$ and $a^2$.
Assuming that all vectors $\vc{c}$ are normalized, class-0 vectors correspond to the dashed arc between $a^1$ and $a^2$, while the complementary solid arc on the disc represents class-1 vectors. 
Vectors $\dot{a}^i$ and $\bar{a}^i$ for $i=1,2$ in Figure~\ref{fig:P2} respectively represent class-0 and class-1 training data close to the boundary vectors $a^i$.
Using this data, the goal is to obtain a classifier such as $L$ that separates the dashed and solid arc, and thereby yielding the desired approximation. 

In the remainder of this section, we identify vectors of class-0 and class-1 \textit{close} to the boundary of the polar cone.
Identifying class-0 vectors in the vicinity of the boundary is achieved by the conic combination of vectors generating the polar cone, as demonstrated next.
In our derivation, we use the fact that $\mathcal{C}\pol$ is the cone generated by the constraint vectors of $\mathcal{C}$, i.e., $\{\vc{a}^i\}_{i \in I}$; see \cite{Rockafellar:97}.
Further, for any weight vector $\vc{w} \in \Re^{m+n}$, we denote the corresponding normalized linear combination of these constraint vectors by $\vc{d}(\vc{w}) = \frac{\sum_{i \in I} w_i \vc{a}^i}{||\sum_{i \in I} w_i \vc{a}^i||}$ when $\vc{w} \neq \vc{0}$, and $\vc{d}(\vc{w}) = \vc{0}$ otherwise.

\begin{proposition} \label{prop:class-0 vicinity}
	Vector $\vc{c} \in \Re^n$ belongs to class-0 if and only if $\vc{c} = \vc{u} + \epsilon \vc{d}(\vc{w})$ for some $\vc{u}$ on the boundary of $\mathcal{C}\pol$, $\epsilon > 0$ and $w_i \geq 0$ for all $i \in I$. 
\end{proposition}

\proof{Proof.}
For the direct implication, assume that $\vc{c} = \vc{u} + \epsilon \vc{d}(\vc{w})$ for some $\vc{u}$ on the boundary of $\mathcal{C}\pol$, $\epsilon > 0$ and $w_i \geq 0$ for all $i \in I$.
Then, $\vc{c}$ can be viewed as a conic combination of the constraint vectors $\{\vc{a}^i\}_{i \in I}$, implying that $\vc{c} \in \mathcal{C}\pol$.
For the reverse implication, assume that $\vc{c} \in \mathcal{C}\pol$.
Therefore, we can write that $\vc{c} = \sum_{i \in I} \bar{w}_i \vc{a}^i$ with $\black{\bar{w}_i} \geq 0$.
The result holds by setting $\epsilon = ||\sum_{i \in I} \bar{w}_i \vc{a}^i||$, $\vc{w} = \bar{\vc{w}}$, and $\vc{u} = \vc{0}$ as the origin is a point on the boundary of $\mathcal{C}\pol$ by definition.
\Halmos 
\endproof

In view of Proposition~\ref{prop:class-0 vicinity}, we note that, if the resulting vector $\vc{c}$ is not on the boundary of the polar cone, its \textit{closeness} to the boundary of the polar cone is controlled by $\epsilon$. In other words, it can be viewed as the \textit{step length} in the direction of $\vc{d}(\vc{w})$ as we move away from vector $\vc{u}$ on the boundary of $\mathcal{C}\pol$.
This statement implies that, given a vector $\vc{u}$ on the boundary of the polar cone and a normalized \textit{direction vector} $\vc{d}(\vc{w})$, smaller values of $\epsilon$ lead to vectors $\vc{c}$ that are closer to the boundary of the polar cone. 
More specifically, for a given distance $r > 0$, we can pick any $\epsilon \leq r$ to ensure that the distance of $\vc{c}$ from $\vc{u}$ calculated as $||\vc{c} - \vc{u}|| = ||\vc{u} + \epsilon \vc{d}(\vc{w}) - \vc{u}|| = \epsilon$ is no larger than $r$.

The following corollary presents a special case for using the results of Proposition~\ref{prop:class-0 vicinity} to generate class-0 vectors arbitrarily close to the boundary of the polar cone.
We will use this method to efficiently generate class-0 vectors in the training set for our computational experiments given in Section~\ref{sec:computation}.

\begin{corollary} \label{cor:class-0}
	Let $i \in I$ be the index of a constraint of $\mathcal{C}$.
	Then, for any $\bar{\epsilon} > 0$ and $v_j \black{\leq} 0$ for all $j \in N$, the vector $\vc{c} = \vc{a}^i + \bar{\epsilon} \bar{\vc{v}}$, where $\bar{\vc{v}} = \frac{\vc{v}}{||\vc{v}||}$ if $\vc{v} \neq \vc{0}$, and $\bar{\vc{v}} = \vc{0}$ otherwise, belongs to class-0.
\end{corollary}

\proof{Proof.}
We consider two cases.
First, assume that the constraint with index $i$ is non-redundant in $\mathcal{C}$.
Then, $\vc{a}^i$ is on the boundary of $\mathcal{C}\pol$.
The result follows from Proposition~\ref{prop:class-0 vicinity} by setting $\vc{u} = \vc{a}^i$, $\epsilon = \bar{\epsilon}$, $w_j = 0$ for $j = 1, \dotsc, m$, and $w_{m+k} = -v_k$ for all $k \in N$.
Second, assume that the constraint with index $i$ is redundant in $\mathcal{C}$.
Proposition~\ref{prop:class-0 vicinity} implies that $\vc{a}^i = \vc{u} + \tilde{\epsilon} \vc{d}(\tilde{\vc{w}})$ for some $\vc{u}$ on the boundary of $\mathcal{C}\pol$, $\tilde{\epsilon} > 0$ and $\tilde{w}_j \geq 0$ for all $j \in I$ with $\tilde{\vc{w}} \neq \vc{0}$.
There are two cases. If $\vc{v} = \vc{0}$, then $\vc{c} = \vc{a}^i = \vc{u} + \tilde{\epsilon} \vc{d}(\tilde{\vc{w}})$, proving the result because of Proposition~\ref{prop:class-0 vicinity}.
If $\vc{v} \neq \vc{0}$, then we write that $\vc{c} = \vc{a}^i + \bar{\epsilon} \bar{\vc{v}} = \vc{u} + \epsilon \, \vc{d}(\bar{\vc{w}})$ where 
$\bar{w}_j = \frac{\tilde{\epsilon} \tilde{w}_j}{||\sum_{i \in I}\tilde{w}_i \vc{a}^i||}$ for $j = 1, \dotsc, m$, $\bar{w}_{j} = \frac{\tilde{\epsilon} \tilde{w}_j}{||\sum_{i \in I}\tilde{w}_i \vc{a}^i||} - \frac{\bar{\epsilon}v_{j - m}}{||\vc{v}||}$ for all $j = m+1, \dotsc, m+n$, and $\epsilon = ||\sum_{i \in I} \bar{w}_i \vc{a}^i||$.
The result follows from Proposition~\ref{prop:class-0 vicinity}.
\Halmos 
\endproof

Unlike the above results for class-0 vectors that can be directly obtained from a conic combination of vectors describing $\mathcal{C}\pol$, the concept of class-1 vectors close to the boundary of $\mathcal{C}\pol$ requires further development.
In particular, we need to understand how to perturb a given vector $\vc{u}$ on the boundary of $\mathcal{C}\pol$ to obtain vectors in the vicinity of $\vc{u}$ but outside of the polar cone.
In the following derivations, we define $\mathcal{C}^J = \mathcal{C} \cap \{\vc{x} \in \Re^n \, | \, \vc{a}^j \vc{x} = 0, \forall j \in J \}$ for any $J \subseteq I$.
We set that $\mathcal{C}^{\emptyset} = \mathcal{C}$.

\begin{proposition} \label{prop:class1}
	Let $\vc{u}$ be a point on the boundary of $\mathcal{C}\pol$, i.e., there exists $J \subseteq I$ such that $\vc{u} = \sum_{j \in J} v_j \vc{a}^j$ where $v_j > 0$ for all $j \in J$.	 
	Let $K \subseteq I$ be a subset of constraint indices such that $\mathcal{C}^J \nsubseteq \mathcal{C}^K$.	
	Then, $\vc{u} - \epsilon \vc{d}(\vc{w})$ belongs to class-1 for any $\epsilon > 0$ and any $\vc{w} \in \Re^{m+n}$ such that $w_{k} > 0$ for all $k \in K$, and $w_{k} \geq 0$ for all $k \in I \setminus J$.
\end{proposition}

\proof{Proof.}
First, we verify the correctness of the opening statement of the proposition, i.e., there exists $J \subseteq I$ such that $\vc{u} = \sum_{j \in J} v_j \vc{a}^j$ where $v_j > 0$ for all $j \in J$.	
By definition, $\vc{u}$ can be obtained as a conic combination of extreme rays of $\mathcal{C}\pol$ as the vector is on its boundary.
The fact that $\mathcal{C}\pol$ is the cone generated by the constraint vectors of $\mathcal{C}$ implies that the set of extreme rays of $\mathcal{C}\pol$ is a subset of $\{\vc{a}^i\}_{i \in I}$.
As a result, the constraint indices with positive coefficients in the conic combination will form $J$.
Note that if $\vc{u}$ is the origin, we have $J = \emptyset$.
For the second part of the proposition, pick any $\epsilon > 0$ and any $\vc{w} \in \Re^{m+n}$ such that $w_{k} > 0$ for all $k \in K$, and $w_{k} \geq 0$ for all $k \in I \setminus J$. 
To prove that $\vc{u} - \epsilon \vc{d}(\vc{w})$ belongs to class-1, i.e., $\max \{\left(\vc{u} - \epsilon \vc{d}(\vc{w})\right) \vc{x} \, | \, \vc{x} \in \mathcal{C} \} = \infty$, it suffices to show that $\left(\vc{u} - \epsilon \vc{d}(\vc{w})\right) \vc{x} > 0$ for some $\vc{x} \in \mathcal{C}$.
It follows from the assumption $\mathcal{C}^J \nsubseteq \mathcal{C}^K$ that there exists a point $\vc{x}^* \in \mathcal{C}^J \subseteq \mathcal{C}$ such that $\vc{x}^* \notin \mathcal{C}^K$, i.e., $\vc{a}^{k^*} \vc{x}^* < 0$ for some $k^* \in K$.
We write that
\begin{align*}
	\left(\vc{u} - \epsilon \vc{d}(\vc{w})\right) \vc{x}^* & = \vc{u} \vc{x}^* - \epsilon \vc{d} (\vc{w}) \vc{x}^* \\
	& = \sum_{j \in J} v_j \vc{a}^j \vc{x}^* - \epsilon \sum_{i \in I} \bar{w}_i \vc{a}^i \vc{x}^* \\
	& = \sum_{j \in J} v_j \vc{a}^j \vc{x}^* - \epsilon \sum_{i \in I \setminus \{k^*\}} \bar{w}_i \vc{a}^i \vc{x}^* - \epsilon \bar{w}_{k^*} \vc{a}^{k^*} \vc{x}^* > 0,
\end{align*}
where $\bar{w}_i = \frac{w_i}{||\sum_{i \in I} w_i \vc{a}^i||}$ for $i \in I$. In the above equations, the second equality is implied by the definitions of $\vc{u}$ and $\vc{d}(\vc{w})$ as a normalized combination of constraint vectors, and the last inequality follows from the facts that $\vc{a}^j \vc{x}^* = 0$ for all $j \in J$, $\vc{a}^i \vc{x}^* \leq 0$ and $\bar{w}_i \geq 0$ for all $i \in I \setminus J$, and $\vc{a}^{k^*} \vc{x}^* < 0$ with weight $\bar{w}_{k^*} > 0$ as $k^* \in K$.  
\Halmos 
\endproof

Proposition~\ref{prop:class1} gives a method to generate class-1 vectors arbitrarily close to the boundary of the polar cone $\mathcal{C}\pol$.
\black{In this proposition, similarly to the results of Proposition~\ref{prop:class-0 vicinity}, $\epsilon$ controls the \textit{closeness} of the generated class-1 vectors to the boundary of the polar cone.}
We next show that all class-1 vectors can be generated through a similar method relative to the boundary of $\mathcal{C}\pol$; hence providing a necessary and sufficient condition to generate class-1 vectors analogous to those of class-0 vectors given in Proposition~\ref{prop:class-0 vicinity}.

\begin{proposition} \label{prop:class1_converse}
	Let $\vc{c}$ be a vector of class-1.
	Then, there exists 
	\begin{itemize}
		\item [(i)] a vector $\vc{u}$ on the boundary of $\mathcal{C}\pol$, i.e., $\vc{u} = \sum_{j \in J} v_j \vc{a}^j$ for some $J \subseteq I$ where $v_j > 0$ for all $j \in J$,
		\item [(ii)] a subset $K \subseteq I$ of constraint indices such that $\mathcal{C}^J \nsubseteq \mathcal{C}^K$,
		\item [(iii)] and parameters $\epsilon > 0$ and $\vc{w} \in \Re^{m+n}$ such that $w_{k} > 0$ for all $k \in K$, and $w_{k} \geq 0$ for all $k \in I \setminus J$,
	\end{itemize}	 
	such that $\vc{c} = \vc{u} - \epsilon \vc{d}(\vc{w})$.
\end{proposition} 

\proof{Proof.}
Since $\vc{c}$ is a vector of class-1, it does not belong to $\mathcal{C}\pol$ by Proposition~\ref{prop:polar}.
It is easy to verify that $\mathcal{C}\pol$ is full-dimensional because $\mathcal{C}\pol$ is the cone generated by vectors $\vc{a}^i$ for $i \in I$ which includes all negative unit vectors $-\vc{e}^j$ for $j \in N$.
Pick a point $\vc{b}$ in the interior of $\mathcal{C}\pol$, i.e., $\vc{b} = \sum_{l \in L} t_l \vc{a}^l$ for some $L \subseteq I$ where $t_l > 0$ for all $l \in L$, and vectors $\vc{a}^l$ are linearly independent.
Let point $\vc{u}$ denote the intersection of the line connecting $\vc{c}$ and $\vc{b}$ with the boundary of $\mathcal{C}\pol$.
It follows from the proof of Proposition~\ref{prop:class1} that $\vc{u} = \sum_{j \in J} v_j \vc{a}^j$ for some $J \subseteq I$ where $v_j > 0$ for all $j \in J$, satisfying condition (i).
Define $K = L \setminus J$.
Note that $K$ must include an index $k^*$ whose corresponding vector $\vc{a}^{k^*}$ does not belong to the linear space spanned by the vectors in $J$, since otherwise $\vc{b}$ would be on a face of $\mathcal{C}\pol$, and not in its interior.
This implies that there exists a vector $\bar{\vc{x}}$ in the linear space orthogonal to the space spanned by the vectors in $J$, but not orthogonal to $\vc{a}^{k^*}$.
In particular, $\vc{a}^{j} \bar{\vc{x}} = 0$ for $j \in J$, and $\vc{a}^{k^*} \bar{\vc{x}} <0 $ as $\mathcal{C}\pol$ is a convex cone, and thus the normal vector of one of its facets makes an obtuse angle to any vector in the cone that is not on that facet.
As a result, we have that $\mathcal{C}^J \nsubseteq \mathcal{C}^K$, satisfying condition (ii).
We can write $\vc{c}$ as $\vc{u} - ||\vc{c} - \vc{u}||\frac{(\vc{b} - \vc{u})}{||\vc{b} - \vc{u}||}$ by construction.
Define $\vc{d}(\vc{w}) = \frac{\vc{b} - \vc{u}}{||\vc{b} - \vc{u}||}$ where $\vc{w} \in \Re^{m+n}$ is defined as $w_k = t_k - v_k$ for $k \in J \cap L$, $w_k = - v_k$ for $k \in J \setminus L$, $w_k = t_k$ for $k \in K$, and $w_k = 0$ for $k \in I \setminus K \cup L$.
Finally, setting $\epsilon = ||\vc{c} - \vc{u}||$ will complete the proof as condition (iii) is satisfied.
\Halmos 
\endproof	

The difference between generating class-0 vectors through Proposition~\ref{prop:class-0 vicinity} and class-1 vectors through Proposition~\ref{prop:class1} is in the efficiency of the approach.
While we can efficiently obtain class-0 vectors through any arbitrary conic combination of vectors $\vc{a}^{i}$, we still need to determine indices $J$ and $K$ in Proposition~\ref{prop:class1} to produce class-1 vectors---a process that can be computationally expensive for a large data set.
For our classification approach, it is critical to generate training data sets efficiently.
In other words, we need to select vectors whose class is identified immediately without solving the CGLP, since otherwise the process of determining the class of a set of arbitrary vectors would be computationally demanding.
\black{The efficiency of generating training data sets is particularly important when implementing the classification approach inside of \BB framework, where training is done \textit{online}, at each layer of the \BB tree, during the optimization process; see Section~\ref{sec:consistency} for more details on such implementations.}  
We next give several results that help generate class-1 vectors more efficiently by applying findings of Propositions~\ref{prop:class1} and \ref{prop:class1_converse} to special cases.
We start with a case where the negative of a conic combination of constraint vectors yields a class-1 vector as long as one of the combination weights is non-zero.

\begin{corollary} \label{cor:class1_special_1}
	Let $k^* \in I$ be a constraint index such that $\mathcal{C}^{\{k^*\}} \neq \mathcal{C}$.
	Then, for any weight vector $\vc{w} \in \Re^{m+n}_+$ such that $w_{k^*} > 0$, the vector $-\vc{d}(\vc{w})$ belongs to class-1. 
\end{corollary}

\proof{Proof.}
Pick $\vc{u} = \vc{0}$ on the boundary of $\mathcal{C}\pol$.
According to Proposition~\ref{prop:class1}, we have $J = \emptyset$.
Note that $\mathcal{C}^J = \mathcal{C}$.
Set $K = \{k^*\}$.
We have that $\mathcal{C}^K \subseteq \mathcal{C}^J$ by definition, as the former set is a face of the latter set.
Therefore, the assumption that $\mathcal{C}^K \neq \mathcal{C}^J$ implies that $\mathcal{C}^J \nsubseteq \mathcal{C}^K$.
The result follows from Proposition~\ref{prop:class1} by setting $\epsilon = 1$. 
\Halmos 
\endproof	

It follows from Corollary~\ref{cor:class1_special_1} that the negative of any constraint that does not contain the entire cone $\mathcal{C}$ is a class-1 vector.
This is a simple and quick approach to generate class-1 vectors in practice.
This derivation can be streamlined further under the full-dimensionality assumption, as described next.

\begin{corollary} \label{cor:class1_special_2}
	Assume that $\mathcal{C}$ is full-dimensional.
	Then, for any nonzero weight vector $\vc{w} \in \Re^{m+n}_+$, the vector $-\vc{d}(\vc{w})$ belongs to class-1. 
\end{corollary}

\proof{Proof.}
Since $\mathcal{C}$ is full-dimensional, it is not contained in the face defined by any of its constraints, i.e., $\mathcal{C}^{\{k\}} \neq \mathcal{C}$ for all $k \in I$.
The result follows from Corollary~\ref{cor:class1_special_1} by setting $k^*$ to be the index of a nonzero component of $\vc{w}$.
\Halmos 
\endproof	

The next corollary prescribes a method to produce class-1 vectors for a non-zero cone.

\begin{corollary} \label{cor:class1_special_3}
	Under the assumption that $\mathcal{C} \neq \{0\}$, the vector $-\vc{d}(\vc{w})$ belongs to class-1 for any weight vector $\vc{w} \in \Re^{m+n}_+$ such that $w_{m+i} > 0$ for all $i \in N$. 
\end{corollary}

\proof{Proof.}
Pick $\vc{u} = \vc{0}$ on the boundary of $\mathcal{C}\pol$.
According to Proposition~\ref{prop:class1}, we have $J = \emptyset$, and hence $\mathcal{C}^J = \mathcal{C}$.
Set $K = \{m+1, \dotsc, m+n\}$.
It follows that $\mathcal{C}^K = \{0\}$ as $\bar{x}_k = 0$ for all $k \in K$ for any point $\bar{\vc{x}} \in \mathcal{C}^K$.
Since $\mathcal{C}^J \neq \mathcal{C}^K$, and considering the fact that $\mathcal{C}^K \subseteq \mathcal{C}^J$, we conclude that $\mathcal{C}^J \nsubseteq \mathcal{C}^K$.
The result follows from Proposition~\ref{prop:class1} by setting $\epsilon = 1$.
\Halmos 
\endproof	

The next result proposes an efficient method to generate vectors of class-1 for non-redundant constraints.

\begin{corollary} \label{cor:class1_special_4}
	Let $j^* \in I$ be the index of a non-redundant constraint of $\mathcal{C}$.
	Then, the vector $\vc{a}^{j^*} - \epsilon \vc{d}(\vc{w})$ belongs to class-1 for any $\epsilon > 0$ and any $\vc{w} \in \Re^{m+n}$ that satisfies any of the following conditions:
	\begin{itemize}
		\item [(i)] $w_{k} \geq 0$ for all $k \in I \setminus \{j^*\}$, and $w_{k} > 0$ for all $k = m+i$ with $i \in N$,
		\item [(ii)] $w_{k} \geq 0$ for all $k \in I \setminus \{j^*\}$, and $w_{k} > 0$ for all $k = m+i$ with $i \in N$ such that $a^{j^*}_i \geq 0$,
		\item [(iii)] $w_{k} \geq 0$ for all $k \in I \setminus \{j^*\}$, and $w_{k} > 0$ for all $k = m+i$ with $i \in N$ such that $a^{j^*}_i \leq 0$.
	\end{itemize}
\end{corollary}

\proof{Proof.}
Since $\vc{a}^{j^*} \vc{x} \leq 0$ is not redundant, $\vc{a}^{j^*}$ cannot be represented as a conic combination of distinct constraint vectors of $\mathcal{C}$.
As a result, $\vc{a}^{j^*}$ is a vector on the boundary of $\mathcal{C}\pol$.
Set $\vc{u} = \vc{a}^{j^*}$ and $J = \{j^*\}$.
\begin{itemize}
	\item [(i)] Set $K = \{m+i\}_{i \in N}$.
	It follows that $\mathcal{C}^K = \{0\}$.
	Since $\vc{a}^{j^*} \vc{x} \leq 0$ is non-redundant, the restriction of $\mathcal{C}$ at this constraint must be non-zero.
	Therefore, $\mathcal{C}^J \nsubseteq \mathcal{C}^K$.
	The result follows from Proposition~\ref{prop:class1}.
	\item [(ii)] Set $K = \{m+i\}$ for all $i \in N$ such that $a^{j^*}_i \geq 0$.
	Non-redundancy of $\vc{a}^{j^*} \vc{x} \leq 0$ implies that there exists a non-zero point $\bar{\vc{x}} \in \mathcal{C}$ such that $\vc{a}^{j^*} \bar{\vc{x}} = 0$, i.e., $\bar{\vc{x}} \in \mathcal{C}^J$.
	It follows from $\bar{\vc{x}} \neq \vc{0}$ that there exists a coordinate $i^* \in N$ such that $\bar{x}_{i^*} > 0$ and $a^{j^*}_{i^*} \geq 0$, since otherwise the point would not satisfy $\vc{a}^{j^*} \bar{\vc{x}} = 0$.
	Note that $\mathcal{C}^K$ is the restriction of $\mathcal{C}$ at the hyperplane $\mathcal{H} = \{\vc{x} \in \Re^n \, | \, x_i = 0, \forall i \in N \text{ such that }a^{j^*}_i \geq 0\}$. Therefore, $\bar{\vc{x}} \notin \mathcal{C}^K$ because $i^*$ is one of the indices of variables that are fixed at zero in the definition of $\mathcal{H}$.
	We conclude that $\mathcal{C}^J \nsubseteq \mathcal{C}^K$ as $\bar{\vc{x}} \in \mathcal{C}^J$ and $\bar{\vc{x}} \notin \mathcal{C}^K$.
	The result follows from Proposition~\ref{prop:class1}.
	\item [(iii)] This result follows from an argument similar to that of condition (ii) by setting $K = \{m+i\}$ for all $i \in N$ such that $a^{j^*}_i \leq 0$.
\end{itemize}
\Halmos 
\endproof

\section{Implementation in Branch-and-Cut.} \label{sec:consistency}

The framework of Section~\ref{sec:classification} can be applied to any CGLP at the nodes of the \BB tree to produce different classes of cutting planes, from disjunctive cuts \citep{Bal85} to reformulation-and-linearization technique (RLT) cuts \citep{SheAda94} through an appropriate projection.
To illustrate, in this paper, we focus on a CGLP that generates a new class of cutting planes called \textit{consistency cuts}.
The application to other classes of CGLP follows similarly; see Remark~\ref{remark:general consistency} at the end of this section.
It is shown in \citep{davarnia:ra:ho:2022}, both theoretically and computationally, that consistency cuts can be more effective than the classical cutting planes such as RLT cuts in reducing the total solution time of the branch-and-cut process by directly targeting the size reduction in the \BB tree.
Inspired by these results, we aim to assess the effectiveness of the proposed machine learning method in further improving this solution time.
First, we give a brief introduction to the concept of consistency in integer programming (IP) and provide results for constructing the CGLP; see \citep{davarnia:ra:ho:2022} for a detailed account.
Then, we show how to transform the CGLP structure to make it amenable to the machine learning framework of the previous section.

In the following, since the CGLP model has the central role in our paper, we use variables $(\vc{x},\vc{y})$ for the CGLP and variables $(\vc{\alpha}, \vc{\beta})$ for the original IP, which is often reversed in the IP literature.
Consider a 0--1 set $\S = \{\vc{\alpha} \in \{0,1\}^n \, | \, A\vc{\alpha} \leq \vc{b}\}$ for some matrix $A$ and vector $\vc{b}$ of proper dimensions.
We refer to the LP relaxation of $\S$ as $\S_{\LP}$.
For any $J\subseteq N=\{1,\ldots,n\}$, let $\vc{\alpha}_J$ be the tuple containing the variables in $\{\alpha_j\;|\;j\in J\}$.
A {\em partial assignment} to $\vc{\alpha}$ is a 0--1 assignment of values to $\vc{\alpha}_J$ for some $J\subseteq N$.
A partial assignment $\vc{\alpha}_N$ is referred to as a full solution.
The following definition of LP-consistency ensues.

\begin{definition}
	We say that a 0--1 set \black{$\S = \{\vc{\alpha} \in \{0,1\}^n \, | \, A\vc{\alpha} \leq \vc{b}\}$} is \textit{LP-consistent} if any partial assignment $\vc{\alpha}_J = \vc{v}_J$ with $J \subseteq N$ such that $\S_\LP \cap \{\vc{\alpha} \in \Re^n: \vc{\alpha}_J = \vc{v}_J\} \neq \emptyset$ can be extended to a full feasible solution of $\S$.
\end{definition}

This definition implies that an LP-consistent set completely eliminates \textit{backtracking} during the \BB search because the LP relaxation at a node of the \BB tree will be infeasible if it does not contain any integer solution, which makes the algorithm prune that node right away without creating new branches.
It is known that modifying a general set to make it LP-consistent is NP-hard \cite{davarnia:ra:ho:2022}.
As a result, a generalized variant of LP-consistency is proposed that controls the level of backtracking elimination through the concept of \textit{ranks}.
This structure provides a substantial practical flexibility, as reducing the rank of LP-consistency leads to more computationally-affordable implementations.

\begin{definition} \label{def:partialLP-consistency}
	Let \black{$\S = \{\vc{\alpha} \in \{0,1\}^n \, | \, A\vc{\alpha} \leq \vc{b}\}$} be a 0--1 set, and consider a subset $I \subseteq N$.
	Then $\S$ has \textit{partial \mbox{LP-consistency} of rank $r$ over $I$}, where $0 < r \leq n - |I|$, if for every partial assignment $\vc{\alpha}_I = \vc{v}_I$ such that $\S_{\LP} \cap \{\vc{\alpha} \in \Re^n: \vc{\alpha}_I = \vc{v}_I\} \neq \emptyset$, and every $J\subseteq N \setminus I$ with $|J|= r$, there exists a 0--1 value assignment $\vc{\alpha}_J = \vc{v}_J$ such that $\S_{\LP} \cap \{\vc{\alpha} \in \Re^n: \vc{\alpha}_{I\cup J} = \vc{v}_{I \cup J}\} \neq \emptyset$.	
\end{definition}

Intuitively, partial LP-consistency of rank $r$ means that there is no backtracking in the next $r$ levels of the \BB tree.
Algorithm~\ref{alg:partialLP-consistency} provides a systematic procedure to achieve partial LP-consistency.

\begin{algorithm}[!b]                    				  	
	\caption{Algorithm to achieve partial LP-consistency}          
	\label{alg:partialLP-consistency}			                        
	\begin{algorithmic}[1]                				    
		\REQUIRE A 0--1 set \black{$\S = \{\vc{\alpha} \in \{0,1\}^n \, | \, A\vc{\alpha} \leq \vc{b}\}$}, a subset $I \subseteq N$, and a positive number $r \leq n - |I|$
		\ENSURE A set $\widehat{S}$ that has partial LP-consistency of rank $r$ over $I$
		\STATE Initialize $\widehat{\S} = \S$
		\STATE{\black{Generate the nonlinear augmented system $(A\vc{\alpha} - b)\prod_{j \in J_1} \alpha_j \prod_{j \in J \setminus J_1} (1-\alpha_j) \leq 0$ for all $J\subseteq N$ with $|J|=r$, and $J_1\subseteq J$}}
		\STATE Linearize the above system by replacing $\alpha_i^2$ with $\alpha_i$, and $\prod_{k \in K} \alpha_k$ with new variable $\beta_{K}$ for each $K \subseteq N$ that appears in the system. Denote the resulting linear constraint set by $\mathcal{R}(\S_{\LP})$, and its projection onto the space of variables $\vc{\alpha}_I$ by $\mathcal{R}(\S_{\LP})|_I$.
		\STATE Add to $\widehat{\S}$ the inequalities in $\mathcal{R}(\S_{\LP})|_I$.		
	\end{algorithmic}
\end{algorithm}

Since generating all the inequalities describing the projection $\mathcal{R}(\S_{\LP})|_I$ can be computationally prohibitive, a CGLP is designed to produce a small subset of these inequalities, which are referred to as \textit{consistency cuts}, through separation.

\begin{proposition} \label{prop:CGLP}
	\black{Consider a 0--1 set $\S = \{\vc{\alpha} \in \{0,1\}^n \, | \, A\vc{\alpha} \leq \vc{b}\}$ whose constraints include the bounds on variables.}
	Select a subset $I \subseteq N$, and a positive number $0 < r \leq n - |I|$.
	Let $\mathcal{R}(\S_{\LP}) = \{ A_I \vc{\alpha}_I + A_J \vc{\alpha}_J + B \vc{\beta} \leq \vc{d} \}$ be the lifted linear system generated in Algorithm~\ref{alg:partialLP-consistency}. 
	Consider the following CGLP
	\begin{subequations}
		\begin{align}
			w^*(v_I) = \max \quad & (A_I \vc{v}_I - \vc{d})  \vc{x}   \label{eq:CGLP-1}\\
			\text{s.t.} \quad & A_J \vc{x} = \vc{0}  \label{eq:CGLP-2}\\
			& B \vc{x} = \vc{0}  \label{eq:CGLP-3}\\
			& \vc{x} \geq \vc{0} \label{eq:CGLP-5}
		\end{align}
	\end{subequations}
	where variables $\vc{x}$ represent the dual weight vector associated with constraints of $\mathcal{R}(\S_{\LP})$.
	Define $\hat{\S}$ to be the set obtained by adding to $\S$ constraints of the form \mbox{$\bar{\vc{x}} A_I \vc{\alpha}_I \leq \bar{\vc{x}} \vc{d}$}, where $\bar{\vc{x}}$ is an optimal ray of \eqref{eq:CGLP-1}--\eqref{eq:CGLP-5}, for all $\vc{v}_I\in\{0,1\}^{|I|}$ that yield the optimal value $w^*(v_I) = \infty$.
	
	Then $\hat{\S}$ has partial \mbox{LP-consistency} of rank $r$ or higher over $I$.
	\Halmos
\end{proposition}

We refer the interested reader to the discussion in Appendix C and the references therein for a detailed account on the possibility of achieving LP-consistency of ranks higher than $r$ through the above proposition.

\begin{remark} \label{remark:consistency}
	The consistency cuts produced in Proposition~\ref{prop:CGLP} have an important computational advantage compared to traditional cutting planes such as RLT cuts.
	When $w^*(v_I) > 0$, the purpose of adding the valid inequality $\bar{\vc{x}} A_I \vc{\alpha}_I \leq \bar{\vc{x}} \vc{d}$ is to separate the partial assignment $\vc{\alpha}_I = \vc{v}_I$.
	Alternatively, this solution can also be separated by adding a so-called \textit{clausal inequality} of the form $\sum_{i \in I: v_i = 0} \alpha_i + \sum_{i \in I: v_i = 1} (1-\alpha_i) \geq 1$.
	As a result, the consistency cuts can be produced based only on the optimal value of the CGLP \eqref{eq:CGLP-1}--\eqref{eq:CGLP-5} without the need for an optimal solution.
	This property is particularly useful for the machine learning approach of Section~\ref{sec:classification}, where the optimal value of the CGLP is approximated through classification.
\end{remark}

\begin{remark} \label{remark:pruning}
	Consistency cuts are globally valid, implying that they can be generated and added at any point during the \BB process, including at the root node.
	In the event that they are generated locally at a node of the \BB tree based on the particular partial assignment in that node, they can be alternatively viewed as a pruning technique. Instead of generating a clausal inequality at a node of the \BB tree, one can immediately prune that node.
\end{remark}

The CGLP of Proposition~\ref{prop:CGLP} can be solved at the nodes of the \BB tree to produce consistency cuts.
Such implementation, however, can be computationally cumbersome due to the large number of node candidates to invoke the CGLP at.
To circumvent this difficulty, the machine learning method can be used to approximate the optimal value of the CGLP without the need to solve it, thereby reducing the implementation time.
To this end, the results obtained for $\mathcal{C}$ in Section~\ref{sec:classification} can be applied to the projection cone \eqref{eq:CGLP-2}--\eqref{eq:CGLP-5} by breaking the equality constraints into two inequalities with opposite signs.
However, this approach has two computational drawbacks: (i) it increases the size of the problem which can result in a larger training data set and longer training time; and (ii) such representation limits our options to generate class-1 vectors as the assumptions of Corollaries~\ref{cor:class1_special_1} and \ref{cor:class1_special_2} would not be satisfied.
To mitigate these drawbacks, we next show how the equality constraints in the CGLP of Proposition~\ref{prop:CGLP} can be represented as inequalities without the need to break them.

\begin{proposition} \label{prop:new_CGLP}
	Consider the setting of Proposition~\ref{prop:CGLP}, and define 
	\begin{subequations}
		\begin{align}
			\bar{w}^*(v_I) = \max \quad &  (A_I \vc{v}_I - \vc{d}) \vc{x}   \label{eq:New_CGLP-1}\\
			\text{s.t.} \quad & A_J \tr \vc{x} \geq \vc{0}  \label{eq:New_CGLP-2}\\
			& B \tr \vc{x} \geq \vc{0}  \label{eq:New_CGLP-3}\\
			& \vc{x} \geq \vc{0}. \label{eq:New_CGLP-5}
		\end{align}
	\end{subequations} 
	Then, $w^*(v_I) = \bar{w}^*(v_I)$ for all $\vc{v}_I \in \{0,1\}^{|I|}$. 
\end{proposition}

\proof{Proof.}
Define $\bar{\mathcal{R}}(\S_{\LP}) = \{ A_I \vc{\alpha}_I + A_J \vc{\alpha}_J + B \vc{\beta} \leq \vc{d}, \, -\vc{\alpha}_J \leq \vc{0}, \, -\vc{\beta} \leq \vc{0} \}$.
It is easy to verify that $\bar{\mathcal{R}}(\S_{\LP}) = \mathcal{R}(\S_{\LP})$ since the added non-negativity bounds are implied by the set of constraints generated in 
Algorithm~\ref{alg:partialLP-consistency}.
As a result, $\bar{\mathcal{R}}(\S_{\LP})|_I = \mathcal{R}(\S_{\LP})|_I$.
The CGLP corresponding to $\mathcal{R}(\S_{\LP})|_I$ is given by \eqref{eq:CGLP-1}--\eqref{eq:CGLP-5}.
We need to show that the CGLP corresponding to $\bar{\mathcal{R}}(\S_{\LP})|_I$ is obtained by \eqref{eq:New_CGLP-1}--\eqref{eq:New_CGLP-5}.
Let $\vc{x}$, $\vc{y}$ and $\vc{z}$ be the dual vectors corresponding to the constraints in $\bar{\mathcal{R}}(\S_{\LP})$.
The CGLP associated with $\bar{\mathcal{R}}(\S_{\LP})|_I$ is as follows.
\begin{subequations}
	\begin{align}
		\bar{w}^*(v_I) = \max \quad & (A_I \vc{v}_I - \vc{d}) \vc{x}   \label{eq:N_CGLP-1}\\
		\text{s.t.} \quad & \left[
		\def\arraystretch{1.5}
		\begin{array}{c|c|c}
			A_J \tr \, & \, -I \, & \, 0 \\
			\hline
			B \tr & 0 & -I
		\end{array}
		\right]
		\left[
		\def\arraystretch{1.5}
		\begin{array}{c}
			\vc{x} \\
			\hline
			\vc{y} \\
			\hline
			\vc{z}
		\end{array}
		\right] = 
		\left[
		\def\arraystretch{1.5}
		\begin{array}{c}
			\vc{0} \\
			\hline
			\vc{0}
		\end{array}
		\right] \label{eq:N_CGLP-2}\\
		& \vc{x}, \vc{y}, \vc{z} \geq \vc{0}, \label{eq:N_CGLP-5}
	\end{align}
\end{subequations} 
where the first and second row blocks in \eqref{eq:N_CGLP-2} represent the constraints associated with the projection of variables $\vc{\alpha}_J$ and $\vc{\beta}$, respectively.
The column blocks of the coefficient matrix in \eqref{eq:N_CGLP-2} represent variables $\vc{x}$, $\vc{y}$ and $\vc{z}$, respectively.
This problem structure allows for removing variables $\vc{y}$ and $\vc{z}$ from the feasible region and replacing the equality constraints with $\geq$ inequalities, yielding the CGLP \eqref{eq:New_CGLP-1}--\eqref{eq:New_CGLP-5}.
\Halmos
\endproof

In view of Proposition~\ref{prop:new_CGLP}, for any component in the vectors $\vc{\alpha}_J$ and $\vc{\beta}$ that obtains a non-negativity bound as a result of the inequalities generated in Algorithm~\ref{alg:partialLP-consistency}, i.e., the bound is included in $A_I \vc{\alpha}_I + A_J \vc{\alpha}_J + B \vc{\beta} \leq \vc{d}$, we can replace the dual variables $\vc{y}$ and $\vc{z}$ corresponding to that inequality with the $\vc{x}$ coordinate in the proof of Proposition~\ref{prop:new_CGLP}.
Consequently, following the closing arguments in the proof, this $\vc{x}$ variable will be removed from the constraints in \eqref{eq:New_CGLP-1}--\eqref{eq:New_CGLP-5}, reducing the size of the CGLP.
We further note that the above arguments can be used to transform equality constraints to inequalities for any CGLP that has a \textit{lift-and-project} type structure such as those producing disjunctive, split, and RLT cuts.
We present computational results for using the machine learning approach to generate consistency cuts inside of \BB in Section~\ref{subsec:IP}.

\begin{remark} \label{remark:general consistency}
	As deduced from Algorithm~\ref{alg:partialLP-consistency}, the consistency cuts presented in this section can be obtained from a special projection of the CGLP associated with the RLT cuts as one of the most general forms of the CGLP.
	Such projections can be used in any other CGLP model to produce a variant of the consistency cuts corresponding to that CGLP.
	As a result, the machine learning framework can be applied in a similar way to any form of the CGLP that is used during branch-and-cut.
\end{remark}

We conclude this section by highlighting the possible approaches to incorporate our proposed framework in the branch-and-cut process.
In general, there are two ways that the outcome of the approximate CGLP can be used at the layers of the \BB process.
The first approach is a \textit{greedy} application, where the outcome of the classification replaces that of the CGLP entirely.
This approach can lead to generating invalid cutting planes at a node associated with a falsely positive classification value.
In this case, the solution of the branch-and-cut process is not guaranteed to be optimal, but it can be used as a fast heuristic to obtain feasible solutions, and hence a primal bound.
The second approach is a \textit{conservative} application, where the outcome of the classification is used as a filtering mechanism to identify the best node candidates that are most likely to admit cutting planes.
At these filtered node candidates, the CGLP is solved to determine whether a cutting plane can be added.
This approach saves time by avoiding solving the CGLP at the nodes that have lower likelihood for admitting cutting planes, while preserving optimality guarantee of the end solution.

\section{Computational Results.} \label{sec:computation}

In this section, we present preliminary computational experiments to show the potential of the proposed machine learning approach in improving the CGLP solution time.
These results will be divided in two parts.
The first part, as given in Section~\ref{subsec:chm}, evaluates the outcome of the machine learning method when applied to the CHM problem. 
These results not only represent an application of our proposed framework to combinatorial problems beyond cut-generation, they will also allow us to study the performance of the approach on benchmark instances in \textit{direct} comparison to the outcome of modern solvers in a controlled environment without being impacted by \BB elements.
The second part, as given in Section~\ref{subsec:IP}, assesses the effectiveness of our approach to produce cutting planes in comparison with the conventional CGLP implementation inside of \BB for benchmark IP instances.
These models are coded in Python v2.7 on a computer with specifications Intel Core i7 1.8 GHz computer processing unit (CPU) and 8 GB of random access memory (RAM).
For the machine learning experiments, the scikit-learn package Version 1.1 has been used. 
The optimization problems are solved by the CPLEX solver Version 12.8.0 at its default setting, unless stated otherwise.

\subsection{Convex Hull Membership Problem.} \label{subsec:chm}
As discussed in Section~\ref{sec:classification}, our goal is to use classification methods in machine learning to separate class-0 and class-1 vectors $\vc{c}$ corresponding to the CGLP defined by $z^*(\vc{c}) = \max \{\vc{c} \vc{x} \, | \, \vc{x} \in \mathcal{C} \}$ where $\mathcal{C} = \{\vc{x} \in \Re^n \, | \, \vc{a}^i \vc{x} \leq 0, \, \forall i \in I \}$ and $I = \{1, \dotsc, m, m+1, \dotsc, m+n\}$.
To generate test instances for the CGLP, we use the data characterization given in \cite{kalantari:zh:2022} for the CHM problem, as solving the CHM reduces to that of a CGLP; see Remark~\ref{remark:CHM}.
We consider three categories for the problem size: (i) $n = 300 $, $m = 100 $; (ii) $n= 1000 $, $m = 200$; (iii) $n = 2000 $, $m = 500$.
For each size category, we create 5 random instances with the following specifications.
Adapting the method used in \cite{kalantari:zh:2022} to generate point coordinates for the CHM problem, we create the CGLP's coefficient matrix randomly from a discrete uniform distribution between $[-100, 100]$ with \black{the conditions that (i) the coefficients of a constraint cannot be all positive since otherwise $\mathcal{C}$ would contain the origin only; and (ii) the coefficients of a constraint cannot be all negative since otherwise the constraint would be implied by the nonnegativity constraints.}
We normalize the resulting vectors that form the polar cone $\mathcal{C}^{\pol}$ to make their magnitudes uniform, while maintaining the constraint structure in $\mathcal{C}$.   

We aim to approximate the indicator function of the CGLP by a classifier that is trained over a pre-determined set of vectors in both classes. 
While increasing the number of data points in the training set can improve the training accuracy, it also leads to an increase in the training time.
Based on our experiments (see Appendix B), setting the number of data points proportional to the number of constraints of the CGLP model, including the variable bounds, leads to a reasonable trade-off, which is chosen for these experiments.
More specifically, the size of the training set for each test instance is chosen to be $4(m+n)$.
These training data sets are generated based on the results developed in Section~\ref{sec:classification}.
In particular, we create $(m+n)$ class-0 data points through using Corollary~\ref{cor:class-0} for each constraint $i \in I$ of $\mathcal{C}$ by setting components of vector $\vc{v}$ randomly from a discrete uniform distribution between $[0, 10000]$ for each $i \in I$, and setting $\bar{\epsilon} = 0.001$.
We create another $(m+n)$ class-0 vectors according to the above method.
We create $(m+n)$ class-1 training data points through using Corollary~\ref{cor:class1_special_1} to include the negative of each constraint vector with index $k^* \in I$ by setting the weight vector $\vc{w} = \vc{e}^{k^*}$.
Note that the assumption of this corollary is satisfied as long as the entire cone $\mathcal{C}$ is not restricted to the hyperplane defined by $\vc{a}^{k^*} \vc{x} = 0$.
We also create another $(m+n)$ class-1 training data points by applying Corollary~\ref{cor:class1_special_4}(i) to each constraint index $j^* \in I$ with $w_{m+i} = 1$ for all $i \in N$, $w_{k} = 0$ for all other components, and $\epsilon = 0.001$.
Since this corollary requires a non-redundancy assumption on constraints of $\mathcal{C}$, we study sensitivity of the classifier to changes in the representation of $\mathcal{C}$ by adding redundancy; see Appendix B.
That section also includes the sensitivity analysis for the classifier accuracy with respect to the proximity of the training data to the boundary of the polar cone.

For the test set, we generate class-1 vectors randomly with components chosen from a uniform distribution between $[-100, 100]$ and then normalize them.
The class of these vectors is verified through solving the CGLP.
Class-0 vectors in the test set are generated through a conic combination of constraint vectors $\vc{a}^i$ for $i \in I$ with weights being randomly selected from a uniform distribution between $[0, 100]$.
The size of the test set is equal to $\sim 30\%$ of the size of the training set for each problem.

We use different classifiers in machine learning to approximate the indicator function of the CGLP.
For a detailed account on the machine learning methods that can be used for function approximation, we refer the readers to \cite{hammer:ge:2003, wasserman:2006}.
Table~\ref{tab: chm summary} reports the summary of the classification results for different problem sizes and classification methods, namely the random forest (RF), the $k$-nearest neighbor (KNN), logistic regression (LR), support vector machines (SVM), and neural network (NN).
The detailed reports that include the classification results for individual instances, the main algorithmic choices and parameter settings for each method, and the elements of confusion matrix are given in Appendix A.

In Table~\ref{tab: chm summary}, columns ``Size" and ``Method" represent the size category and the classification method, respectively. 
Column ``data" distinguishes the results for the training and the test set for each instance.
For the training set, the true class of the vectors are known by construction through the results of Section~\ref{sec:classification}.
For the test set, the true class of the vectors are determined based on the optimal value of the CGLP.
Column ``Accuracy" reports the accuracy of the trained classifier calculated as $\frac{\text{TN} + \text{TP}}{\text{TN} + \text{TP} + \text{FP} + \text{FN}}$, where TN, TP, FP, and FN are the elemtns of the confusion matrix reported in corresponding tables in Appendix A.
These results are obtained by taking the average over the individual instances as shown in Appendix A.
The entry under column ``ML Time" in the row corresponding to the training set contains the average time (in seconds) that it takes to train the classifier.
The entry under column ``ML Time" in the row corresponding to the test set shows the average time it takes for the trained classifier to identify the class of a given vector in the test set.
Column ``CGLP Time" reports the average time that it takes to solve the CGLP for a given objective function vector in the test set.

As observed in Table~\ref{tab: chm summary}, the LR, SVM and NN methods yield high training and test accuracy compared to the RF and KNN.
Considering both accuracy and run time, LR is the superior classifier among others.
Comparing the classification time of the above algorithms with the solution time of the CGLP, we observe the LR, SVM and NN methods respectively achieve an average of $10^4$, $10^2$ and $10^3$ orders of magnitude improvement over the CGLP time for the first problem size category.
This improvement factor becomes even larger as the problem size increases. 
These results show a significant time save when using a classification approach to identify the optimal value of the CGLP with a remarkably high accuracy.
As noted in Section~\ref{sec:introduction}, this approach is most advantageous when there is a large number of objective function vectors that need to be optimized over a fixed cone in the CGLP, which is the case for the CHM problem.
In the \BB settings, for instance, at layer $l$ of the \BB tree, at most $2^l$ objective function vectors (one for each node candidate at that layer) need to be solved through CGLP over the same projection cone of the underlying disjunctions.
Considering the training time for the size category $100 \times 300$, the LR method saves time at any $l \geq 5$.
We will showcase the impact of such time saves inside of \BB in Section~\ref{subsec:IP}.

It is worth noting that the high training accuracy observed in these tables is the result of the artificially-generated training data sets.
Unlike the conventional machine learning settings that could flag such high accuracy as \textit{overfitting}, this occurrence is desirable in function approximation settings.
In particular, in function approximation, there is no randomness in the underlying population. 
The training set is produced artificially with the goal of representing critical points of the domain of the function, which is deterministic.
The test set is the entire function domain, which is the $n$-dimensional space in our model.
Consequently, we seek to find a classifier that achieves high accuracy for the training set as a measure for the quality of the approximation of the underlying function. 
We further note that the higher accuracy level for the test set compared to the training set for some of the instances in Tables~\ref{tab: LR}--\ref{tab: NN} can be attributed to the fact that the training set is designed in such a way that most of its data points are close to the boundary of the polar cone, which makes the classification problem more difficult, leading to a lower training accuracy. The test set, on the other hand, includes any vector in the space whose class can be more accurately identified by the classifier; see Table~\ref{tab: distance} in Appendix B for a sensitivity analysis on the impact of the distance of the test data sets from the boundary of the polar cone on the accuracy of the test set.

\begin{table}[]
	\caption{Summary of classification results for the CHM problem}
	\label{tab: chm summary}
	\begin{center}
				\begin{tabular}{|l|l|l|l|l|l|}
					\hline
					Size                       & Method           & Data   & Accuracy & ML Time    & CGLP Time \\ \hline
					\multirow{10}{*}{100$\times$300}  & \multirow{2}{*}{RF} & train  & 0.99     & 0.34    & -       \\ \cline{3-6} 
					&                    & test  & 0.81     & 1.67e-5 & 0.01    \\ \cline{2-6} 
					& \multirow{2}{*}{KNN} & train & 0.73     & 0.03    & -       \\ \cline{3-6} 
					&                    & test  & 0.6     & 0.002 & 0.01    \\ \cline{2-6} 
					& \multirow{2}{*}{LR} & train & 0.9     & 0.09    & -       \\ \cline{3-6} 
					&                    & test  & 0.99     & 2.00e-6 & 0.01    \\ \cline{2-6} 
					& \multirow{2}{*}{SVM} & train & 0.87     & 0.68    & -       \\ \cline{3-6} 
					&                    & test  & 1     & 0.0003 & 0.01    \\ \cline{2-6} 
					& \multirow{2}{*}{NN} & train & 0.93     & 13.7    & -       \\ \cline{3-6} 
					&                    & test  & 0.98     & 2.6e-6 & 0.01    \\ \hline
					\multirow{10}{*}{200$\times$1000} & \multirow{2}{*}{RF} & train & 0.99     & 1.89    & -       \\ \cline{3-6} 
					&                    & test  & 0.84     & 1.12e-5 & 0.09    \\ \cline{2-6} 
					& \multirow{2}{*}{KNN} & train & 0.77     & 0.45    & -       \\ \cline{3-6} 
					&                    & test  & 0.61     & 0.02 & 0.09    \\ \cline{2-6} 
					& \multirow{2}{*}{LR} & train & 0.95     & 1.07    & -       \\ \cline{3-6} 
					&                    & test  & 0.99     & 4.28e-6    & 0.09    \\ \cline{2-6} 
					& \multirow{2}{*}{SVM} & train & 0.98     & 13.90    & -       \\ \cline{3-6} 
					&                    & test  & 0.99      & 0.002 & 0.09    \\ \cline{2-6} 
					& \multirow{2}{*}{NN} & train & 0.99     & 66.24    & -       \\ \cline{3-6} 
					&                    & test  & 0.99     & 5.13e-6 & 0.09    \\ \hline
					\multirow{10}{*}{500$\times$2000} & \multirow{2}{*}{RF} & train & 0.99     & 4.58    & -       \\ \cline{3-6} 
					&                    & test  & 0.79     & 1.17e-5 & 0.54    \\ \cline{2-6} 
					& \multirow{2}{*}{KNN} & train & 0.92     & 1.83    & -       \\ \cline{3-6} 
					&                    & test  & 0.58     & 0.07 & 0.54    \\ \cline{2-6} 
					& \multirow{2}{*}{LR} & train & 0.99     & 3.67    & -       \\ \cline{3-6} 
					&                    & test  & 0.99     & 5.46e-6 & 0.54    \\ \cline{2-6} 
					& \multirow{2}{*}{SVM} & train & 0.99     & 66.9    & -       \\ \cline{3-6} 
					&                    & test  & 0.96     & 0.006 & 0.54    \\ \cline{2-6} 
					& \multirow{2}{*}{NN} & train  & 0.99     & 111.4    & -       \\ \cline{3-6} 
					&                    & test  & 0.98     & 6.22e-6 & 0.54    \\ \hline
				\end{tabular}
		\end{center}
	\end{table}

	\subsection{0--1 Programs.} \label{subsec:IP}
	
	To assess the potential performance of the machine learning approach in the \BB process, we conduct preliminary computational experiments on benchmark 0--1 programs.
	We emphasize that our goal in this paper is to evaluate the performance of the machine learning approach \textit{in comparison with} that of the CGLP during B\&B.
	As common in such studies in the MIP community, to have a fair and meaningful comparison between these two approaches, we turn off the presolve, heuristics, and cutting plane features in the solver, leading to a controlled \BB algorithm as the basis of our implementation. 
	Based on our experiments, these solver features affect the \BB trees generated from these two approaches differently, creating bias in the end results and invalidating the requirements needed for a controlled test environment. 
	We will compare the result of adding consistency cuts to the base model at the nodes of the \BB tree when the corresponding CGLP is solved exactly by the solver and when it is approximated by the classification method.
	The LP relaxations at each node of the \BB tree are solved by the CPLEX solver at its default settings.
	For the \BB process, we use a fixed branching order based on the order of variable indices $1, \dotsc, n$.
	Considering a fixed variable order is a common practice when studying consistency properties, such as \textit{adaptive consistency}, in constraint programming; see \cite{Hooker1998,BESSIERE200629,Balafrej2016} for a detailed account on this concept.
	For these experiments, we set the time limit to $86400$ seconds (i.e., 24 hours).

	As discussed in Section~\ref{sec:consistency}, it is established in the literature that the consistency cuts can be substantially more effective than the classical cuts such as RLT in reducing the \BB tree size and solution time; \black{see Appendix C for a comparison between the performance of consistency cuts and the strong branching approach.} 
	In the sequel, we evaluate the effectiveness of the classification approach for the consistency framework in further improving its solution time. 
	We conduct our experiments on random instances that are created similarly to the benchmark instances studied in \cite{davarnia:ra:ho:2022}.
	To implement the consistency cuts, we solve the CGLP of Proposition~\ref{prop:CGLP} at each node of the \BB tree.
	For the classification approach, we use the LR method, which showed the best performance among other methods in individual tests as reported in Section~\ref{subsec:chm}.
	At each layer of the \BB tree, we generate the training data set similarly to the method used in Section~\ref{subsec:chm}, as detailed below.
	Let $\mathcal{C} = \{\vc{x} \in \Re^n \, | \, \vc{a}^i \vc{x} \leq 0, \, \forall i \in I \}$ represent the feasible region of the CGLP of Proposition~\ref{prop:new_CGLP}.
	For each constraint $i \in I$, we generate two sets of class-0 vectors according to Corollary~\ref{cor:class-0} by selecting the components of vector $\vc{v}$ randomly from a discrete uniform distribution between $[0, 10000]$ for each $i \in I$, and setting $\bar{\epsilon} = 0.001$.
	To have a balanced training set, we also generate two sets of class-1 vectors. 
	For the first set, we use Corollary~\ref{cor:class1_special_1} to include the negative of each constraint vector with index $k^* \in I$ by setting the weight vector $\vc{w} = \vc{e}^{k^*}$.
	For the second set, we apply the result of Corollary~\ref{cor:class1_special_4}(i) to each constraint index $j^* \in I$ with $w_{m+i} = 1$ for all $i \in N$, $w_{k} = 0$ for all other components, and $\epsilon = 0.001$.
	The number of these class-0 and class-1 vectors is equal to the number of constraints of the CGLP model including the variable bounds.
	For the LR hyperparameters, we set the solver as \textit{liblinear}, the inverse of regularization strength parameter $C$= 0.06, and the threshold probability $P$=0.6.

	The computational results are reported in Table~\ref{tab:greedy_Synthetic_1} for the classical multi-knapsack problem \black{structures} studied in Table~6.1 in \cite{davarnia:ra:ho:2022}.
	These problem instances are of a general 0--1 programming (multi-knapsack) form $\max\{ \vc{c} \vc{\alpha} : A\vc{\alpha} \leq \vc{b}, \vc{\alpha} \in \{0,1\}^n \}$, where the parameters are randomly generated from a discrete uniform distribution on $[-100,100]$. 
	In Table~\ref{tab:greedy_Synthetic_1}, the columns ``$n$" and ``$m$" represent the number of variables and constraints, respectively. 
	For each problem size, we consider five random instances reported in the third column. 
	The columns under ``B\&B" show the number of \BB nodes and the time to solve each problem.
	For these experiments, we use the breadth-first branching strategy that allows for implementing the ML approach over the fixed projection cone in the same layer of the \BB tree.
	We note here that because of this branching strategy, the results obtained for the consistency cuts in this paper are different from those reported in \cite{davarnia:ra:ho:2022} where the default branching strategy of the CPLEX solver is used.
	The symbol ``--" indicates that the problem was not solved to optimality within the time limit.
	
	Columns 6--9 under ``CGLP" show the result of solving the CGLP \eqref{eq:CGLP-1}--\eqref{eq:CGLP-5} to add \black{consistency cuts of rank one at the nodes of the \BB tree; see Appendix D for a computational study of using the ML approach to produce consistency cuts of higher ranks.} 
	The column ``nodes" contains the \BB tree size, and the next column shows the percentage of reduction in the \BB tree size obtained by adding consistency cuts when compared with the base \BB approach.
	The column ``time" reports the total solution time for the consistency approach, and the next column shows the percentage of time reduction achieved by adding consistency cuts when compared with the base \BB approach. 
	These results confirm that the consistency cuts can significantly reduce the \BB tree size and the total solution time when implemented inside of the \BB algorithm.
	
	Columns 10--13 under ``LR" are defined similarly as above for the machine learning approach where the outcome of the LR classifier is used as an approximation for the optimal value of the CGLP.
	Depending on the LR results, clausal inequalities will be added as a consistency cut according to the greedy method discussed in Section~\ref{sec:consistency}; see Remark~\ref{remark:consistency} for more details.
	The time reported in column 12 contains the entire \BB solution time, including the training time for the LR method at each layer of the \BB tree.
	We note that all the instances have been solved to optimality.
	In this table, the bold numbers under a column mark the approach that achieves the best result for that criterion.
	It is evident from these results that the LR method has a significant impact on reducing the total solution time compared to both the base \BB and the consistency approaches. 
	This gap widens as the problem size increases.
	In terms of the \BB tree size, it is clear that the largest reduction is obtained by the consistency approach as it solves the CGLP exactly, whereas the LR method classifies only a subset of all node candidates that admit a consistency cut.
	Despite this trade-off, the LR approach is successful in reducing the \BB tree size substantially compared with the base \BB approach, and within a close interval to that of the consistency approach.

	\begin{table}[!t]
		\caption{Classification results for general multi-knapsack 0--1 programs}
		\label{tab:greedy_Synthetic_1}
		\resizebox{1.00\textwidth}{!}{
			\begin{tabular}{|l|l|l|ll|llll|llll|}
				\hline
				\multirow{2}{*}{$n$}  & \multirow{2}{*}{$m$}  & \multirow{2}{*}{\#} & \multicolumn{2}{c|}{\BB}                  & \multicolumn{4}{c|}{CGLP}                                                                               & \multicolumn{4}{c|}{LR}                                                            \\ \cline{4-13} 
				&                     &                           & \multicolumn{1}{l|}{nodes}   & time                & \multicolumn{1}{l|}{nodes}          & \multicolumn{1}{l|}{$\Delta$ (\%)}       & \multicolumn{1}{l|}{time}     & $\Delta$ (\%)                     & \multicolumn{1}{l|}{nodes}         & \multicolumn{1}{l|}{$\Delta$ (\%)}       & \multicolumn{1}{l|}{time}              & $\Delta$ (\%)                     \\ \hline
				\multirow{5}{*}{40} & \multirow{5}{*}{50} & 1                                & \multicolumn{1}{l|}{34473}   & 3142.60             & \multicolumn{1}{l|}{\textbf{2383}}  & \multicolumn{1}{l|}{\textit{93}} & \multicolumn{1}{l|}{1103.69}  & \textit{64}               & \multicolumn{1}{l|}{3109}          & \multicolumn{1}{l|}{\textit{90}} & \multicolumn{1}{l|}{\textbf{803.37}}   & \textit{74}               \\ \cline{3-13} 
				&                     & 2                                 & \multicolumn{1}{l|}{82417}   & 7496.14             & \multicolumn{1}{l|}{5625}           & \multicolumn{1}{l|}{\textit{93}} & \multicolumn{1}{l|}{3548.07}  & \textit{52}               & \multicolumn{1}{l|}{\textbf{4677}} & \multicolumn{1}{l|}{\textit{94}} & \multicolumn{1}{l|}{\textbf{920.39}}   & \textit{87}               \\ \cline{3-13} 
				&                     & 3                                & \multicolumn{1}{l|}{78965}   & 7167.15             & \multicolumn{1}{l|}{\textbf{5225}}  & \multicolumn{1}{l|}{\textit{93}} & \multicolumn{1}{l|}{2197.50}  & \textit{69}               & \multicolumn{1}{l|}{17587}         & \multicolumn{1}{l|}{\textit{77}} & \multicolumn{1}{l|}{\textbf{2090.72}}  & \textit{70}               \\ \cline{3-13} 
				&                     & 4                                   & \multicolumn{1}{l|}{124327}  & 11377.96            & \multicolumn{1}{l|}{\textbf{5863}}  & \multicolumn{1}{l|}{\textit{95}} & \multicolumn{1}{l|}{2682.24}  & \textit{76}               & \multicolumn{1}{l|}{14113}         & \multicolumn{1}{l|}{\textit{88}} & \multicolumn{1}{l|}{\textbf{1784.86}}  & \textit{84}               \\ \cline{3-13} 
				&                     & 5                            & \multicolumn{1}{l|}{126209}  & 11700.49            & \multicolumn{1}{l|}{\textbf{6197}}  & \multicolumn{1}{l|}{\textit{95}} & \multicolumn{1}{l|}{3120.76}  & \textit{73}               & \multicolumn{1}{l|}{6821}          & \multicolumn{1}{l|}{\textit{94}} & \multicolumn{1}{l|}{\textbf{1126.00}}  & \textit{90}               \\ \hline
				\multirow{5}{*}{45} & \multirow{5}{*}{55} & 1                             & \multicolumn{1}{l|}{1050459} & 119705.93           & \multicolumn{1}{l|}{\textbf{63677}} & \multicolumn{1}{l|}{\textit{93}} & \multicolumn{1}{l|}{31929.71} & \textit{73}               & \multicolumn{1}{l|}{161395}        & \multicolumn{1}{l|}{\textit{84}} & \multicolumn{1}{l|}{\textbf{16446.59}} & \textit{86}               \\ \cline{3-13} 
				&                     & 2                          & \multicolumn{1}{l|}{57621}   & 5449.79             & \multicolumn{1}{l|}{\textbf{3519}}  & \multicolumn{1}{l|}{\textit{93}} & \multicolumn{1}{l|}{2116.51}  & \textit{61}               & \multicolumn{1}{l|}{5907}          & \multicolumn{1}{l|}{\textit{89}} & \multicolumn{1}{l|}{\textbf{1613.00}}  & \textit{70}               \\ \cline{3-13} 
				&                     & 3                    & \multicolumn{1}{l|}{94281}        &  8424.51                   & \multicolumn{1}{l|}{\textbf{4693}}      & \multicolumn{1}{l|}{\textit{95}}   & \multicolumn{1}{l|}{2965.73}         & \textit{64}                 & \multicolumn{1}{l|}{5171}              & \multicolumn{1}{l|}{\textit{94}}   & \multicolumn{1}{l|}{\textbf{1513.17}}         & \textit{82}                 \\ \cline{3-13} 
				&                     & 4                     & \multicolumn{1}{l|}{450175}  & 44544.77            & \multicolumn{1}{l|}{\textbf{22645}} & \multicolumn{1}{l|}{\textit{94}} & \multicolumn{1}{l|}{13808.84} & \textit{69}               & \multicolumn{1}{l|}{31817}         & \multicolumn{1}{l|}{\textit{92}} & \multicolumn{1}{l|}{\textbf{3994.38}}  & \textit{91}               \\ \cline{3-13} 
				&                     & 5                         & \multicolumn{1}{l|}{630449}  & 65030.67            & \multicolumn{1}{l|}{\textbf{44875}} & \multicolumn{1}{l|}{\textit{92}} & \multicolumn{1}{l|}{26270.06} & \textit{59}               & \multicolumn{1}{l|}{252237}        & \multicolumn{1}{l|}{\textit{59}} & \multicolumn{1}{l|}{\textbf{25428.00}} & \textit{60}               \\ \hline
				\multirow{5}{*}{45} & \multirow{5}{*}{60} & 1                     & \multicolumn{1}{l|}{989029}  & 116212.23           & \multicolumn{1}{l|}{\textbf{54721}} & \multicolumn{1}{l|}{\textit{94}} & \multicolumn{1}{l|}{35255.78} & \textit{69}               & \multicolumn{1}{l|}{180123}        & \multicolumn{1}{l|}{\textit{81}} & \multicolumn{1}{l|}{\textbf{18588.75}} & \textit{84}               \\ \cline{3-13} 
				&                     & 2                     & \multicolumn{1}{l|}{257217}  & 25121.01            & \multicolumn{1}{l|}{\textbf{17321}} & \multicolumn{1}{l|}{\textit{93}} & \multicolumn{1}{l|}{9948.30}  & \textit{60}               & \multicolumn{1}{l|}{83747}         & \multicolumn{1}{l|}{\textit{67}} & \multicolumn{1}{l|}{\textbf{8792.25}}  & \textit{65}               \\ \cline{3-13} 
				&                     & 3                        & \multicolumn{1}{l|}{249869}  & 23936.52            & \multicolumn{1}{l|}{\textbf{14063}} & \multicolumn{1}{l|}{\textit{94}} & \multicolumn{1}{l|}{7917.20}  & \textit{66}               & \multicolumn{1}{l|}{35897}         & \multicolumn{1}{l|}{\textit{85}} & \multicolumn{1}{l|}{\textbf{4355.84}}  & \textit{81}               \\ \cline{3-13} 
				&                     & 4                        &  \multicolumn{1}{l|}{-}        &     \textgreater{}86400                & \multicolumn{1}{l|}{\textbf{56847}}      & \multicolumn{1}{l|}{\textit{-}}   & \multicolumn{1}{l|}{31665.33}         & \textit{\textgreater{}63}                 & \multicolumn{1}{l|}{81637}              & \multicolumn{1}{l|}{\textit{-}}   & \multicolumn{1}{l|}{\textbf{29929.28}}         & \textit{\textgreater{}65}                 \\ \cline{3-13} 
				&                     & 5                       & \multicolumn{1}{l|}{320967}  & 31300.88            & \multicolumn{1}{l|}{\textbf{17925}} & \multicolumn{1}{l|}{\textit{94}} & \multicolumn{1}{l|}{12102.74} & \textit{61}               & \multicolumn{1}{l|}{56093}         & \multicolumn{1}{l|}{\textit{82}} & \multicolumn{1}{l|}{\textbf{6487.29}}  & \textit{79}               \\ \hline
				\multirow{5}{*}{50} & \multirow{5}{*}{60} & 1                          & \multicolumn{1}{l|}{329277}  & 31947.91            & \multicolumn{1}{l|}{12393} & \multicolumn{1}{l|}{\textit{96}} & \multicolumn{1}{l|}{11852.70} & \textit{62}               & \multicolumn{1}{l|}{\textbf{12365}}         & \multicolumn{1}{l|}{\textit{96}} & \multicolumn{1}{l|}{\textbf{3025.14}}  & \textit{90}               \\ \cline{3-13} 
				&                     & 2                          & \multicolumn{1}{l|}{531833}  & 54939.26            & \multicolumn{1}{l|}{\textbf{22885}} & \multicolumn{1}{l|}{\textit{95}} & \multicolumn{1}{l|}{20104.97} & \textit{63}               & \multicolumn{1}{l|}{76523}         & \multicolumn{1}{l|}{\textit{85}} & \multicolumn{1}{l|}{\textbf{9224.31}}  & \textit{83}               \\ \cline{3-13} 
				&                     & 3                      & \multicolumn{1}{l|}{-}       & \textgreater{}86400 & \multicolumn{1}{l|}{\textbf{93375}} & \multicolumn{1}{l|}{\textit{-}}  & \multicolumn{1}{l|}{70096.13} & \textit{\textgreater{}18} & \multicolumn{1}{l|}{103369}        & \multicolumn{1}{l|}{\textit{-}}  & \multicolumn{1}{l|}{\textbf{11917.21}} & \textit{\textgreater{}86} \\ \cline{3-13} 
				&                     & 4                           & \multicolumn{1}{l|}{174357}  & 16540.17            & \multicolumn{1}{l|}{\textbf{8533}}  & \multicolumn{1}{l|}{\textit{95}} & \multicolumn{1}{l|}{7478.74}  & \textit{54}               & \multicolumn{1}{l|}{23853}         & \multicolumn{1}{l|}{\textit{86}} & \multicolumn{1}{l|}{\textbf{4164.15}}  & \textit{74}               \\ \cline{3-13} 
				&                     & 5                         & \multicolumn{1}{l|}{-}       & \textgreater{}86400 & \multicolumn{1}{l|}{53671} & \multicolumn{1}{l|}{\textit{-}}  & \multicolumn{1}{l|}{44288.28} & \textit{\textgreater{}48} & \multicolumn{1}{l|}{\textbf{43507}}         & \multicolumn{1}{l|}{\textit{-}}  & \multicolumn{1}{l|}{\textbf{6007.04}}  & \textit{\textgreater{}93} \\ \hline
			\end{tabular}
		}
	\end{table}

	Lastly, we apply our framework to a few benchmark problems from MIPLib \cite{miplib1996} that are studied in \cite{davarnia:ra:ho:2022}.
	We chose pure 0--1 programs that generated search trees that are large enough for a meaningful comparison, but small enough for the algorithms to run without a memory error. 
	The results are reported in Table~\ref{tab:benchmark} where the first column shows the name of the test problem in the MIPLib, and the other columns are defined similarly to those of Table~\ref{tab:greedy_Synthetic_1}.
	It follows from the results of Table~\ref{tab:benchmark} that the general pattern discussed in reference to Table~\ref{tab:greedy_Synthetic_1} holds for these benchmark instances; see Appendix A for more detailed reports that show the classification results for each layer of the \BB tree for some of the instances studied in this section.
	These results also confirm that the machine learning approach tends to be more effective when the \BB tree size is large, as it fails to outperform the consistency approach for the problem $\mathtt{sentoy}$ that has a small \BB tree size.
	Overall, the experiments presented in this section suggest that the machine learning approach is promising in improving both memory- and time-efficiency of adding cutting planes inside of a branch-and-cut framework.

	Because this work introduces a new path toward integrating machine learning into the branch-and-cut process, we believe that the presented framework can be used as a foundation to pursue this direction and expand its scope.
	We hope that the preliminary computational results presented in this paper can motivate and encourage broader and more extensive experiments to unlock the full potential of the machine learning approach in cut-generating efforts when solving challenging MIPs through modern solvers.
	\black{Several directions of future research are of interest, such as applying the ML framework to different problem structures, adapting the approach for different classes of cutting planes, and incorporating the method in combination with modern solver's features such as presolve, cutting planes, and dynamic branching.}

	\begin{table}[!t]
		\caption{Classification results for MIPLIB benchmark instances}
		\label{tab:benchmark}
		\resizebox{1.00\textwidth}{!}{
			\begin{tabular}{|l|l|l|ll|llll|llll|}
				\hline
				\multirow{2}{*}{Class}  & \multirow{2}{*}{$n$}  & \multirow{2}{*}{$m$} & \multicolumn{2}{c|}{\BB}                  & \multicolumn{4}{c|}{CGLP}                                                                               & \multicolumn{4}{c|}{LR}                                                            \\ \cline{4-13} 
				&                     &                           & \multicolumn{1}{l|}{nodes}   & time                & \multicolumn{1}{l|}{nodes}          & \multicolumn{1}{l|}{$\Delta$ (\%)}       & \multicolumn{1}{l|}{time}     & $\Delta$ (\%)                     & \multicolumn{1}{l|}{nodes}         & \multicolumn{1}{l|}{$\Delta$ (\%)}       & \multicolumn{1}{l|}{time}              & $\Delta$ (\%)                     \\ \hline
				$\mathtt{p0033}$ & 33 & 15                                & \multicolumn{1}{l|}{32117}   & 2961.55             & \multicolumn{1}{l|}{\textbf{385}}  & \multicolumn{1}{l|}{\textit{98}} & \multicolumn{1}{l|}{283.14}  & \textit{90}               & \multicolumn{1}{l|}{1927}          & \multicolumn{1}{l|}{\textit{94}} & \multicolumn{1}{l|}{\textbf{205.16}}   & \textit{93}               \\ \hline 
				$\mathtt{pipex}$ & 41 & 48                                & \multicolumn{1}{l|}{1057}   & 1933.48             & \multicolumn{1}{l|}{\textbf{623}}  & \multicolumn{1}{l|}{\textit{41}} & \multicolumn{1}{l|}{1484.10}  & \textit{23}               & \multicolumn{1}{l|}{687}          & \multicolumn{1}{l|}{\textit{35}} & \multicolumn{1}{l|}{\textbf{947.2}}   & \textit{51}               \\ \hline 
				$\mathtt{sentoy}$ & 60 & 30                                & \multicolumn{1}{l|}{703}   & 391.12             & \multicolumn{1}{l|}{\textbf{271}}  & \multicolumn{1}{l|}{\textit{63}} & \multicolumn{1}{l|}{\textbf{187.27}}  & \textit{53}               & \multicolumn{1}{l|}{275}          & \multicolumn{1}{l|}{\textit{61}} & \multicolumn{1}{l|}{207.93}   & \textit{47}               \\ \hline 
				$\mathtt{stein27}$ & 27 & 118                                & \multicolumn{1}{l|}{6099}   & 2474.1             & \multicolumn{1}{l|}{\textbf{4823}}  & \multicolumn{1}{l|}{\textit{22}} & \multicolumn{1}{l|}{2115.07}  & \textit{15}               & \multicolumn{1}{l|}{5245}          & \multicolumn{1}{l|}{\textit{14}} & \multicolumn{1}{l|}{\textbf{1583.46}}   & \textit{36}               \\ \hline 
				$\mathtt{enigma}$ & 100 & 42                                & \multicolumn{1}{l|}{98545}   & 7413.51             & \multicolumn{1}{l|}{\textbf{21679}}  & \multicolumn{1}{l|}{\textit{78}} & \multicolumn{1}{l|}{2669.5}  & \textit{64}               & \multicolumn{1}{l|}{39417}          & \multicolumn{1}{l|}{\textit{60}} & \multicolumn{1}{l|}{\textbf{1260.7}}   & \textit{83}               \\ \hline 
				$\mathtt{lseu}$ & 89 & 28                                & \multicolumn{1}{l|}{219597}   & 4138.8             & \multicolumn{1}{l|}{\textbf{186963}}  & \multicolumn{1}{l|}{\textit{15}} & \multicolumn{1}{l|}{3724.2}  & \textit{10}               & \multicolumn{1}{l|}{191363}          & \multicolumn{1}{l|}{\textit{13}} & \multicolumn{1}{l|}{\textbf{2627.64}}   & \textit{37}               \\ \hline 
			\end{tabular}
		}
	\end{table}

	\section{Conclusion.} \label{sec:conclusion}
	In this paper, we propose a framework that views CGLP as a classification problem and makes use of machine learning methods to approximate the indicator function associated with the CGLP.
	To apply the classification framework, we develop a methodology to generate training data sets that include both class-0 and class-1 vectors.
	These data sets are then used through computational experiments to evaluate the performance of the developed framework on benchmark combinatorial problems and integer programs.
	The results are presented for different machine learning methods, and sensitivity analysis is conducted to identify the marginal impact of some of the data generation factors on the outcome.

	\ACKNOWLEDGMENT{%
		We thank the anonymous referees and the Associate Editor for their helpful comments
		that contributed to improving the paper.
	}
	%
	%
	%
	
	
	

	\bibliographystyle{informs2014}
	\bibliography{consistency}

\begin{thebibliography}{29}
\providecommand{\natexlab}[1]{#1}
\providecommand{\url}[1]{\texttt{#1}}
\providecommand{\urlprefix}{URL }

\bibitem[{Alvarez et~al.(2017)Alvarez, Louveaux, \protect\BIBand{}
  Wehenkel}]{alvarez:lo:we:2017}
Alvarez AM, Louveaux Q, Wehenkel L (2017) A machine learning-based
  approximation of strong branching. \emph{INFORMS Journal on Computing}
  29:185--195.

\bibitem[{Balafrej et~al.(2016)Balafrej, Bessiere, Paparrizou,
  \protect\BIBand{} Trombettoni}]{Balafrej2016}
Balafrej A, Bessiere C, Paparrizou A, Trombettoni G (2016) \emph{Data Mining
  and Constraint Programming: Foundations of a Cross-Disciplinary Approach},
  chapter Adapting Consistency in Constraint Solving, 226--253 (Springer
  International Publishing).

\bibitem[{Balas(1979)}]{Bal79}
Balas E (1979) Disjunctive programming. \emph{Annals of Discrete Mathematics}
  5:3--51.

\bibitem[{Balas(1985)}]{Bal85}
Balas E (1985) Disjunctive programming and a hierarchy of relaxations for
  discrete optimization problems. \emph{SIAM Journal on Algebraic and Discrete
  Methods} 6:466--485.

\bibitem[{Balas \protect\BIBand{} Perregaard(2002)}]{BALAS:per:2002}
Balas E, Perregaard M (2002) Lift-and-project for mixed {0--1} programming:
  {Recent} progress. \emph{Discrete Applied Mathematics} 123(1):129--154.

\bibitem[{Balcan et~al.(2021)Balcan, Prasad, Sandholm, \protect\BIBand{}
  Vitercik}]{balcan:pr:sa:vi:2021}
Balcan MF, Prasad S, Sandholm T, Vitercik E (2021) Sample complexity of tree
  search configuration: Cutting planes and beyond. \emph{Advances in Neural
  Information Processing Systems}, volume~34, 4015--4027.

\bibitem[{Bertsekas(2017)}]{bertsikas:2017}
Bertsekas DP (2017) \emph{Dynamic Programming and Optimal Control, 4th ed.},
  volume 1 and 2 (Nashua, NH: Athena Scientific).

\bibitem[{Bessiere(2006)}]{BESSIERE200629}
Bessiere C (2006) Chapter 3 - constraint propagation. Rossi F, {van Beek} P,
  Walsh T, eds., \emph{Handbook of Constraint Programming}, volume~2 of
  \emph{Foundations of Artificial Intelligence}, 29--83 (Elsevier).

\bibitem[{Busoniu et~al.(2017)Busoniu, Babuska, Schutter, \protect\BIBand{}
  Ernst}]{busoniu:ba:de:er:2017}
Busoniu L, Babuska R, Schutter BD, Ernst D (2017) \emph{Reinforcement Learning
  and Dynamic Programming Using Function Approximators} (CRC Press).

\bibitem[{Conforti et~al.(2014)Conforti, Cornu\'{e}jols, \protect\BIBand{}
  Zambelli}]{conforti:co:za:2014}
Conforti M, Cornu\'{e}jols G, Zambelli G (2014) \emph{Integer Programming}
  (Springer).

\bibitem[{Davarnia(2021)}]{davarnia:2021}
Davarnia D (2021) Strong relaxations for continuous nonlinear programs based on
  decision diagrams. \emph{Operations Research Letters} 49:239--245.

\bibitem[{Davarnia et~al.(2022)Davarnia, Rajabalizadeh, \protect\BIBand{}
  Hooker}]{davarnia:ra:ho:2022}
Davarnia D, Rajabalizadeh A, Hooker J (2022) Achieving consistency with cutting
  planes. \emph{Mathematical Programming}
  \urlprefix\url{https://doi.org/10.1007/s10107-022-01778-8}.

\bibitem[{Davarnia \protect\BIBand{} van Hoeve(2020)}]{davarnia:va:2020}
Davarnia D, van Hoeve WJ (2020) Outer approximation for integer nonlinear
  programs via decision diagrams. \emph{Mathematical Programming} 187:111--150.

\bibitem[{Friedman(1994)}]{friedman:94}
Friedman J (1994) An overview of predictive learning and function
  approximation. Technical report no. 112, Stanford University.

\bibitem[{Hammer \protect\BIBand{} Gersmann(2003)}]{hammer:ge:2003}
Hammer B, Gersmann K (2003) A note on the universal approximation capability of
  support vector machines. \emph{Neural Processing Letters} 17:43--53.

\bibitem[{Hooker(1998)}]{Hooker1998}
Hooker JN (1998) \emph{Advances in Computational and Stochastic Optimization,
  Logic Programming, and Heuristic Search: Interfaces in Computer Science and
  Operations Research}, chapter Constraint Satisfaction Methods for Generating
  Valid Cuts.

\bibitem[{Kalantari \protect\BIBand{} Zhang(2022)}]{kalantari:zh:2022}
Kalantari B, Zhang Y (2022) Algorithm 1024: Spherical triangle algorithm: A
  fast oracle for convex hull membership queries. \emph{ACM Transactions on
  Mathematical Software} 48:1--32.

\bibitem[{Karmarkar(1984)}]{karmarkar:84}
Karmarkar N (1984) A new polynomial-time algorithm for linear programming.
  \emph{In Proceedings of the sixteenth annual ACM symposium on Theory of
  computing} 302--311.

\bibitem[{Khalil et~al.(2016)Khalil, Bodic, Song, Nemhauser, \protect\BIBand{}
  Dilkina}]{khalil:le:so:ne:di:2016}
Khalil E, Bodic PL, Song L, Nemhauser G, Dilkina B (2016) Learning to branch in
  mixed integer programming. \emph{Proceedings of the Thirtieth AAAI Conference
  on Artificial Intelligence}, 724--731.

\bibitem[{{MIPLIB 2}(1996)}]{miplib1996}
{MIPLIB 2} (1996) The mixed integer programming library.
  \urlprefix\url{http://miplib2010.zib.de/miplib2/miplib2.html}.

\bibitem[{Nemhauser \protect\BIBand{} Wolsey(1999)}]{NemWol99}
Nemhauser GL, Wolsey LA (1999) \emph{Integer and Combinatorial Optimization}
  (New York: Wiley).

\bibitem[{Rivlin(1969)}]{rivlin:1969}
Rivlin T (1969) \emph{An Introduction to the Approximation of Functions} (Dover
  Publications).

\bibitem[{Rockafellar(1997)}]{Rockafellar:97}
Rockafellar T (1997) \emph{Convex analysis} (Princeton University Press).

\bibitem[{Salemi \protect\BIBand{} Davarnia(2022)}]{salemi:da:2022}
Salemi H, Davarnia D (2022) On the structure of decision diagram-representable
  mixed integer programs with application to unit commitment. \emph{Operations
  Research} \urlprefix\url{https://doi.org/10.1287/opre.2022.2353}.

\bibitem[{Sherali \protect\BIBand{} Adams(1994)}]{SheAda94}
Sherali HD, Adams WP (1994) A hierarchy of relaxations and convex hull
  characterizations for mixed-integer zero-one programming problems.
  \emph{Discrete Applied Mathematics} 52:83--106.

\bibitem[{Tang et~al.(2020)Tang, Agrawal, \protect\BIBand{}
  Faenza}]{tang:ag:fa:2020}
Tang Y, Agrawal S, Faenza Y (2020) Reinforcement learning for integer
  programming: Learning to cut. \emph{Proceedings of the 37th International
  Conference on Machine Learning, PMLR}, volume 119, 9367--9376.

\bibitem[{Toth et~al.(2017)Toth, O'Rourke, \protect\BIBand{}
  Goodman}]{toth:or:go:2017}
Toth CD, O'Rourke J, Goodman JE (2017) \emph{Handbook of discrete and
  computational geometry} (Chapman and Hall/CRC).

\bibitem[{Wasserman(2006)}]{wasserman:2006}
Wasserman L (2006) \emph{All of nonparametric statistics} (Springer).

\bibitem[{Zarpellon et~al.(2020)Zarpellon, Jo, Lodi, \protect\BIBand{}
  Bengio}]{zarpellon:jo:lo:be:2020}
Zarpellon G, Jo J, Lodi A, Bengio Y (2020) Parameterizing branch-and-bound
  search trees to learn branching policies.
  \urlprefix\url{https://arxiv.org/abs/2002.05120}.

\end{thebibliography}

\newpage




\begin{APPENDICES}
	
	\section{Omitted Computational Results}
	
	In this section, we present the detailed reports for the classification results for the problems studied in Section~\ref{sec:computation}.
	First, we provide the numerical results for the CHM problem.
	In particular, Tables~\ref{tab: random forest}--\ref{tab: NN} report the classification results for the random forest (RF), the $k$-nearest neighbor (KNN), logistic regression (LR), support vector machines (SVM), and neural network (NN), respectively.
	For each of the classification methods, the main algorithmic choices and parameter settings are reported in the caption of the corresponding tables. 
	These parameter values are selected through tuning and cross validation.
	Once selected, these settings remain the same as ``universal settings" across all problem instances and size categories.

	In Tables~\ref{tab: random forest}--\ref{tab: NN}, columns ``Size" and ``\#" represent the size category and the test instance number, respectively.
	Column ``data" distinguishes the results for the training and the test set for each instance.
	The entries of the \textit{confusion matrix} as a result of the training are given in columns ``TN", ``FP", ``FN" and ``TP".
	Column ``Accuracy" reports the accuracy of the trained classifier calculated as $\frac{\text{TN} + \text{TP}}{\text{TN} + \text{TP} + \text{FP} + \text{FN}}$.
	The entry under columns ``ML Time" and ``CGLP Time" are defined similarly to those in Table~\ref{tab: chm summary} in Section~\ref{subsec:chm}.

	Next, we present detailed reports for some instances of the binary programs studied in Section~\ref{subsec:IP} that include the classification results for each layer of the \BB tree.
	Table~\ref{tab: systhetic layer} shows the results for the first instance, i.e., instance \#1 of size $n = 40$ and $m = 50$, of the multi-knapsack problems of Table~\ref{tab:greedy_Synthetic_1}.
	Similarly, Table~\ref{tab: benchmark layer} shows the results for the first instance, i.e., instance $\mathtt{p0033}$, of the benchmark problems of Table~\ref{tab:benchmark}.
	In these tables, each row contains the classification outcome for each layer of the \BB tree.
	The layer number (depth) is reported in the first column of these tables. 
	At each layer, the LR classifier is trained based on the projection cone of the CGLP that is solved at the nodes of that layer. The training time is reported in the second column of Tables~\ref{tab: systhetic layer} and \ref{tab: benchmark layer}. 
	Then, for each unpruned node in that layer, the classification is performed on the vector defined by the fixed partial assignment in that node used in the objective function of the CGLP; see Section~\ref{subsec:IP} for a detailed account on the derivation procedure and algorithmic settings. The classification results are compared with the true values obtained from solving the CGLP at each node. Columns 3-6 in Tables~\ref{tab: systhetic layer} and \ref{tab: benchmark layer} represent the elements of the confusion matrix for all the nodes in the layer. The accuracy of the classification approach for each layer is reported in the last column.

	In view of the results of Tables~\ref{tab: systhetic layer} and \ref{tab: benchmark layer}, we observe that the overall accuracy often decreases as the \BB tree depth increases. Note that size of the the CGLP formed for each node of the \BB tree decreases as the depth of that node increases due to the fact that the number of unfixed variables used in Algorithm~\ref{alg:partialLP-consistency} to multiply with the constraints of the original problem decreases. This leads to a CGLP model with fewer variables and constraints. Therefore, as observed in the CHM results, the classification approach yields a lower accuracy for smaller problem sizes. Nevertheless, the reduction in accuracy is mostly attributed to the increase in the false negative misclassification, as a result of which a cut is not added to the model. This misclassification does not lead to excluding feasible solutions. Furthermore, the above-mentioned decrease in the size of the underlying CGLPs also leads to a smaller training time as the depth of the \BB tree increases. Lastly, we note that the total number of nodes reported in the confusion matrix entries in Tables~\ref{tab: systhetic layer} and \ref{tab: benchmark layer} is smaller than the size of the \BB tree reported in Tables~\ref{tab:greedy_Synthetic_1} and \ref{tab:benchmark} since the LR method is only applied to the nodes that are not pruned by the solver after creation, which are still counted to calculate the total \BB tree size.
	
	\begin{table}[]
		\caption{Classification results for the \textit{Random Forest} with tree number $k=18$}
		\label{tab: random forest}
		\begin{center}
			\resizebox{0.9\textwidth}{!}{
					\begin{tabular}{|l|l|l|l|l|l|l|l|l|l|}
						\hline
						Size                       & \#           & Data  & TN   & FP  & FN  & TP   & Accuracy & ML Time    & CGLP Time \\ \hline
						\multirow{10}{*}{100$\times$300}  & \multirow{2}{*}{1} & train & 800  & 0   & 3   & 797  & 0.99     & 0.36    & -       \\ \cline{3-10} 
						&                    & test  & 200  & 50  & 51  & 199  & 0.79     & 1.99e-5 & 0.01    \\ \cline{2-10} 
						& \multirow{2}{*}{2} & train & 800  & 0   & 8   & 792  & 0.99     & 0.32    & -       \\ \cline{3-10} 
						&                    & test  & 185  & 65  & 59  & 191  & 0.75     & 2.19e-5 & 0.01    \\ \cline{2-10} 
						& \multirow{2}{*}{3} & train & 800  & 0   & 5   & 795  & 0.99     & 0.33    & -       \\ \cline{3-10} 
						&                    & test  & 231  & 19  & 45  & 205  & 0.87     & 1.39e-5 & 0.01    \\ \cline{2-10} 
						& \multirow{2}{*}{4} & train & 799  & 1   & 1   & 799  & 0.99     & 0.34    & -       \\ \cline{3-10} 
						&                    & test  & 209  & 41  & 41  & 209  & 0.83     & 1.39e-5 & 0.01    \\ \cline{2-10} 
						& \multirow{2}{*}{5} & train & 800  & 0   & 1   & 799  & 0.99     & 0.32    & -       \\ \cline{3-10} 
						&                    & test  & 223  & 27  & 56  & 194  & 0.83     & 1.39e-5 & 0.01    \\ \hline
						\multirow{10}{*}{200$\times$1000} & \multirow{2}{*}{1} & train & 2400 & 0   & 9   & 2391 & 0.99     & 1.99    & -       \\ \cline{3-10} 
						&                    & test  & 607  & 93  & 146 & 554  & 0.82     & 1.42e-5 & 0.09    \\ \cline{2-10} 
						& \multirow{2}{*}{2} & train & 2400 & 0   & 11  & 2389 & 0.99     & 2.02    & -       \\ \cline{3-10} 
						&                    & test  & 588  & 112 & 128 & 572  & 0.82     & 1.14e-5 & 0.08    \\ \cline{2-10} 
						& \multirow{2}{*}{3} & train & 2400 & 0   & 6   & 2394 & 0.99     & 1.59    & -       \\ \cline{3-10} 
						&                    & test  & 690  & 10  & 128 & 572  & 0.90     & 1.14e-6    & 0.08    \\ \cline{2-10} 
						& \multirow{2}{*}{4} & train & 2400 & 0   & 8   & 2392 & 0.99     & 1.81    & -       \\ \cline{3-10} 
						&                    & test  & 586  & 114 & 166 & 534  & 0.8      & 7.14e-6 & 0.09    \\ \cline{2-10} 
						& \multirow{2}{*}{5} & train & 2400 & 0   & 4   & 2396 & 0.99     & 2.04    & -       \\ \cline{3-10} 
						&                    & test  & 688  & 12  & 149 & 2396 & 0.88     & 1.21e-5 & 0.09    \\ \hline
						\multirow{10}{*}{500$\times$2000} & \multirow{2}{*}{1} & train & 5000 & 0   & 7   & 4993 & 0.99     & 6.42    & -       \\ \cline{3-10} 
						&                    & test  & 824  & 676 & 266 & 1234 & 0.68     & 1.70e-5 & 1.28    \\ \cline{2-10} 
						& \multirow{2}{*}{2} & train & 5000 & 0   & 12  & 4988 & 0.99     & 7.19    & -       \\ \cline{3-10} 
						&                    & test  & 1298 & 202 & 271 & 1229 & 0.84     & 1.59e-5 & 0.36    \\ \cline{2-10} 
						& \multirow{2}{*}{3} & train & 5000 & 0   & 4   & 4996 & 0.99     & 3.48    & -       \\ \cline{3-10} 
						&                    & test  & 1310 & 190 & 269 & 1231 & 0.84     & 8.00e-6 & 0.36    \\ \cline{2-10} 
						& \multirow{2}{*}{4} & train & 5000 & 0   & 10  & 4990 & 0.99     & 2.93    & -       \\ \cline{3-10} 
						&                    & test  & 1372 & 128 & 284 & 1216 & 0.86     & 7.33e-6 & 0.36    \\ \cline{2-10} 
						& \multirow{2}{*}{5} & train & 5000 & 0   & 8   & 4992 & 0.99     & 2.91    & -       \\ \cline{3-10} 
						&                    & test  & 1080 & 420 & 245 & 1255 & 0.77     & 1.06e-5 & 0.36    \\ \hline
					\end{tabular}
				}
			\end{center}
		\end{table}

		\begin{table}[]
			\caption{Classification results for the \textit{$K$-Nearest Neighbor} with number of neighbors $k=4$}
			\label{tab: knn}
			\begin{center}
				\resizebox{0.9\textwidth}{!}{
						\begin{tabular}{|l|l|l|l|l|l|l|l|l|l|}
							\hline
							Size                       & \#           & Data  & TN   & FP  & FN  & TP   & Accuracy & ML Time    & CGLP Time \\ \hline
							\multirow{10}{*}{100$\times$300}  & \multirow{2}{*}{1} & train & 800  & 0  & 435  & 365  & 0.72     & 0.03  & -       \\ \cline{3-10} 
							&                    & test  & 250  & 0  & 186  & 64   & 0.62     & 0.002 & 0.01    \\ \cline{2-10} 
							& \multirow{2}{*}{2} & train & 800  & 0  & 426  & 374  & 0.73     & 0.03  & -       \\ \cline{3-10} 
							&                    & test  & 243  & 7  & 196  & 54   & 0.59     & 0.002 & 0.01    \\ \cline{2-10} 
							& \multirow{2}{*}{3} & train & 800  & 0  & 428  & 372  & 0.73     & 0.03  & -       \\ \cline{3-10} 
							&                    & test  & 247  & 3  & 207  & 43   & 0.58     & 0.002 & 0.01    \\ \cline{2-10} 
							& \multirow{2}{*}{4} & train & 800  & 0  & 445  & 355  & 0.72     & 0.03  & -       \\ \cline{3-10} 
							&                    & test  & 249  & 1  & 197  & 53   & 0.60     & 0.002 & 0.01    \\ \cline{2-10} 
							& \multirow{2}{*}{5} & train & 800  & 0  & 419  & 381  & 0.73     & 0.03  & -       \\ \cline{3-10} 
							&                    & test  & 250  & 0  & 193  & 57   & 0.61     & 0.002 & 0.01    \\ \hline
							\multirow{10}{*}{200$\times$1000} & \multirow{2}{*}{1} & train & 2400 & 0  & 966  & 1434 & 0.79     & 0.45  & -       \\ \cline{3-10} 
							&                    & test  & 700  & 0  & 553  & 147  & 0.60     & 0.03  & 0.09    \\ \cline{2-10} 
							& \multirow{2}{*}{2} & train & 2400 & 0  & 1002 & 1398 & 0.79     & 0.46  & -       \\ \cline{3-10} 
							&                    & test  & 699  & 1  & 531  & 169  & 0.62     & 0.03  & 0.08    \\ \cline{2-10} 
							& \multirow{2}{*}{3} & train & 2400 & 0  & 994  & 1406 & 0.79     & 0.46  & -       \\ \cline{3-10} 
							&                    & test  & 676  & 24 & 539  & 161  & 0.59     & 0.02  & 0.08    \\ \cline{2-10} 
							& \multirow{2}{*}{4} & train & 2400 & 0  & 1014 & 1386 & 0.78     & 0.41  & -       \\ \cline{3-10} 
							&                    & test  & 667  & 33 & 555  & 145  & 0.58     & 0.02  & 0.09    \\ \cline{2-10} 
							& \multirow{2}{*}{5} & train & 2400 & 0  & 1004 & 1396 & 0.81     & 0.42  & -       \\ \cline{3-10} 
							&                    & test  & 699  & 1  & 543  & 157  & 0.61     & 0.02  & 0.09    \\ \hline
							\multirow{10}{*}{500$\times$2000} & \multirow{2}{*}{1} & train & 5000 & 0  & 796  & 4204 & 0.92     & 2.84  & -       \\ \cline{3-10} 
							&                    & test  & 1471 & 29 & 1191 & 309  & 0.59     & 0.11  & 1.28    \\ \cline{2-10} 
							& \multirow{2}{*}{2} & train & 5000 & 0  & 802  & 4198 & 0.91     & 2.96  & -       \\ \cline{3-10} 
							&                    & test  & 1448 & 52 & 1183 & 317  & 0.58     & 0.06  & 0.36    \\ \cline{2-10} 
							& \multirow{2}{*}{3} & train & 5000 & 0  & 826  & 4174 & 0.91     & 1.41  & -       \\ \cline{3-10} 
							&                    & test  & 1496 & 4  & 1198 & 302  & 0.59     & 0.06  & 0.36    \\ \cline{2-10} 
							& \multirow{2}{*}{4} & train & 5000 & 0  & 789  & 4211 & 0.92     & 1.12  & -       \\ \cline{3-10} 
							&                    & test  & 1388 & 12 & 1173 & 327  & 0.60     & 0.06  & 0.36    \\ \cline{2-10} 
							& \multirow{2}{*}{5} & train & 5000 & 0  & 799  & 4201 & 0.92     & 1.11  & -       \\ \cline{3-10} 
							&                    & test  & 1407 & 93 & 1202 & 298  & 0.56     & 0.05  & 0.36    \\ \hline
						\end{tabular}
					}
				\end{center}
			\end{table}
			
			\begin{table}[]
				\caption{Classification results for the \textit{Logistic Regression} with \textit{solver}$=$\textit{liblinear}, inverse of regularization strength $C=1$, and probability $P=0.5$}
				\label{tab: LR}
				\begin{center}
					\resizebox{0.9\textwidth}{!}{
							\begin{tabular}{|l|l|l|l|l|l|l|l|l|l|}
								\hline
								Size                       & \#           & Data  & TN   & FP  & FN  & TP   & Accuracy & ML Time    & CGLP Time \\ \hline
								\multirow{10}{*}{100$\times$300}  & \multirow{2}{*}{1} & train & 651  & 149 & 1  & 799  & 0.90     & 0.09    & -       \\ \cline{3-10} 
								&                    & test  & 250  & 0   & 1  & 249  & 0.99     & 0.0     & 0.01    \\ \cline{2-10} 
								& \multirow{2}{*}{2} & train & 642  & 158 & 0  & 800  & 0.90     & 0.09    & -       \\ \cline{3-10} 
								&                    & test  & 250  & 0   & 1  & 249  & 0.99     & 1.99e-6 & 0.01    \\ \cline{2-10} 
								& \multirow{2}{*}{3} & train & 635  & 165 & 0  & 800  & 0.89     & 0.09    & -       \\ \cline{3-10} 
								&                    & test  & 250  & 0   & 0  & 250  & 1        & 2.00e-6 & 0.01    \\ \cline{2-10} 
								& \multirow{2}{*}{4} & train & 644  & 156 & 1  & 799  & 0.90     & 0.09    & -       \\ \cline{3-10} 
								&                    & test  & 250  & 0   & 0  & 250  & 1        & 4.00e-6 & 0.01    \\ \cline{2-10} 
								& \multirow{2}{*}{5} & train & 648  & 152 & 1  & 799  & 0.90     & 0.10    & -       \\ \cline{3-10} 
								&                    & test  & 350  & 0   & 1  & 249  & 0.99     & 2.00e-6 & 0.01    \\ \hline
								\multirow{10}{*}{200$\times$1000} & \multirow{2}{*}{1} & train & 2212 & 188 & 0  & 2400 & 0.96     & 1.16    & -       \\ \cline{3-10} 
								&                    & test  & 700  & 0   & 0  & 700  & 1        & 4.99e-6 & 0.09    \\ \cline{2-10} 
								& \multirow{2}{*}{2} & train & 2189 & 211 & 0  & 2400 & 0.95     & 1.22    & -       \\ \cline{3-10} 
								&                    & test  & 700  & 0   & 2  & 698  & 0.99     & 2.85e-6 & 0.08    \\ \cline{2-10} 
								& \multirow{2}{*}{3} & train & 2186 & 214 & 0  & 2400 & 0.95     & 1.07    & -       \\ \cline{3-10} 
								&                    & test  & 700  & 0   & 1  & 699  & 0.99     & 6.42e-6 & 0.08    \\ \cline{2-10} 
								& \multirow{2}{*}{4} & train & 2173 & 227 & 0  & 2400 & 0.95     & 0.96    & -       \\ \cline{3-10} 
								&                    & test  & 700  & 0   & 0  & 700  & 1        & 3.57e-6 & 0.09    \\ \cline{2-10} 
								& \multirow{2}{*}{5} & train & 2195 & 205 & 0  & 2400 & 0.95     & 0.96    & -       \\ \cline{3-10} 
								&                    & test  & 700  & 0   & 0  & 700  & 1        & 3.57e-6 & 0.09    \\ \hline
								\multirow{10}{*}{500$\times$2000} & \multirow{2}{*}{1} & train & 4942 & 58  & 0  & 5000 & 0.99     & 5.30    & -       \\ \cline{3-10} 
								&                    & test  & 1500 & 0   & 4  & 1496 & 0.99     & 6.33e-6 & 1.28    \\ \cline{2-10} 
								& \multirow{2}{*}{2} & train & 4960 & 40  & 0  & 5000 & 0.99     & 3.62    & -       \\ \cline{3-10} 
								&                    & test  & 1500 & 0   & 3  & 1497 & 0.99     & 7.66e-6 & 0.36    \\ \cline{2-10} 
								& \multirow{2}{*}{3} & train & 4948 & 52  & 0  & 5000 & 0.99     & 3.12    & -       \\ \cline{3-10} 
								&                    & test  & 1500 & 0   & 6  & 1494 & 0.99     & 4.99e-6 & 0.36    \\ \cline{2-10} 
								& \multirow{2}{*}{4} & train & 4975 & 25  & 0  & 5000 & 0.99     & 3.31    & -       \\ \cline{3-10} 
								&                    & test  & 1500 & 0   & 9  & 1491 & 0.99     & 4.33e-6 & 0.36    \\ \cline{2-10} 
								& \multirow{2}{*}{5} & train & 4924 & 76  & 0  & 5000 & 0.99     & 3.03    & -       \\ \cline{3-10} 
								&                    & test  & 1500 & 0   & 9  & 1491 & 0.99     & 3.99e-6 & 0.36    \\ \hline
							\end{tabular}
						}
					\end{center}
				\end{table}

				\begin{table}[]
					\caption{Classification results for the \textit{Support Vector Machine} with \textit{kernel}$=$\textit{linear} and regularization parameter $C=0.4$}
					\label{tab: SVM}
					\begin{center}
						\resizebox{0.9\textwidth}{!}{
								\begin{tabular}{|l|l|l|l|l|l|l|l|l|l|}
									\hline
									Size                       & \#           & Data  & TN   & FP  & FN  & TP   & Accuracy & ML Time    & CGLP Time \\ \hline
									\multirow{10}{*}{100$\times$300}  & \multirow{2}{*}{1} & train & 602  & 198 & 0   & 800  & 0.87     & 0.70   & -       \\ \cline{3-10} 
									&                    & test  & 250  & 0   & 0   & 250  & 1        & 0.0003 & 0.01    \\ \cline{2-10} 
									& \multirow{2}{*}{2} & train & 600  & 200 & 0   & 800  & 0.87     & 0.71   & -       \\ \cline{3-10} 
									&                    & test  & 250  & 0   & 0   & 250  & 1        & 0.0003 & 0.01    \\ \cline{2-10} 
									& \multirow{2}{*}{3} & train & 602  & 198 & 0   & 800  & 0.87     & 0.67   & -       \\ \cline{3-10} 
									&                    & test  & 250  & 0   & 0   & 250  & 1        & 0.0003 & 0.01    \\ \cline{2-10} 
									& \multirow{2}{*}{4} & train & 600  & 200 & 0   & 800  & 0.87     & 0.66   & -       \\ \cline{3-10} 
									&                    & test  & 250  & 0   & 0   & 250  & 1        & 0.0003 & 0.01    \\ \cline{2-10} 
									& \multirow{2}{*}{5} & train & 606  & 194 & 0   & 800  & 0.87     & 0.66   & -       \\ \cline{3-10} 
									&                    & test  & 250  & 0   & 0   & 250  & 1        & 0.0003 & 0.01    \\ \hline
									\multirow{10}{*}{200$\times$1000} & \multirow{2}{*}{1} & train & 2346 & 54  & 0   & 2400 & 0.98     & 15.53  & -       \\ \cline{3-10} 
									&                    & test  & 700  & 0   & 5   & 695  & 0.99     & 0.003  & 0.09    \\ \cline{2-10} 
									& \multirow{2}{*}{2} & train & 2342 & 58  & 0   & 2400 & 0.98     & 13.91  & -       \\ \cline{3-10} 
									&                    & test  & 700  & 0   & 9   & 691  & 0.99     & 0.002  & 0.08    \\ \cline{2-10} 
									& \multirow{2}{*}{3} & train & 2317 & 83  & 0   & 2400 & 0.98     & 13.96  & -       \\ \cline{3-10} 
									&                    & test  & 700  & 0   & 13  & 687  & 0.99     & 0.002  & 0.08    \\ \cline{2-10} 
									& \multirow{2}{*}{4} & train & 2304 & 96  & 0   & 2400 & 0.98     & 13.24  & -       \\ \cline{3-10} 
									&                    & test  & 700  & 0   & 4   & 696  & 0.99     & 0.002  & 0.09    \\ \cline{2-10} 
									& \multirow{2}{*}{5} & train & 2317 & 83  & 0   & 2400 & 0.98     & 12.90  & -       \\ \cline{3-10} 
									&                    & test  & 700  & 0   & 6   & 694  & 0.99     & 0.002  & 0.09    \\ \hline
									\multirow{10}{*}{500$\times$2000} & \multirow{2}{*}{1} & train & 5000 & 0   & 2   & 4998 & 0.99     & 100.63 & -       \\ \cline{3-10} 
									&                    & test  & 1500 & 0   & 84  & 1416 & 0.97     & 0.009  & 1.28    \\ \cline{2-10} 
									& \multirow{2}{*}{2} & train & 5000 & 0   & 5   & 4995 & 0.99     & 73.96  & -       \\ \cline{3-10} 
									&                    & test  & 1500 & 0   & 109 & 1391 & 0.96     & 0.005  & 0.36    \\ \cline{2-10} 
									& \multirow{2}{*}{3} & train & 5000 & 0   & 2   & 4998 & 0.99     & 51.58  & -       \\ \cline{3-10} 
									&                    & test  & 1500 & 0   & 95  & 1405 & 0.96     & 0.004  & 0.36    \\ \cline{2-10} 
									& \multirow{2}{*}{4} & train & 5000 & 0   & 6   & 4994 & 0.99     & 59.99  & -       \\ \cline{3-10} 
									&                    & test  & 1500 & 0   & 118 & 1382 & 0.96     & 0.004  & 0.36    \\ \cline{2-10} 
									& \multirow{2}{*}{5} & train & 4998 & 2   & 3   & 4997 & 0.99     & 52.90  & -       \\ \cline{3-10} 
									&                    & test  & 1500 & 0   & 75  & 1425 & 0.97     & 0.005  & 0.36    \\ \hline
								\end{tabular}
							}
						\end{center}
					\end{table}

					\begin{table}[]
						\caption{Classification results for the \textit{Neural Network} with \textit{hidden layer} $=18$, \textit{activation function}$=$ \textit{identity}, \textit{solver} $=$ \textit{sgd}, \textit{batch size} $=100$, \textit{regularization parameter} $\alpha=0.01$, \textit{learning rate} $=$ \textit{adaptive}, \textit{initial learning rate} $= 0.0007$, and \textit{maximum number of iterations} $=1500$}
						\label{tab: NN}
						\begin{center}
							\resizebox{0.9\textwidth}{!}{
									\begin{tabular}{|l|l|l|l|l|l|l|l|l|l|}
										\hline
										Size                       & \#           & Data  & TN   & FP  & FN  & TP   & Accuracy & ML Time    & CGLP Time \\ \hline
										\multirow{10}{*}{100$\times$300}              & \multirow{2}{*}{1}             & train             & 696              & 104             & 6              & 794              & 0.93                 & 12.88               & -                  \\ \cline{3-10} 
										&                                & test              & 250              & 0               & 6              & 244              & 0.98                 & 4.00e-6             & 0.01               \\ \cline{2-10} 
										& \multirow{2}{*}{2}             & train             & 712              & 88              & 2              & 798              & 0.94                 & 14.36               & -                  \\ \cline{3-10} 
										&                                & test              & 250              & 0               & 5              & 245              & 0.99                 & 4.00e-6             & 0.01               \\ \cline{2-10} 
										& \multirow{2}{*}{3}             & train             & 712              & 88              & 3              & 797              & 0.94                 & 14.37               & -                  \\ \cline{3-10} 
										&                                & test              & 250              & 0               & 6              & 244              & 0.98                 & 1.99e-6             & 0.01               \\ \cline{2-10} 
										& \multirow{2}{*}{4}             & train             & 704              & 96              & 9              & 791              & 0.93                 & 13.25               & -                  \\ \cline{3-10} 
										&                                & test              & 250              & 0               & 8              & 242              & 0.98                 & 2.00e-6             & 0.01               \\ \cline{2-10} 
										& \multirow{2}{*}{5}             & train             & 708              & 92              & 8              & 792              & 0.93                 & 13.71               & -                  \\ \cline{3-10} 
										&                                & test              & 250              & 0               & 8              & 242              & 0.98                 & 1.99e-6             & 0.01               \\ \hline
										\multirow{10}{*}{200$\times$1000}             & \multirow{2}{*}{1}             & train             & 2380             & 20              & 1              & 2399             & 0.99                 & 65.60               & -                  \\ \cline{3-10} 
										&                                & test              & 700              & 0               & 7              & 693              & 0.99                 & 6.42e-6             & 0.09               \\ \cline{2-10} 
										& \multirow{2}{*}{2}             & train             & 2383             & 17              & 2              & 2398             & 0.99                 & 69.05               & -                  \\ \cline{3-10} 
										&                                & test              & 700              & 0               & 16             & 684              & 0.98                 & 4.28e-6             & 0.08               \\ \cline{2-10} 
										& \multirow{2}{*}{3}             & train             & 2376             & 24              & 1              & 2399             & 0.99                 & 62.59               & -                  \\ \cline{3-10} 
										&                                & test              & 700              & 0               & 19             & 681              & 0.98                 & 3.57e-6             & 0.08               \\ \cline{2-10} 
										& \multirow{2}{*}{4}             & train             & 2367             & 33              & 1              & 2399             & 0.99                 & 65.45               & -                  \\ \cline{3-10} 
										&                                & test              & 700              & 0               & 14             & 686              & 0.99                 & 6.42e-6             & 0.09               \\ \cline{2-10} 
										& \multirow{2}{*}{5}             & train             & 2367             & 33              & 4              & 2396             & 0.99                 & 63.31               & -                  \\ \cline{3-10} 
										&                                & test              & 700              & 0               & 9              & 691              & 0.99                 & 5.00e-6             & 0.09               \\ \hline
										\multirow{10}{*}{500$\times$2000}             & \multirow{2}{*}{1}             & train             & 4994             & 6               & 1              & 4999             & 0.99                 & 158.87              & -                  \\ \cline{3-10} 
										&                                & test              & 1500             & 0               & 25             & 1475             & 0.99                 & 7.66e-6             & 1.28               \\ \cline{2-10} 
										& \multirow{2}{*}{2}             & train             & 4998             & 2               & 2              & 4998             & 0.99                 & 111.11              & -                  \\ \cline{3-10} 
										&                                & test              & 1500             & 0               & 31             & 1469             & 0.98                 & 8.33e-6             & 0.36               \\ \cline{2-10} 
										& \multirow{2}{*}{3}             & train             & 4995             & 5               & 2              & 4998             & 0.99                 & 110.06              & -                  \\ \cline{3-10} 
										&                                & test              & 1500             & 0               & 24             & 1476             & 0.99                 & 5.33e-6             & 0.36               \\ \cline{2-10} 
										& \multirow{2}{*}{4}             & train             & 4999             & 1               & 3              & 4997             & 0.99                 & 86.85               & -                  \\ \cline{3-10} 
										&                                & test              & 1500             & 0               & 31             & 1469             & 0.98                 & 4.66e-6             & 0.36               \\ \cline{2-10} 
										& \multirow{2}{*}{5}             & train             & 4987             & 13              & 2              & 4998             & 0.99                 & 92.90               & -                  \\ \cline{3-10} 
										&                                & test              & 1500             & 0               & 36             & 1464             & 0.98                 & 5.33e-6             & 0.36               \\ \hline
									\end{tabular}
								}
							\end{center}
						\end{table}

						\begin{table}[]
							\caption{Classification results for each layer of the \BB tree for the first instance of Table~\ref{tab:greedy_Synthetic_1}}
							\label{tab: systhetic layer}
							\begin{center}
								\resizebox{0.9\textwidth}{!}{
									\begin{tabular}{|l|l|l|l|l|l|l|}
										\hline
										\BB Tree Depth & Train Time & True Negative & False Positive & False Negative & True Positive & Accuracy \\ \hline
										1     & 20.83       & 2              & 0               & 0               & 0              & 1.00     \\ \hline
										2     & 20.42       & 2              & 0               & 0               & 2              & 1.00     \\ \hline
										3     & 19.42       & 2              & 0               & 0               & 2              & 1.00     \\ \hline
										4     & 18.19       & 2              & 0               & 0               & 2              & 1.00     \\ \hline
										5     & 17.75       & 4              & 0               & 0               & 0              & 1.00     \\ \hline
										6     & 16.72       & 7              & 0               & 0               & 1              & 1.00     \\ \hline
										7     & 15.54       & 9              & 1               & 0               & 4              & 0.93     \\ \hline
										8     & 15.10       & 10             & 1               & 1               & 6              & 0.89     \\ \hline
										9     & 13.53       & 11             & 0               & 1               & 10             & 0.95     \\ \hline
										10    & 12.76       & 20             & 0               & 2               & 2              & 0.92     \\ \hline
										11    & 12.53       & 20             & 0               & 3               & 16             & 0.92     \\ \hline
										12    & 11.63       & 20             & 2               & 3               & 17             & 0.88     \\ \hline
										13    & 10.70       & 22             & 0               & 4               & 17             & 0.91     \\ \hline
										14    & 10.44       & 22             & 0               & 4               & 19             & 0.91     \\ \hline
										15    & 9.64        & 23             & 2               & 5               & 18             & 0.85     \\ \hline
										16    & 8.67        & 30             & 0               & 9               & 13             & 0.83     \\ \hline
										17    & 8.09        & 26             & 1               & 14              & 35             & 0.80     \\ \hline
										18    & 7.34        & 34             & 0               & 23              & 22             & 0.71     \\ \hline
										19    & 6.65        & 33             & 0               & 24              & 21             & 0.69     \\ \hline
										20    & 5.99        & 23             & 1               & 47              & 30             & 0.52     \\ \hline
										21    & 5.45        & 38             & 0               & 31              & 12             & 0.62     \\ \hline
										22    & 4.58        & 52             & 0               & 33              & 15             & 0.67     \\ \hline
										23    & 3.40        & 44             & 0               & 54              & 36             & 0.60     \\ \hline
										24    & 2.97        & 38             & 0               & 62              & 10             & 0.44     \\ \hline
										25    & 2.58        & 46             & 0               & 51              & 17             & 0.55     \\ \hline
										26    & 2.27        & 44             & 1               & 46              & 18             & 0.57     \\ \hline
										27    & 1.93        & 53             & 0               & 42              & 0              & 0.56     \\ \hline
										28    & 1.65        & 43             & 0               & 54              & 0              & 0.44     \\ \hline
										29    & 1.40        & 23             & 0               & 56              & 0              & 0.29     \\ \hline
										30    & 1.15        & 41             & 0               & 35              & 0              & 0.54     \\ \hline
										31    & 0.95        & 36             & 0               & 29              & 0              & 0.55     \\ \hline
										32    & 0.80        & 33             & 0               & 22              & 0              & 0.60     \\ \hline
										33    & 0.61        & 17             & 0               & 15              & 0              & 0.53     \\ \hline
										34    & 0.40        & 17             & 0               & 9               & 0              & 0.65     \\ \hline
										35    & 0.33        & 7              & 0               & 8               & 0              & 0.47     \\ \hline
										36    & 0.18        & 4              & 0               & 2               & 0              & 0.67     \\ \hline
										37    & 0.10        & 4              & 0               & 1               & 0              & 0.80     \\ \hline
										38    & 0.04        & 1              & 0               & 0               & 0              & 1.00     \\ \hline
									\end{tabular}
								}
							\end{center}
						\end{table}

						\begin{table}[]
							\caption{Classification results for each layer of the \BB tree for the first instance of Table~\ref{tab:benchmark}}
							\label{tab: benchmark layer}
							\begin{center}
								\resizebox{0.9\textwidth}{!}{
									\begin{tabular}{|l|l|l|l|l|l|l|}
										\hline
										\BB Tree Depth & Train Time & True Negative & False Positive & False Negative & True Positive & Accuracy \\ \hline
										1     & 6.0372      & 2              & 0               & 0               & 0              & 1.00     \\ \hline
										2     & 5.5961      & 3              & 0               & 0               & 1              & 1.00     \\ \hline
										3     & 5.1293      & 3              & 0               & 1               & 2              & 0.83     \\ \hline
										4     & 4.9456      & 4              & 1               & 0               & 3              & 0.88     \\ \hline
										5     & 4.7717      & 5              & 0               & 1               & 2              & 0.88     \\ \hline
										6     & 4.5564      & 7              & 0               & 1               & 1              & 0.89     \\ \hline
										7     & 4.0524      & 6              & 1               & 2               & 4              & 0.77     \\ \hline
										8     & 3.8808      & 7              & 0               & 4               & 3              & 0.71     \\ \hline
										9     & 3.7526      & 10             & 0               & 3               & 7              & 0.85     \\ \hline
										10    & 3.3934      & 12             & 2               & 5               & 6              & 0.72     \\ \hline
										11    & 2.9123      & 14             & 0               & 8               & 8              & 0.73     \\ \hline
										12    & 2.6759      & 17             & 1               & 12              & 6              & 0.64     \\ \hline
										13    & 2.4469      & 12             & 2               & 13              & 10             & 0.59     \\ \hline
										14    & 2.2926      & 19             & 0               & 13              & 4              & 0.64     \\ \hline
										15    & 2.0172      & 23             & 0               & 22              & 13             & 0.62     \\ \hline
										16    & 1.8169      & 31             & 0               & 23              & 21             & 0.69     \\ \hline
										17    & 1.6957      & 29             & 1               & 36              & 24             & 0.59     \\ \hline
										18    & 1.6244      & 32             & 0               & 41              & 20             & 0.56     \\ \hline
										19    & 1.4567      & 35             & 0               & 44              & 12             & 0.52     \\ \hline
										20    & 1.2023      & 31             & 0               & 53              & 20             & 0.48     \\ \hline
										21    & 0.7969      & 26             & 0               & 40              & 14             & 0.50     \\ \hline
										22    & 0.6471      & 24             & 0               & 38              & 10             & 0.47     \\ \hline
										23    & 0.5648      & 28             & 0               & 31              & 7              & 0.53     \\ \hline
										24    & 0.4168      & 26             & 0               & 26              & 8              & 0.57     \\ \hline
										25    & 0.3536      & 23             & 0               & 22              & 4              & 0.55     \\ \hline
										26    & 0.2855      & 20             & 0               & 15              & 0              & 0.57     \\ \hline
										27    & 0.1453      & 17             & 0               & 8               & 0              & 0.68     \\ \hline
										28    & 0.1172      & 14             & 0               & 5               & 0              & 0.74     \\ \hline
										29    & 0.0871      & 7              & 0               & 2               & 0              & 0.78     \\ \hline
										30    & 0.0694      & 5              & 0               & 1               & 0              & 0.83     \\ \hline
									\end{tabular}
								}
							\end{center}
						\end{table}

						\clearpage
						\newpage
						
						\section{Sensitivity Analysis}
						
						As noted in Section \eqref{subsec:chm}, one set of class-1 vectors in the training data set is generated based on the result of Corollary~\ref{cor:class1_special_4}(i), which requires non-redundancy of constraints of $\mathcal{C}$. 
						Since checking redundancy of constraints in practical applications is time-consuming, we perform sensitivity analysis on the redundancy ratio of constraints to evaluate its impact on classification results.
						For this analysis, we study different redundancy ratios as in Figure 22 of \cite{kalantari:zh:2022} by creating redundant vectors through a conic combination of vectors of $\mathcal{C}$ and adding them to the cone.
						We conduct these experiments on the first instance of each size category.
						These results for the top three classifiers LR, SVM and NN are summarized in Figure~\ref{fig:redundancy}.
						As shown in these figures, the classification accuracy for training data declines slightly as the redundancy ratio increases.
						This decline is due to the fact that the training data generated through Corollary~\ref{cor:class1_special_4}(i) are not guaranteed to be a class-1 vector when the base constraint is redundant. 
						These vectors, however, are labeled as class-1 in the training set, leading to a potential misclassification.
						Nevertheless, it is inferred from these experiments that the impact of this source of misclassification is reduced due to the balance in the remaining sets of training data, which results in maintaining high accuracy of the test set across different redundancy ratios and problem sizes.

						Another deciding factor in creating training data set is $\epsilon$, which determines the proximity of class-0 and class-1 vectors to the boundary of the polar cone.
						To investigate the marginal impact of different values of $\epsilon$ on the classification performance, we next conduct sensitivity analysis on the first instance of each problem size category and each classification method.
						These results for the top three classifiers LR, SVM and NN are summarized in Figure~\ref{fig:epsilon}.
						As observed in these figures, the common trend across all classification methods and problem sizes is that, increasing the $\epsilon$ value leads to a higher training accuracy and a lower test accuracy.
						This trend is expected as the larger $\epsilon$ values broaden the spatial gap between class-0 and class-1 training data, which enables the classifier to improve training accuracy at the expense of reducing test accuracy.

						\begin{figure}[!hb]
							\centering
							\includegraphics[scale=0.6,clip,trim={0.25cm 0cm 0.2cm 0.2cm}]{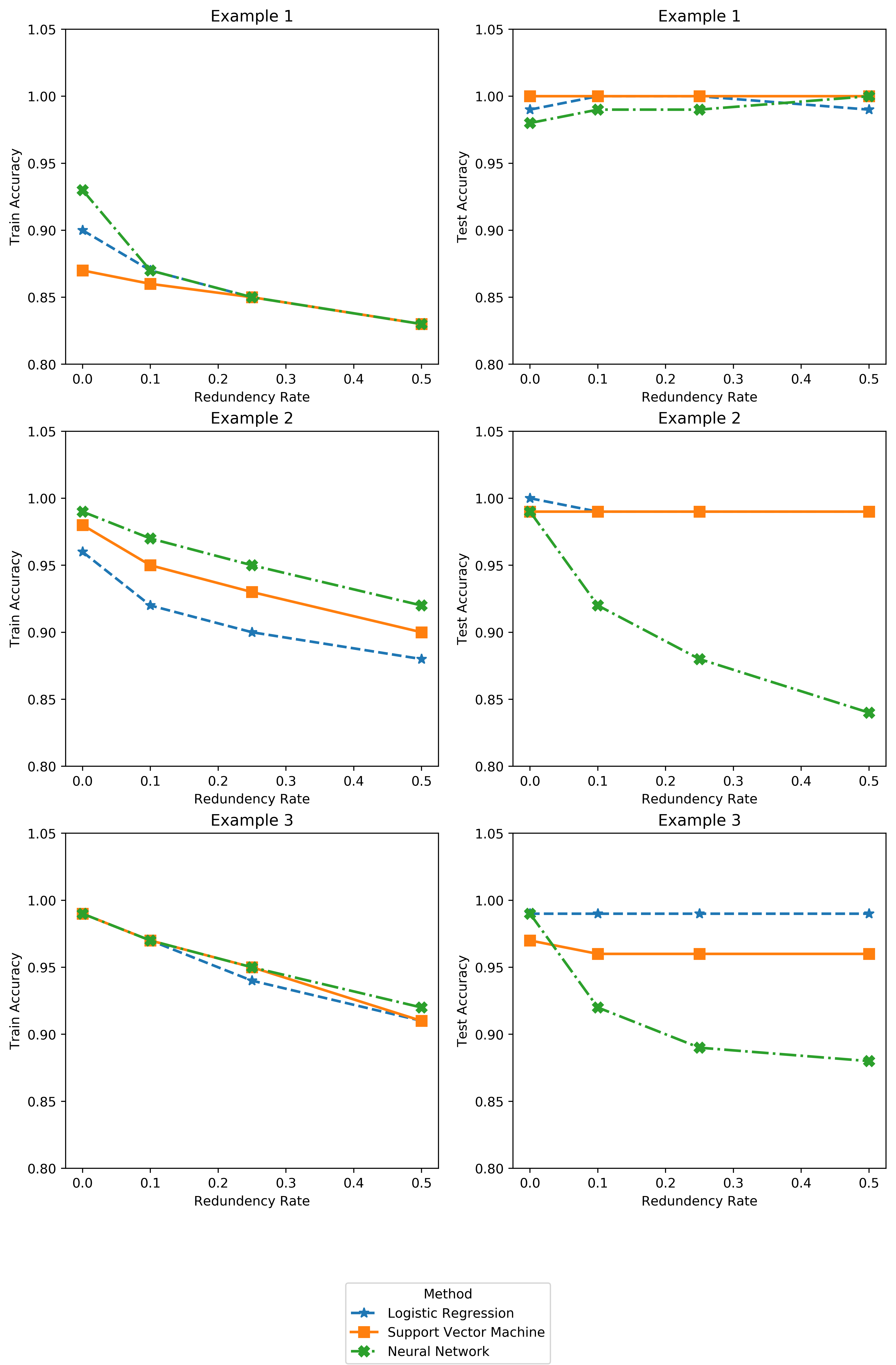}
							\caption{Sensitivity analysis results for redundancy ratio}
							\label{fig:redundancy}
						\end{figure}

						\begin{figure}[!hb]
							\centering
							\includegraphics[scale=0.65,clip,trim={0.25cm 0cm 0.2cm 0.2cm}]{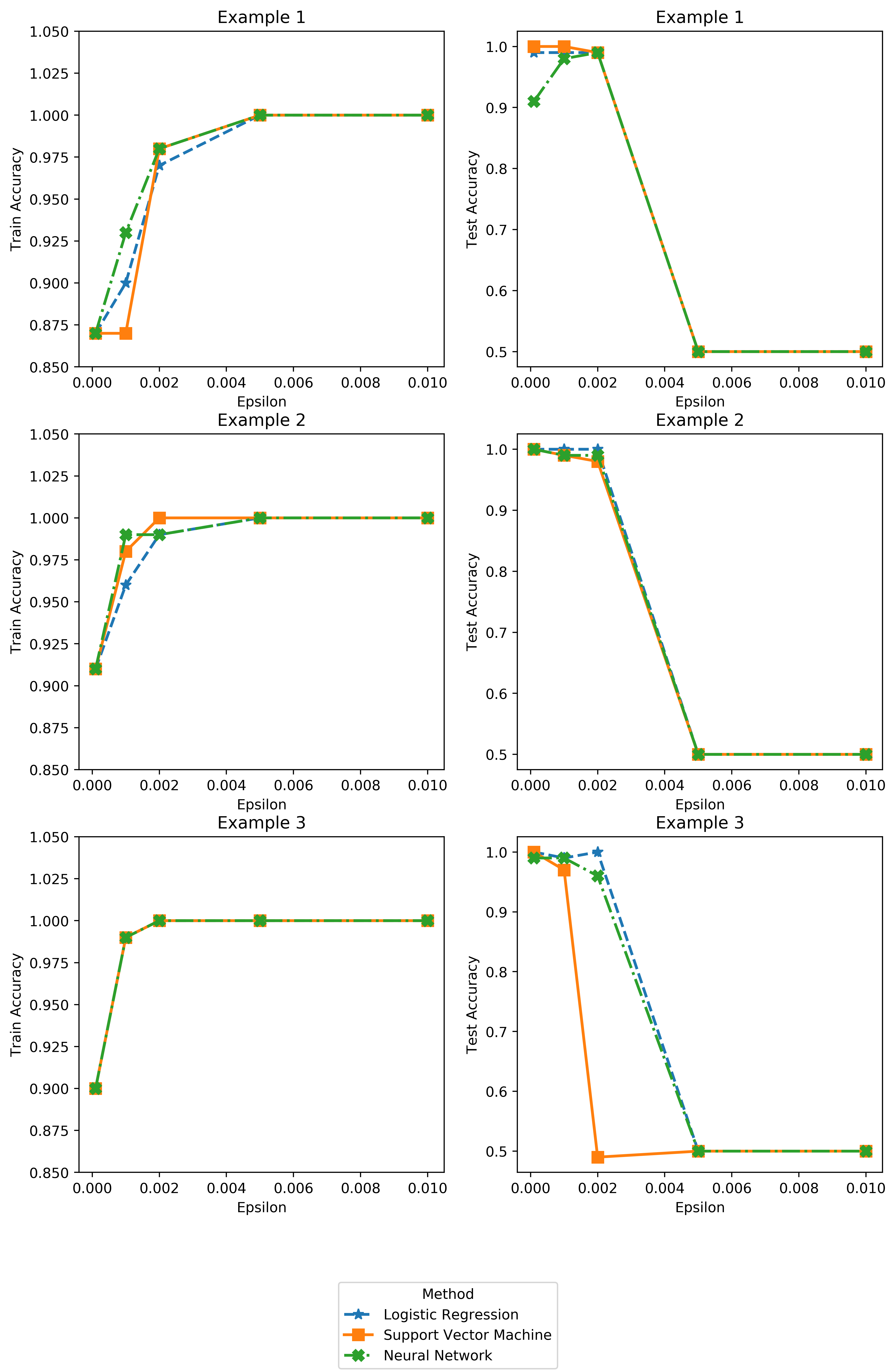}
							\caption{Sensitivity analysis results for values of $\epsilon$}
							\label{fig:epsilon}
						\end{figure}

						As noted in Section 4.1, we have used $4(m+n)$ data points in our training set for the CHM problem. In Figures~\ref{fig:T1} and \ref{fig:T2}, we show the relationship between the size of the training set and the training time and accuracy of the classification for the LR method. We conduct these experiments on one instance of each size category $(100\times300)$, $(200\times1000)$, and $(500\times2000)$. We consider four different sizes for the training set as multiples of the total number of constraints in the cone $\mathcal{C}$, namely $2(m+n)$, $4(m+n)$, $6(m+n)$, and $8(m+n)$.
						All the other parameters are set similarly to those used in the experiments reported in Section 4.1. As observed in these figures, increasing the size of the training set leads to increasing the accuracy at the price of increasing the training time. Based on these results, the choice of $4(m+n)$ provides a reasonable trade-off between the training time and accuracy. 
						
						\begin{figure}[!hb]
							\centering
							\includegraphics[scale=0.75]{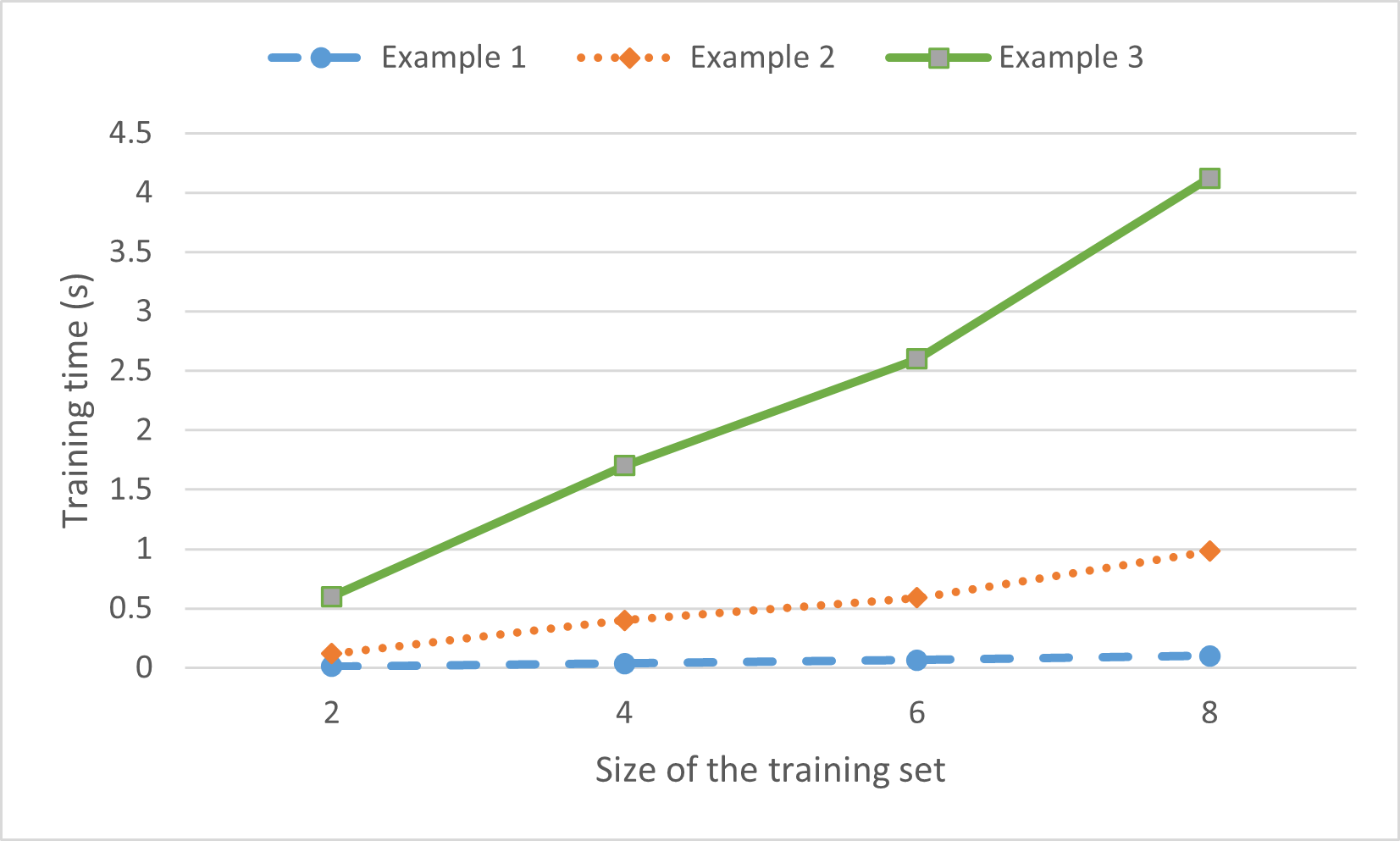}
							\caption{The impact of the size of the training set on the training time for the LR method in CHM problem.}
							\label{fig:T1}
						\end{figure}
						
						\begin{figure}[!hb]
							\centering
							\includegraphics[scale=0.75]{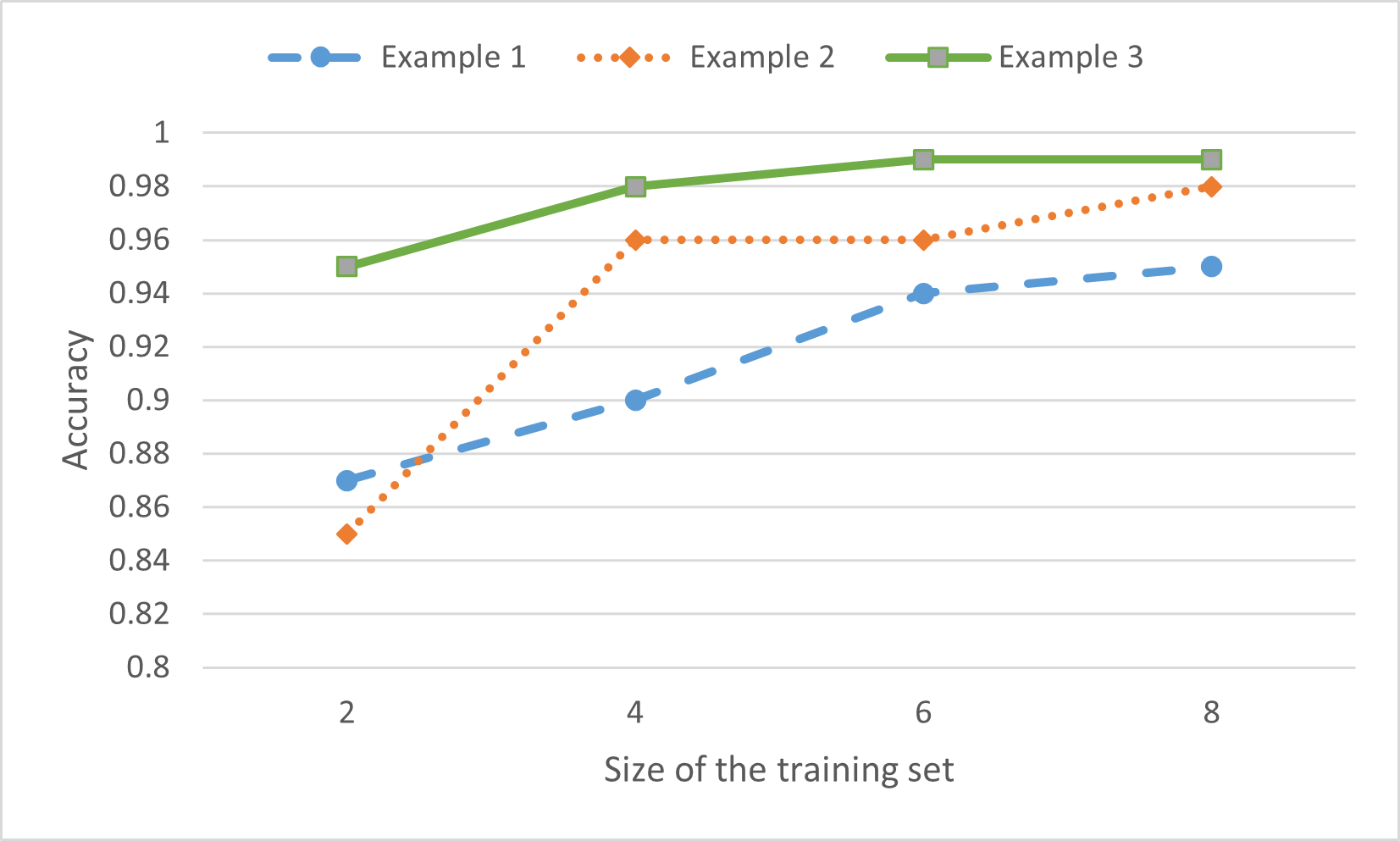}
							\caption{The impact of the size of the training set on the training accuracy for the LR method in CHM problem.}
							\label{fig:T2}
						\end{figure}

						For the next sensitivity analysis, we consider the test accuracy for the LR method applied to the CHM problem for different categories of test instances based on their distance from the boundary of the polar cone. To this end, we use similar settings to those mentioned in Section 4.1 to train a LR classifier for one instance of each size category. Then, we use the resulting classifier to predict the class of data points in the test sets that are generated based on a similar method used to generate the training set where the proximity of the generated vector is determined by the $\epsilon$ value in Corollary 1 for class-0 and Corollary 5(i) for class-1. As noted in Section 4.1, the vectors produced in this process are normalized to have a uniform scaling for all data points in class-1 and class-0. We consider five categories for $\epsilon \in \{0.0005, 0.001, 0.0015, 0.002, 0.0025\}$ to represent different distance values. For each distance category of the problem size $(100\times300)$, $(200\times1000)$, and $(500\times2000)$, we consider the test set size 1600, 4800, and 10000, respectively. The accuracy results, including the elements of the confusion matrix, are given in Table~\ref{tab: distance} and Figure~\ref{fig:E1}. As observed from these plots, the test accuracy increases with the distance of the data points from the boundary of the polar cone.

					\begin{table}[]
						\caption{Classification results for the \textit{Logistic Regression} approach for test instances with different distance values from the boundary of the polar cone}
						\label{tab: distance}
						\begin{center}
									\begin{tabular}{|l|l|l|l|l|l|l|}
										\hline
										Size                       & $\epsilon$           & TN   & FP  & FN  & TP   & Accuracy  \\ \hline
										\multirow{5}{*}{100$\times$300} & 0.0005 & 774  & 26 & 190  & 610  & 0.87 \\ \cline{2-7}  
										& 0.001 & 793  & 7 & 149  & 651  & 0.91       \\ \cline{2-7}  
										& 0.0015 & 800  & 0 & 94  & 706  & 0.95       \\ \cline{2-7} 
										& 0.002 & 800  & 0 & 61  & 739  & 0.98       \\ \cline{2-7} 
										& 0.0025 & 800  & 0 & 15  & 785  & 0.99       \\ \hline
										\multirow{5}{*}{200$\times$1000} & 0.0005 & 2359 & 41 & 394  & 2006 & 0.91     \\ \cline{2-7} 
										&  0.001  & 2385  & 15   & 191  & 2209  & 0.94    \\ \cline{2-7} 
										&  0.0015  & 2400  & 0   & 117  & 2283  & 0.98    \\ \cline{2-7} 
										&  0.002  & 2400  & 0   & 46  & 2354  & 0.99    \\ \cline{2-7} 
										&  0.0025  & 2400  & 0   & 0  & 2400  & 1    \\ \hline 
										\multirow{5}{*}{500$\times$2000} & 0.0005 & 5000 & 0  & 855 & 4145 & 0.92  \\ \cline{2-7} 
										&  0.001  & 5000 & 0   & 504  & 4496 & 0.96    \\ \cline{2-7} 
										&  0.0015  & 5000 & 0   & 57  & 4943 & 0.99    \\ \cline{2-7} 
										&  0.002  & 5000 & 0   & 0  & 5000 & 1    \\ \cline{2-7} 
										&  0.0025  & 5000 & 0   & 0  & 5000 & 1    \\ \hline
									\end{tabular}
							\end{center}
						\end{table}

						\begin{figure}[!hb]
							\centering
							\includegraphics[scale=0.75]{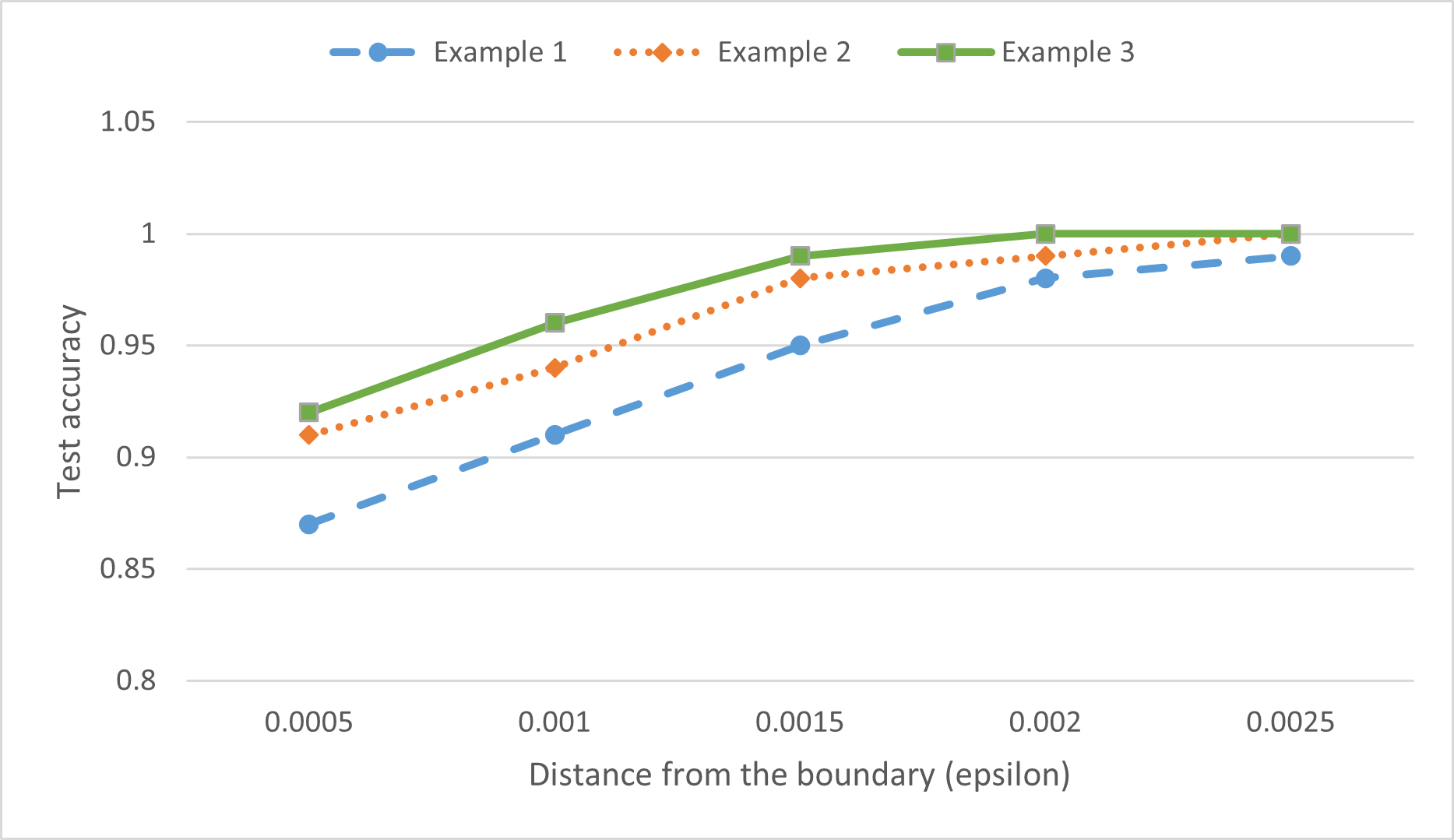}
							\caption{The test accuracy for different distance categories for the LR method in CHM problem.}
							\label{fig:E1}
						\end{figure}

						\section{Consistency Cuts and Strong Branching}
						As mentioned in Section~\ref{sec:consistency}, the goal of using consistency cuts is to achieve (partial) LP-consistency.
						An alternative approach to ensure partial-LP consistency for rank $r$ at any node of the \BB tree associated with the partial assignment $\vc{\alpha}_I = \vc{v}_I$ such that $\S_{\LP} \cap \{\vc{\alpha} \in \Re^n: \vc{\alpha}_I = \vc{v}_I\} \neq \emptyset$ is to check whether there exists $J\subseteq N \setminus I$ with $|J|= r$ such that the LP relaxations described by $\S_{\LP} \cap \{\vc{\alpha} \in \Re^n: \vc{\alpha}_{I\cup J} = \vc{v}_{I \cup J}\}$ is infeasible for every 0--1 value assignment $\vc{\alpha}_J = \vc{v}_J$.
						In this case, the node will be pruned.
						This approach is similar to the so-called \textit{strong branching} technique frequently used in \BB by solvers.
						It is shown in \cite{davarnia:ra:ho:2022} that the consistency cut framework has a computational advantage compared to the strong branching method as outlined next.
						In particular, Proposition 5.5 in \cite{davarnia:ra:ho:2022} shows that applying the CGLP of Proposition~\ref{prop:CGLP} for a certain rank $r$ can achieve LP-consistency of ranks higher than $r$. 
						In other words, the CGLP obtained from the intersection of the multilinear constraints produced by the multiplication with unfixed variables can lead to separating inconsistent faces that remain feasible if we consider the disjunctive model for each unfixed variable individually (which is equivalent to the outcome of the strong branching); see Example 5.3 in \cite{davarnia:ra:ho:2022} for an illustration.
						This property can lead to a significant reduction in the \BB tree size when using the CGLP approach compared to the strong branching approach.
						As a numerical evidence for the above property, we present in Tables~\ref{tab:strong_branching} and \ref{tab:strong_branching_benchmark} computational experiments on the instances of Tables~\ref{tab:greedy_Synthetic_1} and \ref{tab:benchmark}, respectively, which show that the outcome of the consistency approach outperforms that of the strong branching method in terms of both the tree size and the solution time.

						\begin{table}[!t]
							\caption{\black{Comparison of the consistency cuts and the strong branching approach for multi-knapsack problems}}
							\label{tab:strong_branching}
							\resizebox{1.00\textwidth}{!}{
								\begin{tabular}{|l|l|l|ll|llll|llll|}
									\hline
									\multirow{2}{*}{$n$}  & \multirow{2}{*}{$m$}  & \multirow{2}{*}{\#} & \multicolumn{2}{c|}{\BB}                  & \multicolumn{4}{c|}{CGLP}                                                                               & \multicolumn{4}{c|}{Strong Branching}                                                            \\ \cline{4-13} 
									&                     &                           & \multicolumn{1}{l|}{nodes}   & time                & \multicolumn{1}{l|}{nodes}          & \multicolumn{1}{l|}{$\Delta$ (\%)}       & \multicolumn{1}{l|}{time}     & $\Delta$ (\%)                     & \multicolumn{1}{l|}{nodes}         & \multicolumn{1}{l|}{$\Delta$ (\%)}       & \multicolumn{1}{l|}{time}              & $\Delta$ (\%)                     \\ \hline
									\multirow{5}{*}{40} & \multirow{5}{*}{50} & 1                                & \multicolumn{1}{l|}{34473}   & 3142.60             & \multicolumn{1}{l|}{\textbf{2383}}  & \multicolumn{1}{l|}{\textit{93}} & \multicolumn{1}{l|}{\textbf{1103.69}}  & \textit{64}               & \multicolumn{1}{l|}{12755}          & \multicolumn{1}{l|}{\textit{63}} & \multicolumn{1}{l|}{1351.06}   & \textit{57}               \\ \cline{3-13} 
									&                     & 2                                 & \multicolumn{1}{l|}{82417}   & 7496.14             & \multicolumn{1}{l|}{\textbf{5625}}           & \multicolumn{1}{l|}{\textit{93}} & \multicolumn{1}{l|}{\textbf{3548.07}}  & \textit{52}               & \multicolumn{1}{l|}{26373} & \multicolumn{1}{l|}{\textit{68}} & \multicolumn{1}{l|}{4122.8}   & \textit{45}               \\ \cline{3-13} 
									&                     & 3                                & \multicolumn{1}{l|}{78965}   & 7167.15             & \multicolumn{1}{l|}{\textbf{5225}}  & \multicolumn{1}{l|}{\textit{93}} & \multicolumn{1}{l|}{\textbf{2197.50}}  & \textit{69}               & \multicolumn{1}{l|}{25269}         & \multicolumn{1}{l|}{\textit{68}} & \multicolumn{1}{l|}{2436.78}  & \textit{66}               \\ \cline{3-13} 
									&                     & 4                                   & \multicolumn{1}{l|}{124327}  & 11377.96            & \multicolumn{1}{l|}{\textbf{5863}}  & \multicolumn{1}{l|}{\textit{95}} & \multicolumn{1}{l|}{\textbf{2682.24}}  & \textit{76}               & \multicolumn{1}{l|}{36055}         & \multicolumn{1}{l|}{\textit{71}} & \multicolumn{1}{l|}{3754.41}  & \textit{67}               \\ \cline{3-13} 
									&                     & 5                            & \multicolumn{1}{l|}{126209}  & 11700.49            & \multicolumn{1}{l|}{\textbf{6197}}  & \multicolumn{1}{l|}{\textit{95}} & \multicolumn{1}{l|}{\textbf{3120.76}}  & \textit{73}               & \multicolumn{1}{l|}{37863}          & \multicolumn{1}{l|}{\textit{70}} & \multicolumn{1}{l|}{3393.05}  & \textit{71}               \\ \hline
									\multirow{5}{*}{45} & \multirow{5}{*}{55} & 1                             & \multicolumn{1}{l|}{1050459} & 119705.93           & \multicolumn{1}{l|}{\textbf{63677}} & \multicolumn{1}{l|}{\textit{93}} & \multicolumn{1}{l|}{\textbf{31929.71}} & \textit{73}               & \multicolumn{1}{l|}{189083}        & \multicolumn{1}{l|}{\textit{82}} & \multicolumn{1}{l|}{49079.07} & \textit{59}               \\ \cline{3-13} 
									&                     & 2                          & \multicolumn{1}{l|}{57621}   & 5449.79             & \multicolumn{1}{l|}{\textbf{3519}}  & \multicolumn{1}{l|}{\textit{93}} & \multicolumn{1}{l|}{\textbf{2116.51}}  & \textit{61}               & \multicolumn{1}{l|}{17863}          & \multicolumn{1}{l|}{\textit{69}} & \multicolumn{1}{l|}{2288.58}  & \textit{58}               \\ \cline{3-13} 
									&                     & 3                    & \multicolumn{1}{l|}{94281}        &  8424.51                   & \multicolumn{1}{l|}{\textbf{4693}}      & \multicolumn{1}{l|}{\textit{95}}   & \multicolumn{1}{l|}{\textbf{2965.73}}         & \textit{64}                 & \multicolumn{1}{l|}{26399}              & \multicolumn{1}{l|}{\textit{72}}   & \multicolumn{1}{l|}{3875.14}         & \textit{54}                 \\ \cline{3-13} 
									&                     & 4                     & \multicolumn{1}{l|}{450175}  & 44544.77            & \multicolumn{1}{l|}{\textbf{22645}} & \multicolumn{1}{l|}{\textit{94}} & \multicolumn{1}{l|}{\textbf{13808.84}} & \textit{69}               & \multicolumn{1}{l|}{112543}         & \multicolumn{1}{l|}{\textit{75}} & \multicolumn{1}{l|}{18708.34}  & \textit{58}               \\ \cline{3-13} 
									&                     & 5                         & \multicolumn{1}{l|}{630449}  & 65030.67            & \multicolumn{1}{l|}{\textbf{44875}} & \multicolumn{1}{l|}{\textit{92}} & \multicolumn{1}{l|}{\textbf{26270.06}} & \textit{59}               & \multicolumn{1}{l|}{176525}        & \multicolumn{1}{l|}{\textit{71}} & \multicolumn{1}{l|}{27962.9} & \textit{57}               \\ \hline
									\multirow{5}{*}{45} & \multirow{5}{*}{60} & 1                     & \multicolumn{1}{l|}{989029}  & 116212.23           & \multicolumn{1}{l|}{\textbf{54721}} & \multicolumn{1}{l|}{\textit{94}} & \multicolumn{1}{l|}{\textbf{35255.78}} & \textit{69}               & \multicolumn{1}{l|}{326379}        & \multicolumn{1}{l|}{\textit{67}} & \multicolumn{1}{l|}{48809.14} & \textit{58}               \\ \cline{3-13} 
									&                     & 2                     & \multicolumn{1}{l|}{257217}  & 25121.01            & \multicolumn{1}{l|}{\textbf{17321}} & \multicolumn{1}{l|}{\textit{93}} & \multicolumn{1}{l|}{\textbf{9948.30}}  & \textit{60}               & \multicolumn{1}{l|}{56587}         & \multicolumn{1}{l|}{\textit{78}} & \multicolumn{1}{l|}{11304.45}  & \textit{55}               \\ \cline{3-13} 
									&                     & 3                        & \multicolumn{1}{l|}{249869}  & 23936.52            & \multicolumn{1}{l|}{\textbf{14063}} & \multicolumn{1}{l|}{\textit{94}} & \multicolumn{1}{l|}{\textbf{7917.20}}  & \textit{66}               & \multicolumn{1}{l|}{59969}         & \multicolumn{1}{l|}{\textit{76}} & \multicolumn{1}{l|}{9335.40}  & \textit{61}               \\ \cline{3-13} 
									&                     & 4                        &  \multicolumn{1}{l|}{-}        &     \textgreater{}86400                & \multicolumn{1}{l|}{\textbf{56847}}      & \multicolumn{1}{l|}{\textit{-}}   & \multicolumn{1}{l|}{\textbf{31665.33}}         & \textit{\textgreater{}63}                 & \multicolumn{1}{l|}{239325}              & \multicolumn{1}{l|}{\textit{-}}   & \multicolumn{1}{l|}{38016.84}         & \textit{\textgreater{}56}                 \\ \cline{3-13} 
									&                     & 5                       & \multicolumn{1}{l|}{320967}  & 31300.88            & \multicolumn{1}{l|}{\textbf{17925}} & \multicolumn{1}{l|}{\textit{94}} & \multicolumn{1}{l|}{\textbf{12102.74}} & \textit{61}               & \multicolumn{1}{l|}{80241}         & \multicolumn{1}{l|}{\textit{75}} & \multicolumn{1}{l|}{13772.13}  & \textit{56}               \\ \hline
									\multirow{5}{*}{50} & \multirow{5}{*}{60} & 1                          & \multicolumn{1}{l|}{329277}  & 31947.91            & \multicolumn{1}{l|}{\textbf{12393}} & \multicolumn{1}{l|}{\textit{96}} & \multicolumn{1}{l|}{\textbf{11852.70}} & \textit{62}               & \multicolumn{1}{l|}{8945}         & \multicolumn{1}{l|}{\textit{72}} & \multicolumn{1}{l|}{14696.62}  & \textit{54}               \\ \cline{3-13} 
									&                     & 2                          & \multicolumn{1}{l|}{531833}  & 54939.26            & \multicolumn{1}{l|}{\textbf{22885}} & \multicolumn{1}{l|}{\textit{95}} & \multicolumn{1}{l|}{\textbf{20104.97}} & \textit{63}               & \multicolumn{1}{l|}{106367}         & \multicolumn{1}{l|}{\textit{80}} & \multicolumn{1}{l|}{21975.36}  & \textit{60}               \\ \cline{3-13} 
									&                     & 3                      & \multicolumn{1}{l|}{-}       & \textgreater{}86400 & \multicolumn{1}{l|}{\textbf{93375}} & \multicolumn{1}{l|}{\textit{-}}  & \multicolumn{1}{l|}{\textbf{70096.13}} & \textit{\textgreater{}18} & \multicolumn{1}{l|}{103369}        & \multicolumn{1}{l|}{\textit{-}}  & \multicolumn{1}{l|}{76896.19} & \textit{\textgreater{}11} \\ \cline{3-13} 
									&                     & 4                           & \multicolumn{1}{l|}{174357}  & 16540.17            & \multicolumn{1}{l|}{\textbf{8533}}  & \multicolumn{1}{l|}{\textit{95}} & \multicolumn{1}{l|}{\textbf{7478.74}}  & \textit{54}               & \multicolumn{1}{l|}{41845}         & \multicolumn{1}{l|}{\textit{76}} & \multicolumn{1}{l|}{9262.45}  & \textit{44}               \\ \cline{3-13} 
									&                     & 5                         & \multicolumn{1}{l|}{-}       & \textgreater{}86400 & \multicolumn{1}{l|}{\textbf{53671}} & \multicolumn{1}{l|}{\textit{-}}  & \multicolumn{1}{l|}{\textbf{44288.28}} & \textit{\textgreater{}48} & \multicolumn{1}{l|}{43507}         & \multicolumn{1}{l|}{\textit{-}}  & \multicolumn{1}{l|}{50112.72}  & \textit{\textgreater{}42} \\ \hline
								\end{tabular}
							}
						\end{table}

						\begin{table}[!t]
							\caption{\black{Comparison of the consistency cuts and the strong branching approach for MIPLIB problems}}
							\label{tab:strong_branching_benchmark}
							\resizebox{1.00\textwidth}{!}{
								\begin{tabular}{|l|l|l|ll|llll|llll|}
									\hline
									\multirow{2}{*}{Class}  & \multirow{2}{*}{$n$}  & \multirow{2}{*}{$m$} & \multicolumn{2}{c|}{\BB}                  & \multicolumn{4}{c|}{CGLP}                                                                               & \multicolumn{4}{c|}{Strong Branching}                                                            \\ \cline{4-13} 
									&                     &                           & \multicolumn{1}{l|}{nodes}   & time                & \multicolumn{1}{l|}{nodes}          & \multicolumn{1}{l|}{$\Delta$ (\%)}       & \multicolumn{1}{l|}{time}     & $\Delta$ (\%)                     & \multicolumn{1}{l|}{nodes}         & \multicolumn{1}{l|}{$\Delta$ (\%)}       & \multicolumn{1}{l|}{time}              & $\Delta$ (\%)                     \\ \hline
									$\mathtt{p0033}$ & 33 & 15                                & \multicolumn{1}{l|}{32117}   & 2961.55             & \multicolumn{1}{l|}{\textbf{385}}  & \multicolumn{1}{l|}{\textit{98}} & \multicolumn{1}{l|}{\textbf{283.14}}  & \textit{90}               & \multicolumn{1}{l|}{5139}          & \multicolumn{1}{l|}{\textit{84}} & \multicolumn{1}{l|}{887.3}   & \textit{70}               \\ \hline 
									$\mathtt{pipex}$ & 41 & 48                                & \multicolumn{1}{l|}{1057}   & 1933.48             & \multicolumn{1}{l|}{\textbf{623}}  & \multicolumn{1}{l|}{\textit{41}} & \multicolumn{1}{l|}{\textbf{1484.10}}  & \textit{23}               & \multicolumn{1}{l|}{707}          & \multicolumn{1}{l|}{\textit{33}} & \multicolumn{1}{l|}{1527.07}   & \textit{19}               \\ \hline 
									$\mathtt{sentoy}$ & 60 & 30                                & \multicolumn{1}{l|}{703}   & 391.12             & \multicolumn{1}{l|}{\textbf{271}}  & \multicolumn{1}{l|}{\textit{63}} & \multicolumn{1}{l|}{\textbf{187.27}}  & \textit{53}               & \multicolumn{1}{l|}{365}          & \multicolumn{1}{l|}{\textit{48}} & \multicolumn{1}{l|}{218.96}   & \textit{44}               \\ \hline 
									$\mathtt{stein27}$ & 27 & 118                                & \multicolumn{1}{l|}{6099}   & 2474.1             & \multicolumn{1}{l|}{\textbf{4823}}  & \multicolumn{1}{l|}{\textit{22}} & \multicolumn{1}{l|}{2115.07}  & \textit{15}               & \multicolumn{1}{l|}{4819}          & \multicolumn{1}{l|}{\textit{21}} & \multicolumn{1}{l|}{\textbf{1927.38}}   & \textit{22}               \\ \hline 
									$\mathtt{enigma}$ & 100 & 42                                & \multicolumn{1}{l|}{98545}   & 7413.51             & \multicolumn{1}{l|}{\textbf{21679}}  & \multicolumn{1}{l|}{\textit{78}} & \multicolumn{1}{l|}{\textbf{2669.5}}  & \textit{64}               & \multicolumn{1}{l|}{55185}          & \multicolumn{1}{l|}{\textit{44}} & \multicolumn{1}{l|}{3632.37}   & \textit{51}               \\ \hline 
									$\mathtt{lseu}$ & 89 & 28                                & \multicolumn{1}{l|}{219597}   & 4138.8             & \multicolumn{1}{l|}{\textbf{186963}}  & \multicolumn{1}{l|}{\textit{15}} & \multicolumn{1}{l|}{\textbf{3724.2}}  & \textit{10}               & \multicolumn{1}{l|}{193245}          & \multicolumn{1}{l|}{\textit{12}} & \multicolumn{1}{l|}{3806.96}   & \textit{8}               \\ \hline 
								\end{tabular}
							}
						\end{table}

						\section{\black{Machine Learning Approach for Higher Consistency Ranks}}
						
						In this section, we present computational results to evaluate the performance of the ML approach in approximating the CGLP corresponding to partial LP-consistency of higher ranks compared to that of rank one presented in Table~\ref{tab:greedy_Synthetic_1} of  Section~\ref{subsec:IP}.
						To obtain the CGLP for rank $r$, following Algorithm~\ref{alg:partialLP-consistency}, we need to multiply each constraint of the original model with all $2^r$ product combinations of every tuple of size $r$ of unfixed variables at each node of the \BB tree. 
						As a result, increasing $r$ leads to larger problem sizes in terms of both the number of variables and constraints in the CGLP that is solved at each node of the \BB tree, which leads to an extensive computational burden. For example, for most of the instances studied in Table~\ref{tab:greedy_Synthetic_1}, we ran into memory errors when using the CGLP of rank $r=2$.
						As an example, when there are 10 unfixed variables remaining at a layer of the \BB tree, to obtain the CGLP for rank $r=1$, the number of constraints in the linearized system will be multiplied by 20, whereas for CGLP of rank $r=2$, this size will be multiplied by 180.  
						Nonetheless, to provide some insight on the performance of the ML approach for ranks higher than one, we present in Table~\ref{tab:consistency_rank} computational experiments that compare the outcome of using ML for rank $r=1$ and $r=2$ for instances of Section~\ref{subsec:IP}.

						For these experiments, due to extensive computational burden of applying the consistency framework for higher ranks at each node of the \BB tree, as discussed above, we have used the common $(K,L)$ approach described as follows.
						The goal is to reduce the implementation time by (i) applying the CGLP (or its alternative ML approach in this case) in certain layers of the \BB tree only, and (ii) choosing a subset of unfixed variables to multiply with the constraints. We represent the layer candidates for applying the CGLP by $K$, and the variable indices used for multiplication by $L$. It is intuitive to pick top layers of the \BB tree to be included in $K$, and choose variable indices at the bottom layers
						of the tree to be included in $L$; see \cite{davarnia:ra:ho:2022} for further illustration of this approach.
						For the experiments presented in Tables~\ref{tab:consistency_rank} and \ref{tab:consistency_rank_benchmark}, we have chosen $K = 10$ and $L = 10$ for both ML approaches targeting partial LP-consistency of rank one and two to provide a fair comparison in a controlled setting.
						The columns in these tables are defined similarly to those of Tables~\ref{tab:greedy_Synthetic_1} and \ref{tab:benchmark} with a difference that the subcolumns under ``LR, $r=1$'' and ``LR, $r=2$'' include the \BB tree size and solution time for applying the logistic regression approach for partial LP-consistency of rank one and two using the $(K,L)$ approach.
						The implementation and algorithmic settings are similar to those used in Section~\ref{subsec:IP}.
						All instances in these experiments have been solved to optimality for both approaches.
						As observed in Tables~\ref{tab:consistency_rank} and \ref{tab:consistency_rank_benchmark}, the logistic regression method when used to approximate the outcome of the CGLP of rank $r=2$ outperforms that used for rank $r=1$ in both \BB tree size and total solution time in most of the instances.

						\begin{table}[!t]
							\caption{Comparison of the ML approach for different consistency ranks for multi-knapsack problems}
							\label{tab:consistency_rank}
							\resizebox{1.00\textwidth}{!}{
								\begin{tabular}{|l|l|l|ll|llll|llll|}
									\hline
									\multirow{2}{*}{$n$}  & \multirow{2}{*}{$m$}  & \multirow{2}{*}{\#} & \multicolumn{2}{c|}{\BB}                  & \multicolumn{4}{c|}{LR, $r=1$, and $K=L=10$}                                                                               & \multicolumn{4}{c|}{LR, $r=2$, and $K=L=10$}                                                            \\ \cline{4-13} 
									&                     &                           & \multicolumn{1}{l|}{nodes}   & time                & \multicolumn{1}{l|}{nodes}          & \multicolumn{1}{l|}{$\Delta$ (\%)}       & \multicolumn{1}{l|}{time}     & $\Delta$ (\%)                     & \multicolumn{1}{l|}{nodes}         & \multicolumn{1}{l|}{$\Delta$ (\%)}       & \multicolumn{1}{l|}{time}              & $\Delta$ (\%)                     \\ \hline
									\multirow{5}{*}{40} & \multirow{5}{*}{50} & 1                                & \multicolumn{1}{l|}{34473}   & 3142.60             & \multicolumn{1}{l|}{26651}  & \multicolumn{1}{l|}{\textit{23}} & \multicolumn{1}{l|}{2453.62}  & \textit{22}               & \multicolumn{1}{l|}{\textbf{22967}}          & \multicolumn{1}{l|}{\textit{34}} & \multicolumn{1}{l|}{\textbf{2186.80}}   & \textit{31}               \\ \cline{3-13} 
									&                     & 2                                 & \multicolumn{1}{l|}{82417}   & 7496.14             & \multicolumn{1}{l|}{76177}           & \multicolumn{1}{l|}{\textit{8}} & \multicolumn{1}{l|}{7051.42}  & \textit{6}               & \multicolumn{1}{l|}{\textbf{54347}} & \multicolumn{1}{l|}{\textit{35}} & \multicolumn{1}{l|}{\textbf{5116.49}}   & \textit{32}               \\ \cline{3-13} 
									&                     & 3                                & \multicolumn{1}{l|}{78965}   & 7167.15             & \multicolumn{1}{l|}{78717}  & \multicolumn{1}{l|}{\textit{0}} & \multicolumn{1}{l|}{7258.27}  & \textit{$-$1}               & \multicolumn{1}{l|}{\textbf{69523}}         & \multicolumn{1}{l|}{\textit{12}} & \multicolumn{1}{l|}{\textbf{6524.29}}  & \textit{9}               \\ \cline{3-13} 
									&                     & 4                                   & \multicolumn{1}{l|}{124327}  & 11377.96            & \multicolumn{1}{l|}{117233}  & \multicolumn{1}{l|}{\textit{6}} & \multicolumn{1}{l|}{10832.11}  & \textit{5}               & \multicolumn{1}{l|}{\textbf{101065}}         & \multicolumn{1}{l|}{\textit{19}} & \multicolumn{1}{l|}{\textbf{9421.17}}  & \textit{18}               \\ \cline{3-13} 
									&                     & 5                            & \multicolumn{1}{l|}{126209}  & 11700.49            & \multicolumn{1}{l|}{98121}  & \multicolumn{1}{l|}{\textit{24}} & \multicolumn{1}{l|}{9072.97}  & \textit{23}               & \multicolumn{1}{l|}{\textbf{61539}}          & \multicolumn{1}{l|}{\textit{52}} & \multicolumn{1}{l|}{\textbf{5772.71}}  & \textit{51}               \\ \hline
									\multirow{5}{*}{45} & \multirow{5}{*}{55} & 1                             & \multicolumn{1}{l|}{1050459} & 119705.93           & \multicolumn{1}{l|}{1035775} & \multicolumn{1}{l|}{\textit{2}} & \multicolumn{1}{l|}{116557.46} & \textit{3}               & \multicolumn{1}{l|}{\textbf{190993}}        & \multicolumn{1}{l|}{\textit{81}} & \multicolumn{1}{l|}{\textbf{17972.34}} & \textit{85}               \\ \cline{3-13} 
									&                     & 2                          & \multicolumn{1}{l|}{57621}   & 5449.79             & \multicolumn{1}{l|}{53507}  & \multicolumn{1}{l|}{\textit{8}} & \multicolumn{1}{l|}{4265.75}  & \textit{22}               & \multicolumn{1}{l|}{\textbf{32073}}          & \multicolumn{1}{l|}{\textit{45}} & \multicolumn{1}{l|}{\textbf{2971.28}}  & \textit{46}               \\ \cline{3-13} 
									&                     & 3                    & \multicolumn{1}{l|}{94281}        &  8424.51                   & \multicolumn{1}{l|}{62925}      & \multicolumn{1}{l|}{\textit{34}}   & \multicolumn{1}{l|}{5876.83}         & \textit{30}                 & \multicolumn{1}{l|}{\textbf{12029}}              & \multicolumn{1}{l|}{\textit{87}}   & \multicolumn{1}{l|}{\textbf{1149.06}}         & \textit{86}                 \\ \cline{3-13} 
									&                     & 4                     & \multicolumn{1}{l|}{450175}  & 44544.77            & \multicolumn{1}{l|}{398539} & \multicolumn{1}{l|}{\textit{12}} & \multicolumn{1}{l|}{39966.08} & \textit{11}               & \multicolumn{1}{l|}{\textbf{286923}}         & \multicolumn{1}{l|}{\textit{37}} & \multicolumn{1}{l|}{\textbf{27952.67}}  & \textit{38}               \\ \cline{3-13} 
									&                     & 5                         & \multicolumn{1}{l|}{630449}  & 65030.67            & \multicolumn{1}{l|}{630449} & \multicolumn{1}{l|}{\textit{0}} & \multicolumn{1}{l|}{65690.10} & \textit{$-$2}               & \multicolumn{1}{l|}{\textbf{263239}}        & \multicolumn{1}{l|}{\textit{58}} & \multicolumn{1}{l|}{\textbf{25554.13}} & \textit{61}               \\ \hline
									\multirow{5}{*}{45} & \multirow{5}{*}{60} & 1                     & \multicolumn{1}{l|}{989029}  & 116212.23           & \multicolumn{1}{l|}{951755} & \multicolumn{1}{l|}{\textit{4}} & \multicolumn{1}{l|}{104927.14} & \textit{10}               & \multicolumn{1}{l|}{\textbf{779599}}        & \multicolumn{1}{l|}{\textit{22}} & \multicolumn{1}{l|}{\textbf{85763.15}} & \textit{27}               \\ \cline{3-13} 
									&                     & 2                     & \multicolumn{1}{l|}{257217}  & 25121.01            & \multicolumn{1}{l|}{255045} & \multicolumn{1}{l|}{\textit{1}} & \multicolumn{1}{l|}{23902.25}  & \textit{5}               & \multicolumn{1}{l|}{\textbf{207051}}         & \multicolumn{1}{l|}{\textit{20}} & \multicolumn{1}{l|}{\textbf{19344.72}}  & \textit{24}               \\ \cline{3-13} 
									&                     & 3                        & \multicolumn{1}{l|}{249869}  & 23936.52            & \multicolumn{1}{l|}{233177} & \multicolumn{1}{l|}{\textit{7}} & \multicolumn{1}{l|}{21761.88}  & \textit{10}               & \multicolumn{1}{l|}{\textbf{63619}}         & \multicolumn{1}{l|}{\textit{74}} & \multicolumn{1}{l|}{\textbf{5980.96}}  & \textit{76}               \\ \cline{3-13} 
									&                     & 4                        &  \multicolumn{1}{l|}{-}        &     \textgreater{}86400                & \multicolumn{1}{l|}{823355}      & \multicolumn{1}{l|}{\textit{-}}   & \multicolumn{1}{l|}{90986.38}         & \multicolumn{1}{l|}{\textit{-}}                 & \multicolumn{1}{l|}{\textbf{641949}}              & \multicolumn{1}{l|}{\textit{-}}   & \multicolumn{1}{l|}{\textbf{77075.34}}         & \multicolumn{1}{l|}{\textit{-}}                 \\ \cline{3-13} 
									&                     & 5                       & \multicolumn{1}{l|}{320967}  & 31300.88            & \multicolumn{1}{l|}{320363} & \multicolumn{1}{l|}{\textit{0}} & \multicolumn{1}{l|}{31782.85} & \textit{$-$2}               & \multicolumn{1}{l|}{\textbf{182909}}         & \multicolumn{1}{l|}{\textit{45}} & \multicolumn{1}{l|}{\textbf{16944.11}}  & \textit{46}               \\ \hline
									\multirow{5}{*}{50} & \multirow{5}{*}{60} & 1                          & \multicolumn{1}{l|}{329277}  & 31947.91            & \multicolumn{1}{l|}{270375} & \multicolumn{1}{l|}{\textit{18}} & \multicolumn{1}{l|}{25977.05} & \textit{20}               & \multicolumn{1}{l|}{\textbf{166139}}         & \multicolumn{1}{l|}{\textit{50}} & \multicolumn{1}{l|}{\textbf{15501.51}}  & \textit{52}               \\ \cline{3-13} 
									&                     & 2                          & \multicolumn{1}{l|}{531833}  & 54939.26            & \multicolumn{1}{l|}{516531} & \multicolumn{1}{l|}{\textit{3}} & \multicolumn{1}{l|}{52111.54} & \textit{6}               & \multicolumn{1}{l|}{\textbf{93973}}         & \multicolumn{1}{l|}{\textit{82}} & \multicolumn{1}{l|}{\textbf{8679.84}}  & \textit{85}               \\ \cline{3-13} 
									&                     & 3                      & \multicolumn{1}{l|}{-}       & \textgreater{}86400 & \multicolumn{1}{l|}{-} & \multicolumn{1}{l|}{\textit{-}}  & \multicolumn{1}{l|}{-} & \multicolumn{1}{l|}{-} & \multicolumn{1}{l|}{-}        & \multicolumn{1}{l|}{\textit{-}}  & \multicolumn{1}{l|}{-} & \multicolumn{1}{l|}{-} \\ \cline{3-13} 
									&                     & 4                           & \multicolumn{1}{l|}{174357}  & 16540.17            & \multicolumn{1}{l|}{171407}  & \multicolumn{1}{l|}{\textit{2}} & \multicolumn{1}{l|}{15523.94}  & \textit{7}               & \multicolumn{1}{l|}{\textbf{151981}}         & \multicolumn{1}{l|}{\textit{13}} & \multicolumn{1}{l|}{\textbf{13994.98}}  & \textit{16}               \\ \cline{3-13} 
									&                     & 5                         & \multicolumn{1}{l|}{-}       & \textgreater{}86400 & \multicolumn{1}{l|}{-} & \multicolumn{1}{l|}{\textit{-}}  & \multicolumn{1}{l|}{-} & \multicolumn{1}{l|}{-} & \multicolumn{1}{l|}{-}        & \multicolumn{1}{l|}{\textit{-}}  & \multicolumn{1}{l|}{-} & \multicolumn{1}{l|}{-} \\ \hline
								\end{tabular}
							}
						\end{table}

						\begin{table}[!t]
							\caption{Comparison of the ML approach for different consistency ranks for MIPLIB problems}
							\label{tab:consistency_rank_benchmark}
							\resizebox{1.00\textwidth}{!}{
								\begin{tabular}{|l|l|l|ll|llll|llll|}
									\hline
									\multirow{2}{*}{Class}  & \multirow{2}{*}{$n$}  & \multirow{2}{*}{$m$} & \multicolumn{2}{c|}{\BB}                  & \multicolumn{4}{c|}{LR, $r=1$, and $K=L=10$}                                                                               & \multicolumn{4}{c|}{LR, $r=2$, and $K=L=10$}                                                            \\ \cline{4-13} 
									&                     &                           & \multicolumn{1}{l|}{nodes}   & time                & \multicolumn{1}{l|}{nodes}          & \multicolumn{1}{l|}{$\Delta$ (\%)}       & \multicolumn{1}{l|}{time}     & $\Delta$ (\%)                     & \multicolumn{1}{l|}{nodes}         & \multicolumn{1}{l|}{$\Delta$ (\%)}       & \multicolumn{1}{l|}{time}              & $\Delta$ (\%)                     \\ \hline
									$\mathtt{p0033}$ & 33 & 15                                & \multicolumn{1}{l|}{32117}   & 2961.55             & \multicolumn{1}{l|}{22161}  & \multicolumn{1}{l|}{\textit{31}} & \multicolumn{1}{l|}{2191.14}  & \textit{26}               & \multicolumn{1}{l|}{\textbf{13811}}          & \multicolumn{1}{l|}{\textit{57}} & \multicolumn{1}{l|}{\textbf{1658.16}}   & \textit{44}               \\ \hline 
									$\mathtt{pipex}$ & 41 & 48                                & \multicolumn{1}{l|}{1057}   & 1933.48             & \multicolumn{1}{l|}{983}  & \multicolumn{1}{l|}{\textit{7}} & \multicolumn{1}{l|}{1875.01}  & \textit{3}               & \multicolumn{1}{l|}{\textbf{865}}          & \multicolumn{1}{l|}{\textit{18}} & \multicolumn{1}{l|}{\textbf{1643.5}}   & \textit{15}               \\ \hline 
									$\mathtt{sentoy}$ & 60 & 30                                & \multicolumn{1}{l|}{703}   & 391.12             & \multicolumn{1}{l|}{611}  & \multicolumn{1}{l|}{\textit{13}} & \multicolumn{1}{l|}{336.16}  & \textit{14}               & \multicolumn{1}{l|}{\textbf{505}}          & \multicolumn{1}{l|}{\textit{28}} & \multicolumn{1}{l|}{\textbf{273.7}}   & \textit{30}               \\ \hline 
									$\mathtt{stein27}$ & 27 & 118                                & \multicolumn{1}{l|}{6099}   & 2474.1             & \multicolumn{1}{l|}{5001}  & \multicolumn{1}{l|}{\textit{18}} & \multicolumn{1}{l|}{\textbf{2152.37}}  & \textit{13}               & \multicolumn{1}{l|}{\textbf{4879}}          & \multicolumn{1}{l|}{\textit{20}} & \multicolumn{1}{l|}{2251.34}   & \textit{9}               \\ \hline 
									$\mathtt{enigma}$ & 100 & 42                                & \multicolumn{1}{l|}{98545}   & 7413.51             & \multicolumn{1}{l|}{70951}  & \multicolumn{1}{l|}{\textit{28}} & \multicolumn{1}{l|}{5782.14}  & \textit{22}               & \multicolumn{1}{l|}{\textbf{57155}}          & \multicolumn{1}{l|}{\textit{42}} & \multicolumn{1}{l|}{\textbf{4077.15}}   & \textit{45}               \\ \hline 
									$\mathtt{lseu}$ & 89 & 28                                & \multicolumn{1}{l|}{219597}   & 4138.8             & \multicolumn{1}{l|}{217401}  & \multicolumn{1}{l|}{\textit{1}} & \multicolumn{1}{l|}{4183.11}  & \textit{$-$1}               & \multicolumn{1}{l|}{\textbf{199833}}          & \multicolumn{1}{l|}{\textit{9}} & \multicolumn{1}{l|}{\textbf{3641.44}}   & \textit{12}               \\ \hline 
								\end{tabular}
							}
						\end{table}

					\end{APPENDICES}


\end{document}